\documentclass[12pt]{amsart}
\usepackage{xypic}

\usepackage{amsmath, amsthm, amssymb}

\usepackage[colorlinks=true]{hyperref}
\usepackage{mathrsfs}  

\theoremstyle{plain} 
\newtheorem{thm}[equation] {Theorem}
\newtheorem{cor}[equation]{Corollary}
\newtheorem{lem}[equation]{Lemma}
\newtheorem{prop}[equation]{Proposition}

\theoremstyle{definition}
\newtheorem{defn}[equation]{Definition}
\newtheorem{ex}[equation]{Example}

\theoremstyle{remark}
\newtheorem{rem}[equation]{Remark}

\numberwithin{equation}{section}
\numberwithin{table}{section}

\DeclareMathOperator{\Hom}{Hom}
\DeclareMathOperator{\im}{im}
\DeclareMathOperator{\orb}{orb}
\DeclareMathOperator{\rk}{rk}
\DeclareMathOperator{\soc}{soc}
\DeclareMathOperator{\ann}{ann}
\DeclareMathOperator{\spn}{span}
\newcommand{\ra}{\rightarrow}
\newcommand{\field}[1]{\mathbb{#1}}
\newcommand{\charac}[1]{\widehat{#1}}
\newcommand{\addchar}[1]{\charac{#1}}
\newcommand{\F}{\ensuremath{\field{F}}}
\newcommand{\C}{\ensuremath{\field{C}}}
\newcommand{\N}{\ensuremath{\field{N}}}
\newcommand{\Q}{\ensuremath{\field{Q}}}
\newcommand{\R}{\ensuremath{\field{R}}}
\newcommand{\Z}{\ensuremath{\field{Z}}}
\newcommand{\hr}{\ensuremath{\addchar{R}}}
\newcommand{\size}[1]{\lvert #1 \rvert}
\newcommand{\abs}[1]{\lvert #1 \rvert}
\newcommand{\al}{\alpha}
\newcommand{\be}{\beta}
\newcommand{\ga}{\gamma}
\newcommand{\la}{\lambda}
\newcommand{\La}{\Lambda}
\newcommand{\m}{\ensuremath{\mathfrak m}}
\newcommand{\OO}{\ensuremath{\mathcal O}}
\newcommand{\OS}{\ensuremath{{\mathcal O}^\sharp}}
\newcommand{\U}{\ensuremath{\mathcal{U}}}
\newcommand{\qbinom}[2]{\begin{bmatrix}{#1}\\{#2}\end{bmatrix}_q}
\newcommand{\smqbinom}[2]{\left[\begin{smallmatrix}{#1}\\{#2}\end{smallmatrix}\right]_q}

\DeclareMathOperator{\tr}{tr}
\DeclareMathOperator{\Tr}{Tr}
\newcommand{\A}{\ensuremath{\mathcal{A}}}
\newcommand{\ah}{\ensuremath{\charac{A}}}
\newcommand{\PP}{\ensuremath{\mathcal{P}}}

\DeclareMathOperator{\lt}{lt}
\DeclareMathOperator{\rt}{rt}
\newcommand{\glt}{G_{\lt}}
\newcommand{\grt}{G_{\rt}}

\newcommand{\ft}[1]{\widehat{#1}}
\DeclareMathOperator{\hwe}{hwe}
\DeclareMathOperator{\howe}{howe}

\DeclareMathOperator{\wwe}{wwe}
\DeclareMathOperator{\stab}{stab}
\DeclareMathOperator{\GL}{GL}

\newcommand{\wh}{{\textsc{h}}}
\newcommand{\wg}{{\normalfont\textsc{w}}}
\newcommand{\wrk}{{\textsc{r}}}
\DeclareMathOperator{\col}{col}
\DeclareMathOperator{\row}{row}
\DeclareMathOperator{\length}{length}
\DeclareMathOperator{\efflength}{efflng}
\newcommand{\asing}{A^{\rm sing}}

\DeclareMathOperator{\pe}{pe}
\DeclareMathOperator{\ce}{ce}
\DeclareMathOperator{\se}{se}
\DeclareMathOperator{\spec}{spec}
\DeclareMathOperator{\re}{rpe}

\DeclareMathOperator{\rowsp}{rowsp}
\DeclareMathOperator{\colsp}{colsp}
\newcommand{\sW}{\overline{W}}
\newcommand{\eb}{\overline{\eta}}
\newcommand{\ob}{\overline{\omega}}
\DeclareMathOperator{\sAnn}{\overline{Ann}}
\newcommand{\os}{\bar{\sigma}}

\title[Weights with maximal symmetry
]{Weights with Maximal Symmetry and \\
Failures of the MacWilliams Identities 
} 

\author[J. A. Wood]{Jay A. Wood}
\address{
Western Michigan University
}
\email{jay.wood@wmich.edu}
\dedicatory{In memoriam: \\
Melvin G. Rothenberg, \oldstylenums{1934}--\oldstylenums{2023} \\
Erna B. Yackel, \oldstylenums{1939}--\oldstylenums{2022} \\
Robert J. Zimmer, \oldstylenums{1947}--\oldstylenums{2023} }

\begin{document}



\begin{abstract}
This paper examines the $w$-weight enumerators of weights $w$ with maximal symmetry over finite chain rings and matrix rings over finite fields.  
In many cases, including the homogeneous weight, the MacWilliams identities for $w$-weight enumerators fail because there exist two linear codes with the same $w$-weight enumerator whose dual codes have different $w$-weight enumerators.
\end{abstract}

\maketitle
\thispagestyle{empty}

\section{Introduction}

The MacWilliams identities \cite{MR0149978} reveal a relationship between the Hamming weight enumerator ($\hwe$) of a linear code $C$ over a finite field $\F_q$ and the Hamming weight enumerator of its dual code $C^\perp$:
\[ \hwe_{C^\perp}(X,Y) = \frac{1}{\size{C}} \hwe_C(X+(q-1)Y, X-Y) . \]

One way to try to generalize this result is to use any integer-valued weight $w$ on a finite ring $R$ with $1$.   The homogeneous weight, suitably normalized, is an example.  Assume $w(0)=0$ and $w(r) >0$ for $r \neq 0$.  Denote the maximum value of $w$ by $w_{\max}$.  For a vector $x = (x_1, x_2, \ldots, x_n) \in R^n$, set $w(x) = \sum_{i=1}^n w(x_i)$.  Define the $w$-weight enumerator ($\wwe$) of a left $R$-linear code $C \subseteq R^n$ by
\[ \wwe_C(X,Y) = \sum_{c \in C} X^{n w_{\max} - w(c)} Y^{w(c)} . \]
Do MacWilliams identities hold for $w$-weight enumerators?

This paper will show that MacWilliams identities seldom hold for weights having maximal symmetry when the ring is a finite chain ring or a matrix ring over a finite field.  
For example, suppose $R = \Z/4\Z$ and $w$ is an integer-valued weight on $R$ with maximal symmetry (i.e., $w(1)=w(3)$ for this ring).  Then the only weights for which the MacWilliams identities hold are multiples of the Hamming weight and multiples of the homogeneous weight (the Lee weight for this ring), Corollary~\ref{cor:chainm2}.  Similarly, let $R = M_{2 \times 2}(\F_q)$ and $w$ be an integer-valued weight with maximal symmetry (i.e., the value of $w(r)$ depends only on the rank of $r$).  Then the only weights for which the MacWilliams identities hold are multiples of the Hamming weight (any $q$) and multiples of the homogeneous weight ($q=2$ only), Theorem~\ref{main-theorem-for-2by2}.  More generally, the homogeneous weight on $R = M_{k \times k}(\F_q)$, $k \geq 2$, satisfies the MacWilliams identities if and only if $k=q=2$, Theorem~\ref{thm:homog-failures}.

The reason the MacWilliams identities often fail is that there exist $R$-linear codes $C, D \subseteq R^n$ for some $n$ such that $\wwe_C = \wwe_D$ but $\wwe_{C^\perp} \neq \wwe_{D^\perp}$.
This claim presents two challenges:
\begin{itemize}
\item to construct linear codes $C,D \subseteq R^n$ with $\wwe_C = \wwe_D$ in such a way that it is then possible \ldots
\item  to detect differences $\wwe_{C^\perp} \neq \wwe_{D^\perp}$ in the $w$-weight enumerators of their dual codes.
\end{itemize}

A weight $w$ on a finite ring $R$ with $1$ has maximal symmetry when $w(u r u') = w(r)$ for all $r \in R$ and units $u, u'$ in $R$.  Any left linear code $C \subseteq R^n$ can be viewed as the image of an injective homomorphism $\La: M \ra R^n$ of left $R$-modules, for some finite left $R$-module $M$.  The group of units of $R$ acts on $M$ on the left, and the maximal symmetry hypothesis implies that $x \mapsto w(x \La)$ is constant on each orbit of this group action.  Writing $[x]$ for the orbit of $x \in M$, we see that 
\begin{equation}  \label{eqn:wwe-by-orbits}
\wwe_C(X,Y) = \sum_{\text{orbits $[x]$}} \size{[x]} X^{n w_{\max} - w(x\La)} Y^{w(x\La)} .
\end{equation}
In essence, the choice of homomorphism $\La$ determines the `orbit weights' $\omega(x) = w(x \La)$ assigned to each orbit $[x]$ in $M$.  

In order to find another linear code $D \subseteq R^n$ with $\wwe_C = \wwe_D$, one can try to permute the weights assigned to orbits while keeping the sum in \eqref{eqn:wwe-by-orbits} unchanged.  For example, in $R = M_{2 \times 2}(\F_2)$ there is one orbit of size $6$ consisting of matrices of rank $2$ and three orbits, each of size $3$, of matrices of rank $1$.  One can then try to construct linear codes whose orbit weights behave as follows:
\[  \begin{array}{c|ccc|c}
\text{orbit} & [x_1] & [x_2] & [x_3] & [y] \\ \hline
\text{size of orbit} & 3 & 3 & 3 & 6 \\
w(x\La) & a & b & b & c \\
w(x \La') & a & c & c & b 
\end{array}  \]
It is not obvious a priori that such constructions are possible, but Section~\ref{sec:construction} shows that constructions of this type can be carried out for all matrix rings $M_{k \times k}(\F_q)$ over finite fields.  

When $R$ is a finite chain ring, all the orbits have different sizes, so the permutation idea does not work.  However, one can use different modules as the domains of the defining homomorphisms.  
For example, the ring $\Z/8Z$ has three modules of size $8$: $M=\Z/8Z$ itself, $\Z/4Z \oplus \Z/2Z$, and 
$M'=\Z/2Z\oplus \Z/2Z\oplus \Z/2Z$.  Considering $M$ and $M'$, the orbits of $M$ are $\{1,3,5,7\}$, $\{2,6\}$, $\{4\}$, and $\{0\}$, while the orbits of $M'$ are the $8$ subsets of size $1$.  By choosing certain unions of orbits of $M'$, say $\{ 100, 101, 110, 111 \}$, $\{ 011, 010\}$, $\{ 001\}$, $\{000 \}$, of the same size as the orbits of $M$, one can try to construct homomorphisms $\La: M \ra R^n$ and $\La': M' \ra R^n$ achieving the same weights on corresponding orbits.  While this may not seem possible at first glance, Section~\ref{sec:twoFamiliesChain} details how such constructions exist.

In order to show that $\wwe_{C^\perp} \neq \wwe_{D^\perp}$, it is enough to show that $A_j(C^\perp) \neq A_j(D^\perp)$ for some $j>0$; here, $A_j(C^\perp)$ is the number of codewords $v \in C^\perp$ with $w(v)=j$.  The easiest case to understand is when $v \in C^\perp$ has exactly one nonzero entry; such a $v$ is called a \emph{singleton}.    When $C \subseteq R^n$ is the image of a homomorphism $\La: M \ra R^n$, the components of $\La=(\la_1, \la_2, \ldots, \la_n)$ are elements $\la_i \in \Hom_R(M,R)$.  By understanding how many elements $r \in R$ annihilate any given $\la_i$, i.e., $\la_i r =0$, one can write down formulas for the contributions of singletons to $A_j(C^\perp)$, Proposition~\ref{prop:SingletonContribution}.  When $j$ is sufficiently small, only singletons can contribute to $A_j(C^\perp)$, Corollary~\ref{cor:only-singletons}.  This technique turns out to be surprisingly effective in allowing one to prove that $A_j(C^\perp) \neq A_j(D^\perp)$ in a large number of situations.

This paper is divided into three parts.  The first establishes notation and ideas that can apply to any finite ring with $1$.  In particular, codes will usually be linear codes over a finite ring with $1$.  A weight $w$ on $R$ will be assumed to have maximal symmetry and have positive integer values, and $w$ will be extended additively to $R^n$.  The second part examines the construction of linear codes and the analysis of singleton dual codewords over a finite chain ring, while the third part does the same for matrix rings over finite fields.  An appendix provides a short outline of a proof of the MacWilliams identities over finite Frobenius rings using the Fourier transform and the Poisson summation formula.

\part{Generalities}

\section{Preliminaries}
This section will review without proof some terminology and results from \cite{wood:turkey} about characters of finite abelian groups, finite Frobenius rings, additive and linear codes and their dual codes, symmetry groups, and weights.

Let $A$ be a finite abelian group.  A \emph{character} of $A$ is a group homomorphism $\pi: A \ra (\C^\times, \cdot)$ from $A$ to the multiplicative group of nonzero complex numbers.  Denote by $\ah$ the set of all characters of $A$; $\ah$ is a finite abelian group under pointwise multiplication of functions.   The groups $A$ and $\ah$ are isomorphic, but not naturally so; $\size{\ah}=\size{A}$.  The double character group is naturally isomorphic to the group: $(\ah)\charac{\phantom{A}} \cong A$; $a \in A$ corresponds to evaluation at $a$ of $\pi \in \ah$, i.e., $\pi \mapsto \pi(a)$.

For a subgroup $B \subseteq A$, define its \emph{annihilator} by
\[  (\ah:B) = \{ \pi \in \ah: \pi(b)=1, \text{for all $b \in B$}\} . \]
Then $(\ah:B)$ is a subgroup of $\ah$, $(\ah:B) \cong (A/B)\charac{\phantom{A}}$, and $\size{(\ah:B)} = \size{A/B} = \size{A}/\size{B}$.   
Identifying $A \cong (\ah)\charac{\phantom{A}}$, we have $(A:(\ah:B)) = B$. 

Throughout this paper $R$ will denote a finite (associative) ring with $1$; $R$ may be noncommutative.  The group of units (invertible elements) of $R$ is denoted $\U=\U(R)$.  The \emph{Jacobson radical} $J(R)$ of $R$ is the intersection of all maximal left ideals of $R$; $J(R)$ is itself a two-sided ideal of $R$.  The left/right \emph{socle} $\soc({}_R R)$, $\soc(R_R)$ of $R$ is the left/right ideal generated by the minimal left/right ideals of $R$.  A ring $R$ (perhaps infinite) is \emph{Frobenius} if ${}_R J(R) \cong \soc({}_R R)$ and $J(R)_R \cong \soc(R_R)$ \cite[Theorem~(16.14)]{lam:modules}; a theorem of Honold \cite{honold:frobenius} says that one of these isomorphisms suffices for finite rings.

Every finite ring $R$ has an underlying additive abelian group.  Its character group $\hr$ is a bimodule over $R$.  The two scalar multiplications are written in exponential form, with $\pi \in \hr$, $r,s \in R$:
\[  ({}^r\pi)(s) = \pi(sr), \quad \pi^r(s) = \pi(rs) . \]
A finite ring $R$ is Frobenius if and only if $R \cong \hr$ as left (resp., right) $R$-modules, \cite[Theorem~3.10]{wood:duality}.  This implies that a finite Frobenius ring admits a character $\chi$, called a \emph{generating character}, such that $r \mapsto {}^r \chi$ is an isomorphism of left $R$-modules (resp., $r \mapsto \chi^r$ is an isomorphism of right $R$-modules).   A generating character has the property that any one-sided ideal of $R$ that is contained in $\ker \chi$ must be the zero ideal.

An \emph{additive code} of length $n$ over $R$ is an additive subgroup $C \subseteq R^n$.  If $C \subseteq R^n$ is a left, resp., right, $R$-submodule, then $C$ is a left (resp., right) \emph{$R$-linear code}.  One way to present a left $R$-linear code is as the image $C = \im \Lambda$ of a homomorphism $\Lambda: M \ra R^n$ of left $R$-modules, and similarly for right linear codes.

We will write homomorphisms of left $R$-modules
with inputs on the left, so that preservation of scalar multiplication is $(rx)\phi = r(x \phi)$, where $r \in R$, $x \in M$, $M$ a left $R$-module, and $\phi$ a homomorphism of left $R$-modules with domain $M$

Define the \emph{standard dot product} on $R^n$ by
\[  x \cdot y = \sum_{i=1}^n x_i y_i \in R, \]
for $x= (x_1, x_2, \ldots, x_n), y=(y_1, y_2, \ldots, y_n) \in R^n$. 
Given an additive code $C \subseteq R^n$, define dual codes by
\begin{align}  
\mathcal{L}(C) &= \{ y \in R^n: y \cdot x = 0, \text{for all $x \in C$} \},  \notag \\
\label{eq:RightDualCode}
\mathcal{R}(C) &= \{ y \in R^n: x \cdot y = 0, \text{for all $x \in C$} \} .
\end{align}
When $R$ is Frobenius, with generating character $\chi$, also define
\begin{align}  
\mathfrak{L}(C) &= \{ y \in R^n: \chi(y \cdot x) = 0, \text{for all $x \in C$}  \}, \notag \\
\label{eq:FrakDualCode}
\mathfrak{R}(C) &= \{ y \in R^n:  \chi(x \cdot y) = 0, \text{for all $x \in C$} \} .
\end{align}

Using the isomorphisms $r \mapsto {}^r \chi$ and $r \mapsto \chi^r$ of $R$ to  $\hr$, 
there are isomorphisms $R^n \ra \hr^n$ of left, resp., right, $R$-modules given by $x \mapsto {}^x\chi$ and $x \mapsto \chi^x$, where ${}^x\chi(y) = \chi(y \cdot x)$ and $\chi^x(y) = \chi(x \cdot y)$, for $x, y \in R^n$.  Under the isomorphism $x \mapsto {}^x \chi$, $\mathfrak{R}(C)$ is taken to $(\hr^n:C)$, while under the isomorphism $x \mapsto \chi^x$, $\mathfrak{L}(C)$ is taken to $(\hr^n:C)$.

\begin{lem}  \label{lem:DualityFacts}
Suppose $R$ is Frobenius and $C \subseteq R^n$ is an additive code.  Then
\begin{itemize}
\item $\size{C} \cdot \size{\mathfrak{L}(C)} = \size{C} \cdot \size{\mathfrak{R}(C)} = \size{R^n}$;
\item  $\mathfrak{L}(\mathfrak{R}(C)) = C = \mathfrak{R}(\mathfrak{L}(C))$.
\end{itemize}
\end{lem}

\begin{rem}  \label{rem:LinearCaseAnnihilators}
Note that $\mathcal{R}(C) \subseteq \mathfrak{R}(C)$ and $\mathcal{L}(C) \subseteq \mathfrak{L}(C)$.  In general these containments will be proper.  However, if $C$ is a left $R$-linear code, then $\mathcal{R}(C) = \mathfrak{R}(C)$.  Similarly, $\mathcal{L}(C) = \mathfrak{L}(C)$ if $C$ is right $R$-linear.
\end{rem}

A \emph{weight} on $R$ is a function $w: R \ra \C$ from $R$ to the complex numbers $\C$ with $w(0)=0$.  In most of this paper we will study weights having positive integer values, except for $w(0)=0$.  A weight $w$ will be extended additively to $R^n$, so that $w(v) = \sum_{i=1}^n w(v_i) \in \C$, where $v = (v_1, v_2,\ldots, v_n) \in R^n$.

Every weight $w$ on $R$ has two \emph{symmetry groups}, left and right:
\begin{align}
G_{\lt}(w) &= \{ u \in \U: w(ur) = w(r) \text{ for all $r \in R$}\}, \notag \\
G_{\rt}(w) &= \{ u \in \U: w(ru) = w(r) \text{ for all $r \in R$}\}. \label{eq:SymmetryGroups}
\end{align}
A weight $w$ has \emph{maximal symmetry} when $G_{\lt}(w) = G_{\rt}(w) = \U$.

\begin{ex} \label{ex:Hamming}
The most well-known weight on $R$ is the \emph{Hamming weight} $\wh$, defined by $\wh(0)=0$ and $\wh(r) = 1$ for $r \neq 0$.  The Hamming weight has maximal symmetry.
\end{ex}

\begin{ex}  \label{ex:homogeneous}
Another well-known weight on $R$ having maximal symmetry is the \emph{homogeneous weight} $\wg: R \ra \R$.  The homogeneous weight was first introduced in \cite{Constantinescu-Heise:metric} over integer residue rings and generalized to all finite rings and modules in \cite{Honold-Nechaev:weighted-modules} and \cite{greferath-schmidt:combinatorics}.  A homogeneous weight is characterized by the choice of a real number $\zeta >0$ and the following properties \cite{greferath-schmidt:combinatorics}:
\begin{itemize}
\item  $\wg(0)=0$;
\item $G_{\lt}(\wg) = \U$; and
\item  $\sum_{x \in Rr} \wg(x) = \zeta \size{Rr}$ for nonzero principal left ideals $Rr \subseteq R$.
\end{itemize}
The last property says that all nonzero left principal ideals of $R$ have the same average weight $\zeta$.  In fact, the average weight property holds for all nonzero left ideals of $R$ if and only if $R$ is Frobenius \cite[Corollary~1.6]{greferath-schmidt:combinatorics}.

When $R = \F_q$, all the nonzero elements are units, so $\wg(u) = \wg(1)$ for all units $u$.  Thus $\wg$ is a constant multiple (namely, $\wg(1)$) times the Hamming weight.  Note that $\zeta =(q-1) \wg(1)/q$ over $\F_q$.

Greferath and Schmidt \cite[Theorem~1.3]{greferath-schmidt:combinatorics} prove that homogeneous weights exist on any $R$ by giving an explicit formula for $\wg$ in terms of $\zeta$ and the M\"obius function $\mu$  (see \cite[\S 5.5]{MR4249619}) of the poset of principal left ideals of $R$; namely:
\begin{equation}  \label{eq:GS-homogwt}
\wg(r) = \zeta \left( 1 - \frac{\mu(0, Rr)}{\size{\U r}} \right) , \quad r \in R.
\end{equation}
This formula implies that all the values of $\wg$ are rational multiples of $\zeta$.  By choosing $\zeta$ appropriately, one can produce a homogeneous weight on $R$ with integer values.  Another consequence of the formula is that any two homogeneous weights on $R$ are scalar multiples of each other:  if $\wg$ and $\wg'$ are homogeneous weights on $R$ with average weights $\zeta$ and $\zeta'$, respectively, then $\wg' = (\zeta'/\zeta) \wg$.
\end{ex}

\begin{ex}
Let $R = M_{k \times k}(\F_q)$ be the ring of $k \times k$ matrices over the finite field $\F_q$.  The \emph{rank weight} $\wrk$ is defined by $\wrk(r) = \rk(r)$, the usual rank of the matrix $r \in R$.  The rank weight has maximal symmetry.
\end{ex}

\section{MacWilliams identities}  \label{sec:MWIdentities}
In her 1962 doctoral dissertation, Florence Jessie MacWilliams gave a formula relating the Hamming weight enumerator of a linear code over a finite field to the Hamming weight enumerator of its dual code \cite{MR2939359,MR0149978}.  In this section we will describe this work of MacWilliams as well as some of its generalizations.

For any linear code $C \subseteq R^n$ over a finite ring $R$, the \emph{Hamming weight enumerator} is the following homogeneous polynomial of degree $n$:
\begin{equation}  \label{eq:DefnHammingWE}
\hwe_C(X,Y) = \sum_{x \in C} X^{n-\wh(x)} Y^{\wh(x)} ,
\end{equation}
where $\wh$ is the Hamming weight, as in Example~\ref{ex:Hamming}.
The formula relating $\hwe_C$ and $\hwe_{C^\perp}$ is quoted next.

\begin{thm}[MacWilliams identities {\cite{MR2939359,MR0149978}}] \label{thm:MWids-original}
If $C \subseteq \F_q^n$ is a linear code over the finite field $\F_q$, then
\[  \hwe_{C^\perp}(X,Y) = \frac{1}{\size{C}} \hwe_C(X+(q-1)Y, X-Y) . \]
\end{thm}

\begin{rem}
Note in particular that the formula for $\hwe_{C^\perp}$ depends only on $\hwe_C$ and not on a more detailed knowledge of the code $C$.  
By applying the MacWilliams identities to $C^\perp$ and $C=(C^\perp)^\perp$, the roles of $C$ and $C^\perp$ can be reversed. 
\end{rem}

We isolate one consequence of the MacWilliams identities.
\begin{cor}  \label{cor:respect-duality}
 If $C$ and $D$ are two linear codes over $\F_q$ with $\hwe_C = \hwe_D$, then $\hwe_{C^\perp} = \hwe_{D^\perp}$.
\end{cor}

The MacWilliams identities 
for the Hamming weight enumerator
can be generalized in several ways.  One way is to generalize the algebraic structure of the codes.  There are versions of the MacWilliams identities with the Hamming weight enumerator for additive codes over finite abelian groups \cite{delsarte:linear-programming}, as well as for left (or right) linear codes over a finite Frobenius ring \cite[Theorem~8.3]{wood:duality}.   In the latter, one replaces $q$ with $\size{R}$ and $C^\perp$ with $\mathcal{R}(C)$ (with $\mathcal{L}(C)$ if $C$ is right linear).   

Another way to generalize the MacWilliams identities is to generalize the enumerator.  There are two  broad ways of doing this, stemming from two interpretations of the exponents in \eqref{eq:DefnHammingWE}.  Following Gluesing-Luerssen \cite{MR3336966}, one of the generalizations will be called \emph{partition enumerators}; the other will be called \emph{$w$-weight enumerators}.  These enumerators will be defined below, and the Hamming weight enumerator will be an example of both.  
While most of the following material can be formulated for additive codes over finite abelian groups, the discussion here will be restricted to linear codes over finite rings.

Suppose $S$ is a finite set.  A \emph{partition} of $S$ is a collection $\mathcal{P} = \{ P_i \}$ of nonempty subsets of $S$ such that the subsets are pairwise disjoint and cover $S$, i.e., $S = \uplus_i P_i$.  The subsets $P_i$ are called the \emph{blocks} of the partition.

Suppose a finite ring $R$ has a partition $\mathcal{P} = \{ P_i \}_{i=1}^m$.  Define counting functions
$n_i : R^n \ra \N$, $i=1, 2, \ldots, m$, by $n_i(x) = \size{\{j: x_j \in P_i\}}$, for $x = (x_1, x_2, \ldots, x_n) \in R^n$.  The counting functions count how many entries of $x$ belong to each block of the partition.  For a linear code $C \subseteq R^n$ define the \emph{partition enumerator} associated to $C$ and the partition $\mathcal{P}$ to be the following homogeneous polynomial of degree $n$ in the variables $Z_1, Z_2, \ldots, Z_m$:
\begin{equation}
\pe_C^{\mathcal{P}}(Z_1, \ldots, Z_m) = \sum_{x \in C} \prod_{i=1}^m Z_i^{n_i(x)} .
\end{equation}

Examples of such partition enumerators include:
\begin{itemize}
\item the \emph{complete enumerator} ($\ce$) based on the singleton partition $\mathcal{P} = \{ \{r\} \}_{r \in R}$;
\item a \emph{symmetrized enumerator} ($\se$) based on a partition consisting of the orbits of a group action on $R$;
\item  the Hamming (weight) enumerator based on the partition with blocks $\{0\}$ and the set difference $R-\{0\}$.
\end{itemize}
While the literature refers to the examples above as weight enumerators, the first two do not involve weights, so I will use the shorter names indicated.  

Suppose $R$ has two partitions $\mathcal{P} = \{ P_i \}_{i=1}^m$ and $\mathcal{Q} = \{ Q_j\}_{j=1}^{m'}$.  Also suppose $\mathcal{P}$ is a \emph{refinement} of $\mathcal{Q}$, i.e., each block $P_i$ is contained in some (unique) block $Q_j$; write $j = f(i)$.  Write the partition enumerators of a linear code $C \subseteq R^n$, using variables $Z_i$, $i=1,2,\ldots,m$,  for $\mathcal{P}$, and $\mathcal{Z}_j$, $j=1,2,\ldots, m'$, for $\mathcal{Q}$:
\[  \pe_C^{\mathcal{P}}(Z_1, \ldots, Z_m) \quad \text{and} \quad \pe_C^{\mathcal{Q}}(\mathcal{Z}_1, \ldots, \mathcal{Z}_{m'}) . \]
The specialization of variables $Z_i \leadsto \mathcal{Z}_{f(i)}$ allows us to write the $\mathcal{Q}$-enumerator in terms of the $\mathcal{P}$-enumerator:
\begin{equation}  \label{eq:spec-refinement}
\pe_C^{\mathcal{Q}}(\mathcal{Z}_1, \ldots, \mathcal{Z}_{m'}) = \left. \pe_C^{\mathcal{P}}(Z_1, \ldots, Z_m) \right|_{Z_i \leadsto \mathcal{Z}_{f(i)}} .
\end{equation}

The MacWilliams identities are known to generalize to the complete enumerator and certain symmetrized enumerators, over finite fields \cite{MacWilliamsSloane1977} and finite Frobenius rings \cite{wood:duality}.  The MacWilliams identities generalize to so-called reflexive partition enumerators over finite Frobenius rings; see \cite{MR3336966} for details.  
For the symmetrized enumerator, see Theorem~\ref{thm:MWforSE}.  
A short review of the main arguments used for proving the MacWilliams identities over finite Frobenius rings is in Appendix~\ref{sec:appendix}.

For the other type of enumerator,  
suppose 
 $R$ is a finite ring with $1$, and $w$ is a weight on $R$ with positive integer values (except $w(0)=0$).  Denote the largest value of $w$ by $w_{\max}$.   For any left $R$-linear code $C \subseteq R^n$, define the \emph{$w$-weight enumerator} of $C$ by
\begin{equation}  \label{eq:defnOfwwe}
\wwe_C(X,Y) = \sum_{x \in C} X^{n w_{\max} - w(x)} Y^{w(x)} . 
\end{equation}
The $w$-weight enumerator is a homogeneous polynomial of degree $n w_{\max}$ in $X$ and $Y$.  Different codewords in $C$ may have the same weight.  Collecting terms in \eqref{eq:defnOfwwe} leads to
\begin{equation}  \label{eq:A-formOfwwe}
\wwe_C(X,Y) = \sum_{j=0}^{n w_{\max}} A_j^w(C) X^{n w_{\max} - j} Y^j , 
\end{equation}
where $A_j^w(C)$ is the number of codewords of $C$ having weight $j$:
\begin{equation}  \label{eq:DefnOfAsubj}
A_j^w(C) = \size{\{x \in C: w(x) = j\}}.
\end{equation}
We write $A_j(C)$ when $w$ is clear from context.  To save space in examples in later sections we will often write $\wwe$ with $X=1$ and $Y=t$, so that $\wwe_C = \sum_j A_j(C) t^j$.  When $w = \wh$, the Hamming weight, we recover the Hamming weight enumerator $\hwe$.

\begin{rem}  \label{rem:appendZeroFuncl}
A disadvantage of using the notation $\wwe_C = \sum_j A_j(C) t^j$ is that information about the length $n$ of the code is lost.  Of course, if the length of $C$ is known, then the homogeneous form \eqref{eq:defnOfwwe} of $\wwe$ is easily recovered.  
For example, suppose $C$ is a linear code of length $n$, and let $D$ be the linear code of length $n+1$ obtained by appending a zero to each codeword of $C$.  Since $w(0)=0$, there are no changes in the weights of the codewords, so that $A_j(C) = A_j(D)$ for all $j$.  However, $\wwe_D(X,Y) = X \wwe_C(X,Y)$.  
\end{rem}

The partition enumerators and the $w$-weight enumerators are related.  Given a weight $w$ on $R$, let $\mathcal{Q}$ be the partition of $R$ into the orbits $\orb(r)$ of $G_{\lt}(w)$ acting on $R$ on the left, and let $\mathcal{P}$ be the complete partition of $R$; $\mathcal{P}$ is a refinement of every partition, hence a refinement of $\mathcal{Q}$.  Use variables $Z_r$, $r \in R$, for $\mathcal{P}$, and $\mathcal{Z}_{\orb(r)}$ for $\mathcal{Q}$.  Then the specialization of variables $\mathcal{Z}_{\orb(r)} \leadsto X^{w_{\max}-w(r)} Y^{w(r)}$ allows us to write $\wwe_C$ in terms of $\se_C^{\mathcal{Q}}$ for any $R$-linear code $C \subseteq R^n$:
\begin{equation}  \label{eq:spec-se-to-wwe}
\wwe_C(X,Y) = \left. \se_C^{\mathcal{Q}}(\mathcal{Z}_{\orb(r)}) \right|_{\mathcal{Z}_{\orb(r)} \leadsto X^{w_{\max}-w(r)} Y^{w(r)}} .
\end{equation}
This specialization is well-defined by the definition of the symmetry group $G_{\lt}$: the value of $w$ is constant on every left orbit of $G_{\lt}$.  Section~\ref{sec:examples} gives details of this situation over finite chain rings.

One way to view the MacWilliams identities is in terms of the dia\-gram in Figure~\ref{fig:comm-diagram} below.  
\begin{figure}[htb]
\[  \xymatrix{
\{ \text{$(n,M)$-linear codes}\} \ar[rr]^{\mathcal{R}} \ar[d]^{\ce} \ar@/_1.5pc/[dd]_{\se} \ar@/_3pc/[ddd]_{\wwe} && \{\text{$(n,\size{R}^n/M)$-linear codes}\} \ar[d]_{\ce} \ar@/^1.5pc/[dd]^{\se} \ar@/^3pc/[ddd]^{\wwe} \\
\C[Z_r]_n \ar[rr]^{MW} \ar@{~>}[d]^{\spec} && \C[Z_r]_n \ar@{~>}[d]_{\spec} \\
\C[\mathcal{Z}_{\orb(r)}]_n \ar[rr]^{MW} \ar@{~>}[d]^{\spec} && \C[\mathcal{Z}_{\orb(r)}]_n \ar@{~>}[d]_{\spec} \\
\C[X,Y]_{n w_{\max}} \ar@{.>}[rr]^{?} & & \C[X,Y]_{n w_{\max}}
}
\]
\caption{Relations among enumerators}  \label{fig:comm-diagram}
\end{figure}
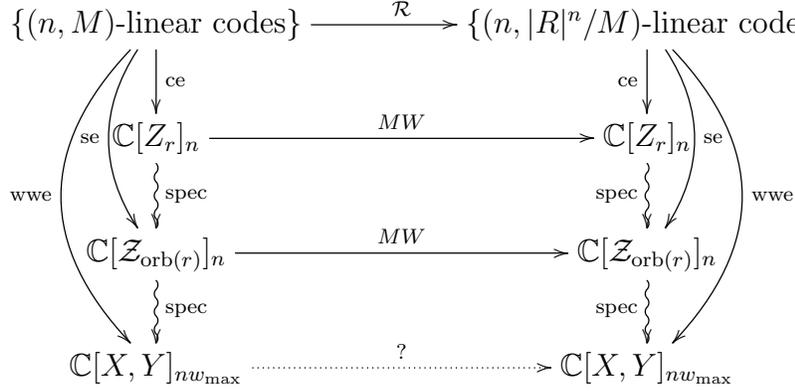
In the diagram, the map $\mathcal{R}$ sends a left $R$-linear code of size $M$ to its right $R$-linear dual code $\mathcal{R}(C)$.  Under favorable circumstances (e.g., $R$ Frobenius), the dual code has size $\size{\mathcal{R}(C)} = \size{R}^n/\size{C}$.  The vertical maps $\ce, \se, \wwe$ associate to a linear code its complete enumerator, $\mathcal{Q}$-symmetrized enumerator, and $w$-weight enumerator, respectively.  The other vertical maps (both called $\spec$) are the specializations of variables described in \eqref{eq:spec-refinement} and \eqref{eq:spec-se-to-wwe}.

Because the MacWilliams identities hold for $\ce$ and $\se$, the horizontal maps $MW$ are the MacWilliams transforms that provide the linear changes of variables.  
The solid arrows in the diagram commute.

The big question under study in this paper is whether there is a horizontal map `?' that makes the diagram commute for the $w$-weight enumerator.  If such a map exists, then the following property holds:  if $\wwe_C = \wwe_D$, then $\wwe_{\mathcal{R}(C)} = \wwe_{\mathcal{R}(D)}$.  We refer to this property by saying 
\emph{the weight $w$ respects duality}.  To formalize:
\begin{defn}
A weight $w$ on a finite ring $R$ \emph{respects duality} if $\wwe_C=\wwe_D$ implies $\wwe_{\mathcal{R}(C)} = \wwe_{\mathcal{R}(D)}$ for all left $R$-linear codes $C,D \subseteq R^n$, $n \geq 1$.
\end{defn}

Corollary~\ref{cor:respect-duality} says that the Hamming weight on $\F_q$ respects duality.  If a weight $w$ does not respect duality, then the MacWilliams identities cannot hold for $\wwe$.

As an example, let us see what happens when a weight is multiplied by a positive constant.
Suppose $w: R \ra \Z$ is a weight on $R$, and let $\tilde{w} = c w$, with $c$ a positive integer.  Denote the $w$-weight enumerators for $w$ and $\tilde{w}$ by $\wwe$ and $\wwe^c$ respectively.

\begin{lem}  \label{lem:scale-weight}
Let $w: R \ra \Z$ be a weight on $R$, and let $\tilde{w} = c w$, with $c$ a positive integer.  Then 
\[  \wwe_C^c(X,Y) = \wwe_C(X^c, Y^c) . \]
\end{lem}

\begin{proof}
For any element $r \in R$, we have 
\[  X^{\tilde{w}_{\max} - \tilde{w}(r)} Y^{ \tilde{w}(r)} = (X^c)^{w_{\max} - w(r)} (Y^c)^{w(r)} . \qedhere \]
\end{proof}

\begin{prop}  \label{prop:hold-weakly}
Let $w: R \ra \Z$ be a weight on $R$, and let $\tilde{w} = c w$, with $c$ a positive integer.  Then $\tilde{w}$ respects duality if and only if $w$ respects duality.
\end{prop}

Even if 
a weight $w$ does not respect duality (and hence $\wwe$ does not satisfy the MacWilliams identities), it is still possible to determine both $\wwe_C$ and $\wwe_{\mathcal{R}(C)}$ by calculating $\se_C$, using the MacWilliams identities for $\se$ to find $\se_{\mathcal{R}(C)}$, and then specializing variables to get $\wwe_C$ and $\wwe_{\mathcal{R}(C)}$.  But one cannot go directly from $\wwe_C$ to $\wwe_{\mathcal{R}(C)}$.

The main objective of this paper is to show that it is rare for a weight 
having maximal symmetry to respect duality, at least over finite chain rings or matrix rings over finite fields.  
In addition to Theorem~\ref{thm:MWids-original} for the Hamming weight enumerator and its generalization to finite Frobenius rings, the MacWilliams identities for $w$-weight enumerators are known to hold for the Lee weight on $\Z/4\Z$ \cite{Hammons-etal} (and see Theorem~\ref{thm:chainMacWids}) and the homogeneous weight on the matrix ring $M_{2 \times 2}(\F_2)$, Theorem~\ref{thm:MWIds-homog-2by2-F2}.  The MacWilliams identities for $w$-weight enumerators are known to fail for the Rosenbloom-Tsfasman weight on matrices \cite{MR1900585}, the Lee weight on $\Z/m\Z$, $m \geq 5$ \cite{MR4119402}, and the homogeneous weight on $\Z/m\Z$ for composite $m \geq 6$ \cite{revuma.2807}.  In all cases, the failure is proved by showing that the weight does not respect duality: there exist linear codes $C$ and $D$ with $\wwe_C = \wwe_D$, yet $\wwe_{\mathcal{R}(C)} \neq \wwe_{\mathcal{R}(D)}$, by virtue of $A_j^w(\mathcal{R}(C)) \neq A_j^w(\mathcal{R}(D))$ for some $j$.

The hypothesis that $w$ has maximal symmetry is important.   There are results about `Lee weights' of different types that can be valid because they secretly tap into the Hamming weight or the Lee weight on $\Z/4\Z$; cf., \cite{MR3168981}.

\section{Linear codes via multiplicity functions}  \label{sec:W-stuff}
In later sections linear codes will be presented as images of homomorphisms of left $R$-modules.  In turn, the homomorphisms will be described in terms of multiplicity functions.  In this section, the use of multiplicity functions to describe linear codes will be summarized briefly.  Addititonal information can be found in \cite[\S 3]{wood:structure} or \cite[\S 7]{wood:turkey}.

As throughout this paper, let $R$ be a finite ring with $1$.  Suppose $M$ is a finite unital left $R$-module; \emph{unital} means $1x=x$ for all $x \in M$.  We call any homomorphism $\la: M \ra R$ of left $R$-modules a \emph{linear functional} on $M$.  Define $M^\sharp = \Hom_R(M,R)$ to be the set of all left linear functionals on $M$.  We will write inputs to linear functionals on the left, so that $(rx)\la = r(x\la)$ for $r \in R$, $x \in M$, and $\la \in M^\sharp$; $M^\sharp$ is a right $R$-module, with $\la r$ given by $x(\la r) = (x \la) r$ for $r \in R$, $x \in M$, and $\la \in M^\sharp$.  The left $R$-module $M$ admits a left action by the group of units $\U$ of $R$ using left scalar multiplication.  Denote the orbit of $x \in M$ by $\orb(x)$ or by $[x]$.  Similarly, the right $R$-module $M^\sharp$ admits a right $\U$-action, with orbits denoted $\orb(\la)$ or $[\la]$.

A left $R$-linear code of length $n$ \emph{parametrized} by $M$ is the image $C = \im \La$ of a homomorphism $\La: M \ra R^n$ of left $R$-modules. The module $M$ is the \emph{information module} of the linear code $C$.  Denote the components of $\La$ by $\La = (\la_1, \la_2, \ldots, \la_n)$, with each $\la_i: M \ra R$ being linear functional on $M$.
We call the $\la_i$ the \emph{coordinate functionals} of the linear code $C$.

Suppose $w$ is a weight on $R$.  Then the weight $w$ and a parametrized code $C$ given by $\La: M \ra R^n$ define a \emph{weight function} $W_{\La}: M \ra \C$ by $W_{\La}(x) = w(x \La) = \sum_{i=1}^n w(x\la_i)$, $x \in M$.

\begin{lem}  \label{lem:EffectsSymmetry}
Suppose $w$ is a weight on $R$ with symmetry groups \eqref{eq:SymmetryGroups}, and suppose $C$ is an $R$-linear code parametrized by $\La: M \ra R^n$.  Then
\begin{enumerate}
\item the weight function $W_{\La}: M \ra \C$ is constant on each left $\glt(w)$-orbit $\orb(x) \subseteq M$;
\item  $W_{\La'} = W_{\La}$, if $\La' = (\la_{\tau(1)} u_1, \ldots, \la_{\tau(n)} u_n)$, where $\tau$ is a permutation of $\{ 1, 2, \ldots, n\}$ and $u_1, \ldots, u_n \in \grt(w)$.
\end{enumerate}
\end{lem}

\begin{proof}
These follow directly from \eqref{eq:SymmetryGroups}.
\end{proof}

\begin{rem}  \label{rem:listToWwe}
If the weight $w$ has nonnegative integer values, the weight function $W_{\La}$ determines the $w$-weight enumerator of the linear code $C$:
\[  \wwe_C = \frac{1}{\size{\ker \La}} \sum_{x \in M} t^{W_{\La}(x)} . \]
When $\La$ is injective, the $w$-weight enumerator can be written in terms of the sizes of orbits:
\begin{equation}  \label{wwe-in-terms-of-orbits}
\wwe_C = \sum_{[x] \subseteq M} \size{[x]} t^{W_{\La}(x)} .
\end{equation}
\end{rem}

Lemma~\ref{lem:EffectsSymmetry} shows that the weight function $W_\La$ depends only on the numbers of coordinate functionals belonging to different $\grt(w)$-orbits in $M^\sharp$.  We formalize this observation next.  Write $F(X,Y)$ for the set of all functions from $X$ to $Y$.

Given an information module $M$, define the orbit spaces:
\begin{align*}
\OO = G_{\lt}(w) \backslash M , \quad
\OS   = M^\sharp/G_{\rt}(w) .
\end{align*}
Lemma~\ref{lem:EffectsSymmetry} shows that $W_\La$ depends only on a function $\eta \in F(\OS, \N)$, where $\eta([\la])$ is the number of coordinate functionals that belong to the orbit $\orb(\la) = [\la]$.  (In \cite[\S 3.6]{MR49:12164}, Peterson and Weldon call $\eta$ the \emph{modular representation} of the linear code $C$.) 

The weight $w$ induces an additive map $W: F(\OS, \N) \ra F(\OO, \C)$.  For $\eta \in F(\OS, \N)$, define $\omega = W(\eta)$, the \emph{list of orbit weights}, by 
\begin{equation}  \label{eqn:defn-omega}
\omega([x]) = \sum_{[\la] \in \OS} w(x \la) \eta([\la]) .
\end{equation}
By \eqref{eq:SymmetryGroups}, the value of $w(x \la)$ is well-defined.  If $w$ has values in $\Z$ or $\Q$, then $W$ has values in $F(\OO, \Z)$ or $F(\OO, \Q)$, accordingly.  In these latter cases, 
tensoring with $\Q$ yields a linear transformation $W: F(\OS, \Q) \ra F(\OO, \Q)$.

By ordering the elements of $\OO$ and $\OS$, one can define a matrix $W$ whose rows are indexed by $\OO$, whose columns are indexed by $\OS$, and whose entry at position $([x], [\la])$ is the well-defined value $w(x \la)$:
\begin{equation}  \label{eq:DefnOfWmatrix}
W_{[x],[\la]} = w(x \la), \quad [x]\in \OO, [\la] \in \OS.
\end{equation}
Treating $\eta$ and $\omega$ as column vectors, \eqref{eqn:defn-omega} is just matrix multiplication: $\omega = W \eta$.

Any element $\eta \in F(\OS, \N)$, called a \emph{multiplicity function}, determines an $R$-linear code as the image of the homomorphism $\La_{\eta}: M \ra R^N$ of left $R$-modules given by sending $x \in M$ to the $N$-tuple $(\ldots, x \la, \ldots)$, where a representative of each orbit $[\la] \in \OS$ is repeated $\eta([\la])$ times; $N = \sum_{[\la] \in \OS} \eta([\la])$.  Said another way, treat elements $\la \in M^\sharp$ as columns of a generator matrix, with representatives of $[\la]$ repeated $\eta([\la])$ times.  The resulting linear code $C_{\eta}$ is well-defined up to monomial equivalence; its list $\omega = W \eta$ of orbit weights is well-defined.  Using $\omega$, one can write down the $w$-weight enumerator of $C_{\eta}$, as in Remark~\ref{rem:listToWwe}:
\[  \wwe_{C_{\eta}} = \sum_{x \in M} t^{w(x \La_{\eta})} = \sum_{[x] \in \OO} \size{[x]} t^{\omega([x])} , \]
assuming that $\La_{\eta}$ is injective.  (If $\La_{\eta}$ is not injective, we must divide by $\size{\ker \La_{\eta}}$.)

\begin{rem}  \label{rem:W0-map}
The definition of $W: F(\OS,\N) \ra F(\OO,\C)$ is valid whether or not $w(0)=0$.  When $w(0)=0$, then $W(\eta)([0]) = \omega([0])= 0$ for any $\eta \in F(\OS,\N)$.  Also, when $w(0)=0$, $W([\delta_{[0]}]) = 0$, where $\delta_{[0]}$ is the indicator function of the orbit $[0] \in \OS$ of the zero-functional.

More formally, define
\begin{align*}
F_0(\OS, \N) &= \{ \eta \in F(\OS,\N): \eta([0]) = 0 \}, \\
F_0(\OO, \C) &= \{ \omega \in F(\OO,\C): \omega([0]) = 0 \} .
\end{align*}
When $w(0)=0$, 
the image of $W$ is contained in $F_0(\OO, \C)$.  We denote the restriction of
$W$ to $F_0(\OS, \N)$ by $W_0: F_0(\OS, \N) \ra F_0(\OO, \C)$; $W_0$ is an additive map.  The map $W_0$ is represented by a matrix whose rows are indexed by the nonzero orbits $[x] \in \OO$, whose columns are indexed by the nonzero orbits $[\la] \in \OS$, and with entries given as in \eqref{eq:DefnOfWmatrix}.  When $w(0)=0$, the map $W$ can never be injective (because $\delta_{[0]} \in \ker W$), but the map $W_0$ often is injective.  When $w$ has values in $\Q$ and $w(0)=0$, the linear transformation $W_0: F_0(\OS, \Q) \ra F_0(\OO, \Q)$ is often invertible.  This will be an important tool in later sections.
\end{rem}

We conclude this section with a short discussion of the effective length of codes. 
The next lemma is a variant of \cite[(6.1)]{wood:structure}. 

\begin{lem}  \label{lem:SumOfWeights}
Suppose $C \subseteq R^n$ is an $R$-linear code.  Then
\begin{equation}  \label{eq:weight-sum-a}
\sum_{c \in C} w(c) = \sum_{i=1}^n \size{\ker \la_i} \sum_{b \in \im\la_i} w(b) .
\end{equation}
\end{lem}

\begin{proof}
Write $\la_1, \ldots, \la_n$ for the coordinate functionals of $C$.  Then, $\sum_{c \in C} w(c) = \sum_{c \in C} \sum_{i=1}^n w(c_i) = \sum_{c \in C} \sum_{i=1}^n w(c \la_i)$.  Now interchange the finite sums, and use that $\la_i$ is a homomorphism, so $\sum_{c \in C} w(c \la_i) = \size{\ker \la_i} \sum_{b \in \im\la_i} w(b)$.
\end{proof}

There are situations where \eqref{eq:weight-sum-a} simplifies.  Define a weight $w$ on $R$ to be \emph{egalitarian} if there exists a constant $\ga$ such that, for any nonzero left ideal $B \subseteq R$, $\sum_{b \in B} w(b) = \ga \size{B}$.  This definition is due to \cite{Honold-Nechaev:weighted-modules}.  The homogeneous weight on a finite Frobenius ring is an example of an egalitarian weight; see Example~\ref{ex:homogeneous}.

Define the \emph{effective length} 
of a linear code $C$ to be $\efflength(C)=\size{\{ i : \la_i \neq 0 \} }$.  
If $C$ is given by a generator matrix, the effective length counts the number of nonzero columns of the generator matrix.  

\begin{prop}  \label{prop:essentialLength}
Suppose $C \subseteq R^n$ is an $R$-linear code  
and $w$ is an egalitarian weight on $R$.  Then
\[  \sum_{c \in C} w(c) = \ga \size{C} \efflength(C) . \]
\end{prop}

\begin{proof}
Use Lemma~\ref{lem:SumOfWeights}.  For $\la_i \neq 0$, set $B = \im\la_i$ and note that $\size{C} = \size{\ker\la_i}\cdot\size{\im\la_i}$.
\end{proof}

\begin{cor}  \label{cor:sameEssentialLength}
Let $w$ be an egalitarian weight on $R$.  Suppose $C, D \subseteq R^n$ are two $R$-linear codes.  If $\wwe_C = \wwe_D$, then the effective lengths of $C$ and $D$ are equal.
\end{cor}

\begin{proof}  
The hypothesis means that $A_j^w(C) = A_j^w(D)$ for all $j$.   But
\begin{align*}
\size{C} = \sum_{j=0}^{n w_{\max}} A_j^w(C)  \quad \text{and} \quad
\sum_{c \in C} w(c) = \sum_{j=0}^{n w_{\max}} j A_j^w(C) ,
\end{align*}
so $\size{C} = \size{D}$ and $\sum_{c \in C} w(c) = \sum_{d \in D} w(d)$.  Apply Proposition~\ref{prop:essentialLength}.
\end{proof}

\section{Singletons in dual codes}
In this section we describe some general results on the contributions to $\wwe_{C^\perp}$ coming from singleton vectors.

As usual, suppose $R$ is a finite ring with $1$ and $C \subseteq R^n$ is a left $R$-linear code.  Its dual code $C^\perp$ is the right dual code $\mathcal{R}(C)$ of \eqref{eq:RightDualCode}.  We assume $w$ is an integer-valued weight on $R$ with $w(r)>0$ for $r \neq 0$ and $w(0)=0$.   
Let $\mathring{w} = \min\{w(r): r \neq 0 \}$, so that $\mathring{w} >0$.  

We say that a vector $v \in R^n$ is a \emph{singleton} if $v$ has exactly one nonzero entry.  
Given a vector $v \in R^n$, recall that the weight of the vector is $w(v) = \sum_{j=1}^n w(v_j)$.  The smallest possible nonzero weight of a vector is $\mathring{w}$, which is attained by any singleton whose nonzero entry $r$ has $w(r)=\mathring{w}$.

We want to write down the contributions of singletons to the $w$-weight enumerator of a linear code, especially to a dual code.  As in \eqref{eq:DefnOfAsubj}, recall that $A_j(C) = \size{\{ x \in C: w(x)=j \}}$.  To track the contributions of singletons we write \[  \asing_j(C) = \size{\{ x \in C: \text{$x$ is a singleton and $w(x)=j$} \}} . \]
Of course, $\asing_j(C) \leq A_j(C)$.  Equality will be addressed in Corollary~\ref{cor:only-singletons} below.

For any $\la \in M^\sharp$ and positive integer $j$, define $\ann_{\rt}(\la, j) = \{ r \in R: \la r = 0 \text{ and } w(r)=j \}$, the set of elements in $R$ of weight $j$ that annihilate $\la$. 

\begin{lem}  \label{lem:SizeAnnWellDefined}
Suppose $w$ has maximal symmetry.  For any $\la \in M^\sharp$ and $u \in \U$, $\ann_{\rt}(\la u, j) = u^{-1} \ann_{\rt}(\la, j)$.
In particular, $\size{ \ann_{\rt}(\la u, j)} = \size{\ann_{\rt}(\la, j)}$ for any $u \in \U$, $\la \in M^\sharp$.
\end{lem}

\begin{proof}
Suppose $w(r)=j$.  By maximal symmetry, $w(ur)=w(r)=j$ for all $u \in \U$.
Because $(\la u)(u^{-1} r) = \la r$, we see that $r \in \ann_{\rt}(\la,j)$ if and only if $u^{-1} r \in \ann_{\rt}(\la u, j)$.
\end{proof}

\begin{prop}  \label{prop:SingletonContribution}
Assume $w$ is an integer-valued weight on $R$ with maximal symmetry.  If $C$ is a linear code determined by a multiplicity function $\eta$, then, for any positive integer $j$,
\[  \asing_j(C^\perp) = \sum_{[\la] \in \OS} \size{\ann_{\rt}(\la, j)} \eta([\la]) . \]
\end{prop}

\begin{proof}
Suppose $C$ has coordinate functionals $\la_1, \ldots, \la_n$.  Let $v$ be a singleton vector with nonzero entry $r \in R$ appearing in position $i$.  Then $v \in C^\perp$ if and only if $\la_i r = 0$.  Thus
\[  \asing_j(C^\perp) = \sum_{i=1}^n \size{\ann_{\rt}(\la_i,j)} , \]
which reduces to the stated formula because of Lemma~\ref{lem:SizeAnnWellDefined}.
\end{proof}

In later sections, formulas for $\size{\ann_{\rt}(\la, j)}$ will be very specific, depending on the nature of the ring $R$.

\begin{lem}  \label{lem:singletons}
Suppose $v \in R^n$ has weight $w(v)$ satisfying $\mathring{w} \leq w(v) < 2\mathring{w}$.  Then $v$ must be a singleton.
\end{lem}

\begin{proof}
Suppose $v$ has at least two nonzero entries, say in positions $j_1, j_2$.  Then $w(v) \geq w(v_{j_1}) + w(v_{j_2}) \geq 2 \mathring{w}$.
\end{proof}

\begin{cor}  \label{cor:only-singletons}
If $\mathring{w} \leq d < 2 \mathring{w}$, then $A_d(C) = \asing_d(C)$.  
\end{cor}

In later sections, Corollary~\ref{cor:only-singletons} will be applied mostly to dual codes, in tandem with Proposition~\ref{prop:SingletonContribution}.

\begin{rem}  \label{rem:NotInLattice}
In order that $\asing_j(C)$ be nonzero, it is necessary that $j = w(r)$ for some $r \in R$.

It is possible that $A_j(C) = \asing_j(C)$ even when $j \geq 2 \mathring{w}$.  For example: when $j = w(r)$ is not equal to a linear combination of the form $\sum_{s: w(s) < j} c_s w(s)$ with $c_s$ being nonnegative integers.
\end{rem}

\part{Finite Chain Rings}

\section{Definitions and a positive result}  \label{sec:ChainRings}
A finite ring $R$ with $1$ is a \emph{chain ring} if its left ideals form a chain under set inclusion.  In particular, $R$ has a unique maximal left ideal, denoted $\m$, so that $R$ is a local ring.  Examples of chain rings include finite fields, $\Z/p^m\Z$ with $p$ prime, Galois rings, $\F_q[X]/(X^m)$; cf., \cite{MR0354768}.  Every finite chain ring is Frobenius \cite[Lemma~14]{MR1634126}.

From \cite[Lemma~1]{ClarkDrake:chain-rings} we know that $\m$ is a principal ideal, say $\m =R\theta = \theta R$, that $\theta^m =0$ for some (smallest) $m \geq 1$, and that every left or right ideal of $R$ is a two-sided ideal of the form $R \theta^j = \theta^j R$, $j=0,1, \ldots, m$.  In particular, $\m$ is a two-sided ideal, so that $R/\m$ is a finite field, say $R/\m \cong \F_q$, of order $q$, a prime power.  Write $(\theta^j)$ for $R \theta^j = \theta^j R$.  Thus, all the ideals of $R$ are displayed here:
\begin{equation}  \label{eq:chainofideals}
R=(\theta^0) \supset (\theta) \supset (\theta^2) \supset \cdots \supset (\theta^{m-1}) \supset (\theta^m)=(0) .
\end{equation}
Each quotient $(\theta^j)/(\theta^{j+1})$ is a one-dimensional vector space over $R/\m$, with basis element $\theta^j + (\theta^{j+1})$.  It follows that
\begin{equation}  \label{eq:chainidealsize}
\size{(\theta^j)} = q^{m-j}, \quad j=0,1, \ldots, m .
\end{equation}
In particular, $\size{R}=q^m$.  The group of units of $R$, denoted $\U = \U(R)$, equals the set difference $R - \m$.  The group $\U$ acts of $R$ on the left and on the right by multiplication.  The orbits of the actions are exactly the set differences $\orb(\theta^j) = (\theta^j) - (\theta^{j+1})$, 
which have size 
\begin{equation}  \label{eq:OrbitSizes}
\size{\orb(\theta^j)} = q^{m - j -1}(q-1), \quad \text{for $j < m$}.  
\end{equation}
In particular, the left orbits of $\U$ equal the right orbits: $\U \theta^j = \theta^j \U$.

From \eqref{eq:chainofideals}, we see that
every element $r \in R$ has the form $r = u \theta^j$ where $u$ is a unit of $R$  and $j$ is uniquely determined by $r$ (the largest $i$ such that $r \in (\theta^i)$).  Note that the annihilator of $(\theta^j)$ is $(\theta^{m -j})$.  

Let $w$ be a weight on $R$ with positive integer values for $r \neq 0$ in $R$.  Assume that $w$ has maximal symmetry, so that $w(ur)=w(ru)=w(r)$ for all $r \in R$ and units $u \in \U$.  This means that $w$ is constant on the $\U$-orbits $\orb(\theta^j) = (\theta^j) - (\theta^{j+1})$.  Define $w_j$ as the common value of $w$ on $\orb(\theta^j)$, so that $w_j = w(u \theta^j) = w(\theta^j u)$ for all units $u \in \U$.  Then $w_0, w_1, \ldots, w_{m-1}$ are positive integers, and $w_m = 0$.

\begin{ex}
Choosing $\zeta = q-1$, we see from Example~\ref{ex:homogeneous} and \eqref{eq:chainofideals} that the homogeneous weight $\wg$ on a chain ring $R$ has the following integer values:
\[  \wg(r) = \begin{cases}
0, & r=0, \\
q, & r \in (\theta^{m-1})-(0), \\
q-1, & r \in R-(\theta^{m-1}).
\end{cases} \]
Then $\wg_0 = \cdots = \wg_{m-2} = q-1$, $\wg_{m-1}=q$, and $\wg_m=0$.
\end{ex}

Do the MacWilliams identities hold for the homogeneous weight enumerator over a finite chain ring $R$?  We will see that the answers depend on $q$ and $m$.

When $m=1$, then $\theta=0$, so that $R$ is a finite field $\F_q$.  As we saw in Example~\ref{ex:homogeneous}, the homogeneous weight $\wg$ on $\F_q$ equals a multiple of the Hamming weight.  By Theorems~\ref{thm:MWids-original} and Lemma~\ref{lem:scale-weight}, the homogeneous weight over finite fields respects duality.

For $m \geq 2$, there is one special case ($m=q=2$, Theorem~\ref{thm:chainMacWids}, below) where the MacWilliams identities hold for the homogeneous weight.  In the remaining cases, we will see in Theorem~\ref{thm:MainThmWeaklyMonotone} that the homogeneous weight does not respect duality.

The results just described apply to the chain rings $\Z/p^m\Z$, so that the MacWilliams identities hold for the homogeneous weight over $\Z/p\Z$, $p$ prime, and over $\Z/4\Z$ \cite[Equation~(9)]{Hammons-etal}, but not for other prime powers.
More generally, over $\Z/m\Z$, $m$ not a prime power (so that $\Z/m\Z$ is not a chain ring), the homogeneous weight does not respect duality \cite[Theorem~6.2]{revuma.2807}.

\begin{thm}  \label{thm:chainMacWids}
The MacWilliams identities hold for the homogeneous weight enumerator over a finite chain ring $R$ with $q=2$ and $m=2$.  If $C \subseteq R^n$ is a linear code and $C^\perp$ is its dual code, then
\[  \howe_{C^\perp}(X,Y) = \frac{1}{\size{C}} \howe_C(X+Y, X-Y) . \]
\end{thm}

\begin{proof}
Appendix~\ref{sec:appendix} outlines of a proof of the MacWilliams identities over finite Frobenius rings and describes the Fourier transform.  Here, we provide details relevant to the chain rings appearing in this theorem. 

We know that $\size{R} = 4$, with $R=\{0, 1, \theta, 1+\theta\}$.  The values of the homogeneous weight, with $\zeta =1$, are:
\[  \begin{array}{c|cccc}
r & 0 & 1 & \theta & 1+\theta \\ \hline
\wg(r) & 0 & 1 & 2 & 1 
\end{array} . \]
The additive group of $R$ could be a cyclic group of order $4$ (in which case $\theta=2$, in order that $(\theta)$ be a maximal ideal) or a Klein $4$-group.   In either case, there exists a generating character $\chi$ of $R$ with the following values.  (What is crucial is that $\chi(\theta)=-1$.)
\[  \begin{array}{lc|cccc}
& r & 0 & 1 & \theta & 1+\theta \\ \hline
\text{cyclic} & \chi(r) & 1 & i & -1 & -i \\
\text{Klein} & \chi(r) & 1 & 1 & -1 & -1
\end{array} \]

Define $f: R \ra \C[X,Y]$ by $f(r) = X^{2-\wg(r)} Y^{\wg(r)}$.  For either choice of $\chi$, the Fourier transform $\ft{f}$ of \eqref{eqn:FTdefn} is the same:
\[  \begin{array}{c|cccc} r & 0 & 1 & \theta & 1+\theta \\ \hline
f(r) & X^2 & XY & Y^2 & XY \\
\ft{f}(r) & (X+Y)^2 & X^2 - Y^2  & (X-Y)^2 & X^2 - Y^2 \\
\end{array} \]
Note that the values of $\ft{f}$ have the form $\ft{f}(r)=(X+Y)^{2-\wg(r)} (X-Y)^{\wg(r)}$.  I.e., $\ft{f}(r) = f(r)|_{X \leftarrow X+Y, Y \leftarrow X-Y}$.
The rest of the argument in Appendix~\ref{sec:appendix} now carries through.
\end{proof}

\section{Modules over a chain ring}

We continue to assume that $R$ is a finite chain ring with maximal ideal $\m = (\theta)$ such that $R/\m \cong \F_q$ and $\theta^m = 0$, for some integer $m \geq 2$.
By forming quotients from \eqref{eq:chainofideals}, we define cyclic $R$-modules $Z_k = R/(\theta^k)$, of order $q^k$, $k=1, 2, \ldots, m$.  The module $Z_1$ is the unique simple $R$-module.  Also define semisimple modules $S_k = Z_1 \oplus \cdots \oplus Z_1$ with $k$ summands, $k= 2, 3,  \ldots, m$; $\size{S_k} = q^k$.  

In the next two sections, we will construct two families of linear codes with the same $w$-weight enumerators.
Here is a brief sketch of the construction.  We will build examples using $Z_k$ and $S_k$, $2 \leq k \leq m$, as the underlying information modules.  On $Z_k$, choose a multiplicity function.  The weights of elements will be constant along nonzero orbits $\OO$, which have sizes $q^k - q^{k-1}, q^{k-1} - q^{k-2}, \ldots, q^2 - q, q-1$.  We then need to build a linear code based on $S_k$ with the same $w$-weight enumerator.  This entails choosing subsets of $S_k$ with sizes matching the sizes $q^k - q^{k-1}, q^{k-1} - q^{k-2}, \ldots, q^2 - q, q-1$.  We do this by fixing a filtration of $S_k$ using linear subspaces.
We can then solve for a multiplicity function for $S_k$.  If the original multiplicity function on $Z_k$ calls for one functional from each class, then the multiplicity function on $S_k$ is reasonably nice.  We then calculate the common weight enumerator for these multiplicity functions.

In this section we study the orbit structure of $Z_k$ as well as the subsets arising from a filtration of $S_k$.  In Section~\ref{sec:twoFamiliesChain}, the multiplicity functions are described and analyzed.

In order to define linear codes over $Z_k$ and $S_k$, let us examine their linear functionals.  
Recall first that the left linear functionals of $R$ itself are given by right multiplications by elements of $R$.  That is, $R^\sharp \cong R$ as right $R$-modules, with $r \in R$ corresponding to the left linear functional $\rho_r \in R^\sharp$ defined by $r' \rho_r = r' r$, $r' \in R$.

If $\la \in Z_k^\sharp$, then the composition with the natural quotient map must equal $\rho_r$ for some $r \in R$.  
\[  \xymatrix{
R \ar@{.>}[rr]^{\rho_r} \ar@{->>}[d] & & R  \\
Z_k = R/(\theta^k) \ar[urr]_{\la}
}
\]
Conversely, $\rho_r: R \ra R$ factors through $Z_k$ if and only if $(\theta^k) \subseteq \ker \rho_r$; i.e., if and only if $\theta^k r =0$.  This occurs when $r \in (\theta^{m-k})$.  Thus $Z_k^\sharp \cong (\theta^{m-k})$ as right $R$-modules.  
In particular, $Z_1^\sharp \cong (\theta^{m-1})$ as right $R$-modules.

As for $S_k = Z_1 \oplus \cdots \oplus Z_1$, we have $S_k^\sharp \cong (\theta^{m-1}) \oplus \cdots \oplus (\theta^{m-1})$ as right $R$-modules.  For $s = (s_1, \ldots, s_k) \in S_k$ and $\mu = \langle \mu_1 \theta^{m-1}, \ldots, \mu_k \theta^{m-1} \rangle \in S_k^\sharp$, $s \mu = \sum_{i=1}^k s_i \mu_i \theta^{m-1} \in R$.  Both $S_k$ and $S_k^\sharp$ are $k$-dimensional vector spaces over $R/\m \cong \F_q$.

In order to exploit maximal symmetry in Lemma~\ref{lem:EffectsSymmetry}, we want to understand the orbit structures of the left actions of $\U$ on $Z_k$ and $S_k$, and, to a lesser extent, the orbit structures of the right actions of $\U$ on $Z_k^\sharp$ and $S_k^\sharp$.

Because $(\theta^k)$ is a two-sided ideal of $R$, $Z_k$ is itself a chain ring.  Its left $R$-submodules are the same as its left ideals:
\begin{equation}  \label{eq:FilterZk}
Z_k \supset R \theta \supset R \theta^2 \supset \cdots \supset R \theta^{k-1} \supset \{0\} .
\end{equation}
Note that \eqref{eq:FilterZk} can be viewed as a filtration of $Z_k$ by $R$-submodules. 

\begin{lem}  \label{lem:Zk-orbits}
The orbits of the left action of $\U$ on $Z_k$ are: 
\[  \orb(1) = \orb(\theta^0), \orb(\theta), \ldots, \orb(\theta^{k-1}), \{0\} = \orb(0) = \orb(\theta^k) . \]
The sizes of the orbits are:  $\size{\orb(0)} = 1$ and
\[  \size{\orb(\theta^i)} = q^{k-i-1}(q-1), \quad i=0, 1, \ldots, k-1 . \]
\end{lem}

\begin{proof}
Apply \eqref{eq:OrbitSizes} to the chain ring $Z_k$.
\end{proof}

Similar to Lemma~\ref{lem:Zk-orbits}, the orbits of the right action of $\U$ on $Z_k^\sharp$ are $\orb(\theta^{m-k}), \orb(\theta^{m-k+1}), \ldots, \orb(\theta^{m-1}), \{0\} = \orb(\theta^m)$.

Because $S_k$ is a vector space over $\F_q \cong R/\m$, the action of $\U$ on $S_k$ reduces to the action of the multiplicative group $\F_q^\times$.  The $\U$-orbits are $\{0\}$ and $L-\{0\}$, for every $1$-dimensional subspace $L \subseteq S_k$.  The same structure applies to the dual vector space $S_k^\sharp$; the $\U$-orbits are $\{0\}$ and the nonzero elements of $1$-dimensional subspaces.

For later use in Section~\ref{sec:twoFamiliesChain}, we will identify subsets (say, $\mathcal{S}_1, \ldots, \mathcal{S}_k$) of $S_k$ that will be the counterparts to the orbits $\orb(1), \ldots, \orb(\theta^{k-1})$ of $Z_k$ (in reverse order).  In particular, we want two features:
\begin{itemize}
\item each $\mathcal{S}_i$ is a union of $\U$-orbits in $S_k$;
\item  $\size{\mathcal{S}_i} = \size{\orb(\theta^{k-i})} = q^{i-1}(q-1)$ for $i=1, \ldots, k$.
\end{itemize}

To define the $\mathcal{S}_i$, 
we define a filtration on $S_k$.  Recall that $S_k$ is a $k$-dimensional vector space over $R/\m \cong \F_q$.  Write elements of $S_k = Z_1\oplus \cdots \oplus Z_1$ as row vectors of length $k$ over $\F_q$.  
(Row vectors will be written as $(x_1, \ldots, x_k)$, while column vectors will be written as $\langle x_1, \ldots, x_k \rangle$.)  
Define vector subspaces of $S_k$: $V_0 = \{0\}$ and $V_i = \{ (0, \ldots, 0, s_{k-i+1}, \ldots, s_k) \in S_k: s_j \in \F_q\}$, for $i=1, 2, \ldots, k$.  In $V_i$, the first $k-i$ entries are zero; the last $i$ entries vary over $\F_q$.  Then, $\dim_{\F_q} V_i = i$ for $i = 0, 1, \ldots, k$, so that $\size{V_i} = q^i$, and 
\begin{equation}  \label{eq:FilterSk}
S_k = V_k \supset V_{k-1} \supset \cdots \supset V_1 \supset V_0 = \{0\}.
\end{equation}
Set $\mathcal{S}_i = V_i - V_{i-1}$, for $i=1, \ldots, k$; then $\size{\mathcal{S}_i} = q^{i-1}(q-1)$.  Set $\mathcal{S}_0 = \{0\}$.  Because the $V_i$ are vector subspaces, the $\mathcal{S}_i$ are unions of $\U$-orbits.

We also want to understand the linear functionals on $S_k$ in terms of the filtration \eqref{eq:FilterSk}. To that end, we examine the dual filtration of $S_k^\sharp$ defined by annihilators $\mathcal{V}_i = \ann(V_i) = \{ \mu \in S_k^\sharp: V_i \mu = 0\}$.  Then, $\dim_{\F_q} \mathcal{V}_i = k - \dim_{\F_q} V_i = k-i$, $\size{\mathcal{V}_i} = q^{k-i}$, and 
\begin{equation}  \label{eq:FilterSksharp}
\{0\} = \mathcal{V}_k \subset \mathcal{V}_{k-1} \subset \cdots \subset \mathcal{V}_1 \subset \mathcal{V}_0 = S_k^\sharp .
\end{equation}
If we view elements of $S_k^\sharp$ as column vectors $\mu = \langle \mu_1 \theta^{m-1}, \ldots, \mu_k \theta^{m-1} \rangle$, then $\mathcal{V}_i$ consists of those $\mu$ whose last $i$ entries equal zero.  The zero entries of $\mathcal{V}_i$ align with the nonzero entries of elements of $V_i$.

The set differences $\mathcal{V}_i - \mathcal{V}_{i+1}$ consist of all $\mu = \langle \mu_1 \theta^{m-1}, \ldots, \mu_k \theta^{m-1} \rangle$ with $\mu_{k-i+1}= \cdots = \mu_k=0$, $\mu_{k-i} \neq 0$, and $\mu_1, \ldots, \mu_{k-i-1} \in \F_q$.  The $\U$-orbit of $\mu$ is the set of all nonzero scalar multiplies of $\mu$.  In each orbit there is exactly one `normalized' representative with $\mu_{k-i}=1$.  For $i=0, 1, \ldots, k-1$, let $B_i$ be the subset of $\mathcal{V}_i - \mathcal{V}_{i+1}$ consisting of all the normalized representatives; i.e.,
\begin{equation}  \label{defn:Bi}
B_i = \{ \mu \in S_k: \mu_{k-i} =1, \mu_{k-i+1}= \cdots = \mu_k=0 \};
\end{equation}
as a special case, set $B_k = \{0\}$.  Then $\size{B_i} = q^{k-i-1}$, except $\size{B_k}=1$.

\begin{lem}  \label{lem:IntersectionCount}
Let $s \in \mathcal{S}_i$.  Then, for $j=0, 1, \ldots, k-1$, 
\[  \size{B_j \cap \ann(s)} = \begin{cases}
q^{k-j-1} , & i \leq j \leq k-1 , \\
0 , & j = i-1, \\
q^{k-j-2} , & 0 \leq j \leq i-2 .
\end{cases} \]
\end{lem}

\begin{proof}
The case $\mathcal{S}_0 = \{0\}$ has $\ann(0) = S_k^\sharp$.  Then $\size{B_j \cap \ann(0)} = \size{B_j} = q^{k-j-1}$ for all $0 \leq j \leq k-1$. 

Now let $1 \leq i \leq k$.  The element $s \in \mathcal{S}_i = V_i-V_{i-1}$ is nonzero and has the form $s=(0, \ldots, 0, s_{k-i+1}, \ldots, s_k)$, with $s_{k-i+1} \neq 0$.  Any $\mu \in \mathcal{V}_{i-1} - \mathcal{V}_{i}$ has the form $\mu = \langle \mu_1 \theta^{m-1}, \ldots, \mu_{k-i+1} \theta^{m-1}, 0, \ldots, 0 \rangle$ with $\mu_{k-i+1} \neq 0$.  Thus $s \mu = s_{k-i+1} \mu_{k-i+1} \theta^{m-1} \neq 0$, so that $\size{B_{i-1} \cap \ann(s)} = 0$.
For $i \leq j \leq k-1$, use the definition of $\mathcal{V}_i$ to see that $\mathcal{V}_j \subseteq \mathcal{V}_i = \ann(V_i) \subseteq \ann(s)$.  So $B_j \subseteq \ann(s)$ and $\size{B_j \cap \ann(s)} = \size{B_j} = q^{k-j-1}$.

Because $s \neq 0$, $\ann(s)$ is a vector subspace of $S_k^\sharp$ with $\dim_{\F_q} \ann(s) = k-1$.  By dimension counting, $\dim(\mathcal{V}_j \cap \ann(s))$ equals $k-j-1$ or $k-j$.  The case $\dim(\mathcal{V}_j \cap \ann(s)) = k-j$ occurs when $\mathcal{V}_j \subseteq \ann(s)$.  This is the case $i \leq j \leq k-1$ above.  When $j=i-1$, $\mathcal{V}_{i-1} \cap \ann(s) = \mathcal{V}_i$, so that $\size{B_{i-1} \cap \ann(s)} = 0$, as we saw above.  Finally, let $0 \leq j \leq i-2$.  Then $\dim(\mathcal{V}_j \cap \ann(s)) = k-j-1$ and $\dim(\mathcal{V}_{j+1} \cap \ann(s)) = k-j-2$.  This implies $\size{(\mathcal{V}_j - \mathcal{V}_{j+1}) \cap \ann(s)} = \size{\mathcal{V}_j  \cap \ann(s)} - \size{\mathcal{V}_{j+1}  \cap \ann(s)} = q^{k-j-2}(q-1)$.  Taking normalized representatives implies $\size{B_j \cap \ann(s)} = q^{k-j-2}$.
\end{proof}

\section{Two families of linear codes with the same \texorpdfstring{$\wwe$}{wwe}}  \label{sec:twoFamiliesChain}

We continue to assume that $R$ is a finite chain ring with maximal ideal $\m = (\theta)$ such that $R/\m \cong \F_q$ and $\theta^m = 0$, for some integer $m \geq 2$.  Let $w$ be an integer-valued weight on $R$ with maximal symmetry.  Denote the common value of $w$ on $\orb(\theta^i)$ by $w_i >0$, $i = 0, 1, \ldots, m-1$, and $w_m = w(0) = 0$.

To begin, we use $w_0, \ldots, w_{m-1}$ to define several numerical quantities.
For $i=0, 1, \ldots, m-1$, define
\begin{equation}  \label{eq:FormulaFora_i}
a_i = \sum_{j=0}^i q^j w_{m-j-1} .
\end{equation}
Also define, for $k=2, 3, \ldots, m$,
\begin{equation}  \label{eq:defnOfDelta_k}
\Delta_k = k q^{k-1} w_{m-1} - \sum_{i=0}^{k-1} q^{k-i-1} a_i .
\end{equation}

Recall the cyclic $R$-module $Z_k = R/(\theta^k)$ and the semisimple $R$-module $S_k = Z_1 \oplus \cdots \oplus Z_1$ ($k$ summands) from Section~\ref{sec:ChainRings}.  We will construct $R$-linear codes with $Z_k$ and $S_k$ as their underlying information modules.  The linear codes will be images of homomorphisms $\La: Z_k \ra R^n$ and $\Gamma: S_k \ra R^n$ of left $R$-modules.  As explained in Section~\ref{sec:W-stuff}, the linear codes are determined by their multiplicity functions.

\begin{defn}  \label{def:Ck}
Define a left $R$-linear code $C_k$ parametrized by $Z_k$ by using the linear functionals given by right multiplication by each of $\theta^{m-k}, \ldots, \theta^{m-1}$, each repeated $q^{k-1} w_{m-1}$ times, and the zero functional, repeated $\max\{0, -\Delta_k\}$ times; cf., \eqref{eq:defnOfDelta_k}.
\end{defn}
Equivalently, $C_k$ has a generator matrix of size $1 \times (k q^{k-1} w_{m-1} + \max\{0, -\Delta_k\})$, with entries $\theta^{m-k}, \ldots, \theta^{m-1}$, each repeated $q^{k-1} w_{m-1}$ times, plus entries of $0$, repeated $\max\{0, -\Delta_k\}$ times.  Thus, $0$ does not appear if $\Delta \geq 0$, and $0$ appears $-\Delta_k$ times when $\Delta_k < 0$.

\begin{prop}
The linear code $C_k$ of Definition~\textup{\ref{def:Ck}} has length $k q^{k-1} w_{m-1} + \max\{0, -\Delta_k\}$.  Its weight function $W_{\La}$ has values 
\[  W_{\La}(\theta^i) = q^{k-1} w_{m-1} (w_{m-k+i} + \cdots + w_{m-1}) , \]
for $i=0, 1, \ldots, k-1$, and the $w$-weight enumerator of $C_k$ is
\begin{align*}
\wwe_{C_k} &= 1 + \sum_{i=0}^{k-1} \size{\orb(\theta^i)} t^{W_{\La}(\theta^i)} \\
&= 1 + \sum_{i=0}^{k-1} q^{k-i-1}(q-1) t^{q^{k-1} w_{m-1} \sum_{j=i}^{k-1} w_{m-k+j}} .
\end{align*}  
\end{prop}

\begin{proof}
The formula for the length follows directly from Definition~\ref{def:Ck}.  In calculating $W_{\La}(\theta^i)$, remember that $w_m = w(\theta^m) = w(0)=0$.
\begin{align*}
W_{\La}(\theta^i) &= \sum_{j=m-k}^{m-1} q^{k-1} w_{m-1} w(\theta^i \theta^j) \\
&= q^{k-1} w_{m-1} \left( w_{i+m-k} + \cdots + w_{m-1} \right) .
\end{align*}
The functional in $Z_k^\sharp$ given by right multiplication by $\theta^{m-k}$ is injective, so that $\La$ is also injective.  Then $\wwe_{C_k}$ follows from \eqref{wwe-in-terms-of-orbits}.
\end{proof}

We now want to define linear codes $D_k$ parametrized by $\Gamma: S_k \ra R^n$ such that $\wwe_{C_k} = \wwe_{D_k}$.  
The form of $\wwe_{C_k}$ was determined by \eqref{wwe-in-terms-of-orbits}, in particular by the sizes of the orbits $\orb(\theta^i)$ and the value of $W_{\La}$ on those orbits.  
In order to be able to match terms in the equation $\wwe_{C_k} = \wwe_{D_k}$, 
we make use of the subsets $\mathcal{S}_1, \ldots, \mathcal{S}_k$ of $S_k$ defined following \eqref{eq:FilterSk}.  Also recall the sets $B_i$ defined in \eqref{defn:Bi}.  We will design $\Gamma$ so that  
$W_{\Gamma}$ is constant on each $\mathcal{S}_i$, with value equal to the value of $W_{\La}$ on $\orb(\theta^{k-i})$.

\begin{defn}   \label{def:Dk}
Define a left $R$-linear code $D_k$ parametrized by $\Gamma: S_k \ra R^n$ by using the linear functionals in $\cup_i B_i$, with each $\mu \in B_i$ repeated $a_i$ times (cf., \eqref{eq:FormulaFora_i}), $i=0,1,\ldots, k-1$, and the zero-functional in $B_k$ repeated $\max\{\Delta_k, 0\}$ times (cf., \eqref{eq:defnOfDelta_k}).
\end{defn}
Equivalently, $D_k$ has a generator matrix of size $k \times (\sum_{i=0}^{k-1} q^{k-i-1} a_i + \max\{\Delta_k, 0\})$, with columns given by $\mu \in \cup_i B_i$, with each $\mu \in B_i$ repeated $a_i$ times, $i=0,1, \ldots, k$, the zero-column in $B_k$ repeated $\max\{\Delta_k, 0\}$ times.  Thus, zero-columns do not appear if $\Delta \leq 0$, and the zero-column appears $\Delta_k$ times when $\Delta_k > 0$.

We now express the values of $W_{\Gamma}$ in terms of $a_0, a_1, \ldots, a_{k-1}$.
\begin{prop}
The weight function $W_{\Gamma}: S_k \ra \Z$ is constant on each $\mathcal{S}_i$, with  $W_{\Gamma}(\mathcal{S}_0) = 0$ and, for $i=1, \ldots, k$,
\[  W_{\Gamma}(\mathcal{S}_i) =  \left(  q^{k-i} a_{i-1} + \sum_{j=0}^{i-2} q^{k-j-2}(q-1) a_j  \right) w_{m-1} .  \]
\end{prop}

\begin{proof}
If $i=0$, then $\mathcal{S}_0 = \{0\}$, so that $W_{\Gamma}(\mathcal{S}_0) = 0$.  Now let $i=1, \ldots, k$, and let $s \in \mathcal{S}_i$.  Using Lemma~\ref{lem:IntersectionCount}, we see that 
\begin{align*}
W_{\Gamma}(s) &= \sum_{j=0}^{k-1} a_j \sum_{\mu \in B_j} w(s \mu) 
= \sum_{j=0}^{i-1} a_j \sum_{\mu \in B_j} w(s \mu) \\
&= \left(  q^{k-i} a_{i-1} + \sum_{j=0}^{i-2} q^{k-j-2}(q-1) a_j \right) w_{m-1} .
\end{align*}
This formula depends only on $i$, so $W_{\Gamma}$ is constant on each $\mathcal{S}_i$.
\end{proof}

In order that $\wwe_{C_k} = \wwe_{D_k}$, we need $W_{\La}(\orb(\theta^{k-i})) = W_{\Gamma}(\mathcal{S}_i)$, for $i=1, \ldots, k$.  That is, canceling a common factor, we need 
\begin{equation}  \label{eq:equateWs}
q^{k-1} (w_{m-i} + \cdots + w_{m-1}) = q^{k-i} a_{i-1} + \sum_{j=0}^{i-2} q^{k-j-2}(q-1) a_j   ,
\end{equation}
for $i=1, \ldots, k$.  The number of additive terms on each side of the equation is $i$.  Starting with $i=1$ and working upwards, we solve a triangular system recursively for $a_0, a_1, \ldots, a_{k-1}$.

\begin{lem}  \label{lem:solutionsToEqn}
The solutions of \eqref{eq:equateWs} are
\[  a_i = w_{m-1} + q w_{m-2} + \cdots + q^i w_{m-i-1} , \]
for $i=0, 1, \ldots, k-1$.  This formula matches \eqref{eq:FormulaFora_i}.
\end{lem}

\begin{proof}
Exercise, by induction.  The terms simplify by telescoping.
\end{proof}

\begin{thm}  \label{thm:SameWWE}
For each $k = 2, 3, \ldots, m$,  the codes $C_k$ and $D_k$ have the same length and satisfy $\wwe_{C_k} = \wwe_{D_k}$.
\end{thm}

\begin{proof}
The definition of $\Delta_k$ in \eqref{eq:defnOfDelta_k} guarantees that the codes have the same length.  The $a_i$ of \eqref{eq:FormulaFora_i} were defined so that Lemma~\ref{lem:solutionsToEqn} holds.  Thus $W_{\La}(\orb(\theta^{k-i})) = W_{\Gamma}(\mathcal{S}_i)$, for $i=1, \ldots, k$.  Because $\size{\orb(\theta^{k-i})} = \size{\mathcal{S}_i}$, the equality of the $w$-weight enumerators follows from \eqref{wwe-in-terms-of-orbits}.
\end{proof}

\begin{rem}
It is possible to generalize the constructions of $C_k$ and $D_k$ by allowing more general expressions for the multiplicities of the linear functionals appearing in $C_k$.  In Definition~\ref{def:Ck}, $\theta^{m-k}, \ldots, \theta^{m-1}$ could be repeated $q^{k-1} w_{m-1} b_{m-k}, \ldots, q^{k-1} w_{m-1} b_{m-1}$ times, respectively.  One can then mimic \eqref{eq:equateWs} and Lemma~\ref{lem:solutionsToEqn} to solve for the multiplicities $a_i$ used in defining $D_k$.  One must be careful in choosing the $b$'s in order that the $a$'s come out nonnegative.  Sufficiently large $b$'s should work.  The factors of $q^{k-1} w_{m-1}$ are present so that the $a$'s are integers.  The present work does not need this level of generality.
\end{rem}

In later sections, we will use the formula below for $\Delta_k$, which is expressed in terms of some numerical quantities defined next.
\begin{defn}  \label{defn:epsilons}
For $i=1, 2, \ldots, m$, define $\epsilon_i = w_i - w_{i-1}$.  For example, $\epsilon_1 = w_1 - w_0$, while $\epsilon_m= - w_{m-1}$ (as $w_m=0$).  For $i=1, \ldots, m-2$, define $\epsilon_i' = \epsilon_i$, and define $\epsilon_{m-1}' = (q-1)w_{m-1} - qw_{m-2} = q \epsilon_{m-1} + \epsilon_m$.
\end{defn} 
Also define these polynomial expressions in $q$:  $p_0 = 0$, and 
\begin{equation}  \label{eq:DefnOfpl}
p_i = 1 + 2q + 3q^2 + \cdots + i q^{i-1}, \quad i=1,2, \ldots . 
\end{equation}
Then for $i = 1, 2, \ldots $, the $p_i$ are positive and satisfy the formulas
\begin{equation}  \label{eq:pkFormula}
p_i - q p_{i -1} = 1 + q + q^2 + \cdots + q^{i -1} .
\end{equation}

\begin{prop}  \label{prop:FormulaForDelta}
For $k=2, 3, \ldots, m$, 
\[  \Delta_k = p_{k-1} \epsilon_{m-1}' + \sum_{j = 2}^{k-1} q^j p_{k-j} \epsilon_{m-j}' . \]
\end{prop}

\begin{proof}
Start with \eqref{eq:defnOfDelta_k}, and replace $a_i$ using \eqref{eq:FormulaFora_i}:  
\[ \Delta_k = k q^{k-1} w_{m-1} - \sum_{i=0}^{k-1} q^{k-i-1} \sum_{j=0}^i q^j w_{m-1-j} . \]
Interchange the order of summation and use \eqref{eq:pkFormula}:
\begin{align*}
\Delta_k &= k q^{k-1} w_{m-1} - \sum_{j=0}^{k-1}  \sum_{i=j}^{k-1} q^{k-i+j-1}  w_{m-1-j} \\
&= k q^{k-1} w_{m-1} - \sum_{j=0}^{k-1}  q^j (p_{k-j} - q p_{k-j-1})  w_{m-1-j} \\
&= k q^{k-1} w_{m-1} - \sum_{j=0}^{k-1}  q^j p_{k-j} w_{m-1-j} + \sum_{j=0}^{k-1}  q^{j+1} p_{k-j-1} w_{m-1-j} .
\end{align*}
Re-index the last summation with $\ell = j+1$ (the $\ell=k$ term vanishes), separate some initial terms, and combine the rest using Definition~\ref{defn:epsilons}:
\begin{align*}
\Delta_k 
&= k q^{k-1} w_{m-1} - \sum_{j=0}^{k-1}  q^j p_{k-j} w_{m-1-j} + \sum_{\ell=1}^k  q^{\ell} p_{k-\ell} w_{m-\ell} \\
&= k q^{k-1} w_{m-1} - p_k w_{m-1} + q p_{k-1} \epsilon_{m-1} + \sum_{j=2}^{k-1} q^j p_{k-j} \epsilon_{m-j}.
\end{align*}
Simplify the coefficient of $w_{m-1}$ and again use Definition~\ref{defn:epsilons}:
\begin{align*}
\Delta_k 
&=  p_{k-1} \epsilon_m + q p_{k-1} \epsilon_{m-1} + \sum_{j=2}^{k-1} q^j p_{k-j} \epsilon_{m-j} \\
&= p_{k-1} \epsilon_{m-1}' + \sum_{j=2}^{k-1} q^j p_{k-j} \epsilon_{m-j}' . \qedhere
\end{align*}
\end{proof}

\section{Analysis of dual codewords of low weight}
Our ultimate objective is to prove that for some $k = 2, 3, \ldots, m$, the codes $C_k$ and $D_k$ of Theorem~\ref{thm:SameWWE} have dual codes with different weight enumerators: $\wwe_{C_k^\perp} \neq \wwe_{D_k^\perp}$.  We will try to do this in the most direct way---by showing that $C_k^\perp$ and $D_k^\perp$ have different numbers of codewords of the smallest possible weight.  With that in mind, let's develop some notation.

Recall that we are assuming the chain ring $R$ is equipped with an integer-valued weight $w$ of maximal symmetry.  The common value on $\orb(\theta^i)$ is denoted $w_i >0$, and $w_m = w(0)=0$.  Let $\mathring{w} = \min\{w_i: i=0,1, \ldots, m-1\}$, so that $\mathring{w} >0$.  Define $\mathring{I} = \{i: w_i = \mathring{w}\}$, the set of exponents of $\theta$ that achieve the minimum value of the weight; $\mathring{I}$ is nonempty.

We now turn our attention to the linear codes $C_k$ and $D_k$ of Theorem~\ref{thm:SameWWE} and codewords of weight $d < 2\mathring{w}$ in their dual codes.  All such codewords must be singletons by Lemma~\ref{lem:singletons}.  We will abuse notation slightly by using the phrase `singleton $\theta^i$' to mean a singleton whose nonzero entry is a unit multiple of $\theta^i$.

\begin{lem}  \label{lem:singletonContribution}
Suppose an integer $d$ satisfies $\mathring{w} \leq d < 2\mathring{w}$.  Let $I_d = \{ i: w_i = d \}$.  If $0 \in I_d$, then, for $k = 2, 3, \ldots, m$,
\[  A_d(C_k^\perp) - A_d(D_k^\perp) = -\size{\orb(1)} \Delta_k - \sum_{\substack{ i \in I_d\\ 0<i<k}} (k-i) q^{k-1}w_{m-1} \size{\orb(\theta^i)} . \]
If $0 \not\in I_d$ and $I_d$ is nonempty, then, for $k = 2, 3, \ldots, m$,
\[  A_d(C_k^\perp) - A_d(D_k^\perp) =  - \sum_{\substack{ i \in I_d \\ i<k}} (k-i) q^{k-1}w_{m-1} \size{\orb(\theta^i)} . \]
\end{lem}

\begin{proof}
Nonzero contributions to $A_d(C_k^\perp) - A_d(D_k^\perp)$ are made by singletons of weight $d$ in $C_k^\perp$ or $D_k^\perp$.  The nonzero entry of a singleton of weight $d$ must be a unit multiple of $\theta^i$ with $i \in I_d$.  In order for a singleton to belong to a dual code, its nonzero entry---located in position $j$, say---must annihilate column $j$ of the generator matrix of the primal code.  

If $0 \in I_d$, a singleton $1$ annihilates only zero-columns.  (Remember that the nonzero entry of a singleton $1$ is a unit.)  The number of zero-columns is determined by $\Delta_k$: if $\Delta_k >0$, then $D_k$ has $\Delta_k$ zero-columns; if $\Delta_k < 0$, then $C_k$ has $-\Delta_k$ zero-colums.  The net contribution to $A_d(C_k^\perp) - A_d(D_k^\perp)$ is $-\size{\orb(1)} \Delta_k$, the number of units times the number of zero-columns.

If $0<i \in I_d$, then a singleton $\theta^i$ annihilates all the columns of $D_k$.  Such singletons contribute $- \size{\orb(\theta^i)} \cdot \length(D_k)$ to $A_d(C_k^\perp) - A_d(D_k^\perp)$.  
On the other hand, when $i \geq k$, a singleton $\theta^i$ annihilates all columns of $C_k$.  Such singletons contribute $\size{\orb(\theta^i)} \cdot \length(C_k)$ to $A_d(C_k^\perp) - A_d(D_k^\perp)$. 
When $i < k$, a singleton $\theta^i$ annihilates all columns of $C_k$ except those with entries $\theta^{m-k}, \ldots, \theta^{m-i-1}$.  Such singletons contribute $\size{\orb(\theta^i)} (\length(C_k) - (k-i) q^{k-1} w_{m-1})$ to $A_d(C_k^\perp) - A_d(D_k^\perp)$.  Because $C_k$ and $D_k$ have the same length, the total contribution by singleton $\theta^i$'s is $0$ when $i \geq k$ and $- (k-i) q^{k-1} w_{m-1} \size{\orb(\theta^i)}$ when $i<k$.
Summing over $i \in I_d$ completes the proof.
\end{proof}

Our main interest is dual codewords of weight $\mathring{w}$.  For $k=2, 3, \ldots, m$, define
\begin{equation}  \label{eq:DefnLittledelta_k}
\delta_k = A_{\mathring{w}}(C_k^\perp) - A_{\mathring{w}}(D_k^\perp) .
\end{equation}
Our aim is to show, whenever possible, for a given weight $w$, that $\delta_k \neq 0$ for some $k$.  We restate Lemma~\ref{lem:singletonContribution} for the case where $d = \mathring{w}$.

\begin{lem}  \label{lem:formulaForLittledelta}
Fix $k = 2, 3, \ldots, m$.  
If $0 \in \mathring{I}$, then
\[  \delta_k = -\size{\orb(1)} \Delta_k - \sum_{\substack{ i \in \mathring{I}\\ 0<i<k}} (k-i) q^{k-1}w_{m-1} \size{\orb(\theta^i)} . \]
If $0 \not\in \mathring{I}$, then 
\[  \delta_k =  - \sum_{\substack{ i \in \mathring{I} \\ i<k}} (k-i) q^{k-1}w_{m-1} \size{\orb(\theta^i)} . \]
\end{lem}

We draw three corollaries, using notation from Definition~\ref{defn:epsilons}.

\begin{cor}  \label{cor:w0big}
If $0 \not\in \mathring{I}$, then $\delta_k <0$ for all 
$k=1+\min{\mathring{I}}, \ldots, m$.  
Weights $w$ on $R$ with $w_0 > \mathring{w}$ do not respect duality.
\end{cor}

\begin{cor}  \label{cor:Only0}
Suppose $\mathring{I} = \{0\}$.  Then the following hold:
\begin{enumerate}
\item If $m \geq 3$ and $j$, $1 \leq j \leq m-1$, is the largest index with $\epsilon_j' \neq 0$, then $\delta_{m-j+1} \neq 0$.
\item  If $m=2$ and $\epsilon_1' = (q-1) w_1 - q w_0 \neq 0$, then $\delta_2 \neq 0$.
\item  If $m=2$, $\epsilon_1' = (q-1) w_1 - q w_0 = 0$, and $q>2$, then $A_{w_1}(C_2^\perp) < A_{w_1}(D_2^\perp)$.
\end{enumerate}
Weights $w$ on $R$ with $w_0 < w_i$ for all $i=1, 2, \ldots, m-1$ do not respect duality, except when $m=2$, $q=2$, and $w_1 = 2 w_0$; cf., Theorem~\textup{\ref{thm:chainMacWids}}.
\end{cor}

\begin{proof}
Suppose $m \geq 3$.  By the hypothesis on $\mathring{I}$, $\epsilon_1' = w_1 -w_0 >0$, so there exists a maximal index $j$, $1 \leq j \leq m-1$, with $\epsilon_j' \neq 0$.  Thus, $\epsilon_{\ell}' =0$ for $\ell = j+1, \ldots, m-1$.  From Proposition~\ref{prop:FormulaForDelta}, $\Delta_{m-j+1} = q^{m-j} \epsilon_j' \neq 0$.   From Lemma~\ref{lem:formulaForLittledelta} and the hypothesis, $\delta_{m-j+1} = - \size{\orb(1)} \Delta_{m-j+1} \neq 0$. 

When $m=2$, $\epsilon_1' = (q-1) w_1 - q w_0$, not $w_1-w_0$.  If $\epsilon_1' \neq 0$, the proof proceeds as above: $\Delta_2 \neq 0$, and $\delta_2 \neq 0$.  

When $m=2$ and $\epsilon_1'=0$, then $\Delta_2=\delta_2 = 0$.  However, $\epsilon_1'=0$ means that $(q-1) w_1 = q w_0$, i.e., $w_1 = (q/(q-1)) w_0 > w_0$.  But note, for integers $q \geq 2$, that $q/(q-1) \leq 2$ with equality holding if and only if $q=2$.  Assuming $q>2$, we have $\mathring{w} = w_0 < w_1 <2 \mathring{w}$.  By Lemma~\ref{lem:singletonContribution} applied to $d = w_1$ and $0 \not\in I_d$, we see that $A_{w_1}(C_2^\perp) < A_{w_1}(D_2^\perp)$.

When $m=2$, $q=2$, and $w_1 = 2 w_0$, we are in the situation of Theorem~\ref{thm:chainMacWids}, where the MacWilliams identities hold.
\end{proof}

\begin{cor}  \label{cor:chainm2}
Let $R$ be a finite chain ring with $m=2$.  Then every weight $w$ on $R$ having maximal symmetry does not respect duality, except for multiples of the Hamming weight \textup{(}any $q$\textup{)} or the homogeneous weight \textup{(}$q=2$ only\textup{)}.
\end{cor}

\begin{proof}
Because $m=2$, there are only $w_0$ and $w_1$.  If $w_0>w_1$, Corollary~\ref{cor:w0big} implies $w$ does not respect duality.  If $w_0=w_1$, $w$ is a multiple of the Hamming weight, and the MacWilliams identities hold \cite[Theorem~8.3]{wood:duality}.  If $w_0<w_1$, then Corollary~\ref{cor:Only0} applies: $w$ does not respect duality, except for multiples of the homogeneous weight if $q=2$.
\end{proof}

\section{Weak monotonicity and final arguments}
When $\{0\} \subsetneq \mathring{I}$, the formula for $\delta_k$ in Lemma~\ref{lem:formulaForLittledelta} is difficult to exploit systematically; the combinatorics can be formidable.  (But not always: see Example~\ref{ex:outlier}.)  In order to make progress, we will assume that the weight $w$ on the chain ring $R$ is \emph{weakly monotone}; i.e., we assume
\begin{equation}  \label{eq:WeaklyMonotone}
  \mathring{w} = w_0 \leq w_1 \leq \cdots \leq w_{m-2} \leq w_{m-1} . 
\end{equation}
This hypothesis implies that $\epsilon_i \geq 0$ for $i=1, 2, \ldots, m-1$ in Definition~\ref{defn:epsilons}.  However, $\epsilon_m = -w_{m-1} <0$.
The weakly monotone hypothesis allows us to state our main result.

\begin{thm}  \label{thm:MainThmWeaklyMonotone}
Let $R$ be a finite chain ring with a weakly monotone weight $w$.  Then 
$w$ does not respect duality, except when
\begin{itemize}
\item $w$ is a multiple of the Hamming weight, or
\item $m=2$, $q=2$, and $w_1=2w_0$.
\end{itemize}
\end{thm}

Theorem~\ref{thm:MainThmWeaklyMonotone} will follow from Theorem~\ref{thm:detailedMainThm}, which we will prove after we establish some technical lemmas.

Equalities are possible in \eqref{eq:WeaklyMonotone}.  Given a weakly monotone weight $w$, define $j_0$ to be the largest index such that 
\begin{equation}  \label{eq:Firstj}
\mathring{w} = w_0 = \cdots = w_{j_0} < w_{j_0+1}.
\end{equation}
Similarly, define $j_1$ to be the smallest index such that
\begin{equation}  \label{eq:Secondj}
w_{j_1 - 1} < w_{j_1} = \cdots = w_{m-1} .
\end{equation}
There are three situations to highlight, depending upon how many nonzero values $w$ takes.
\begin{itemize}
\item If $w$ has only one nonzero value, then $w_0 = \cdots = w_{m-1}$, so that $j_0=m-1$ and $j_1=0$.  The weight $w$ is a multiple of the Hamming weight.
\item  If $w$ has exactly two nonzero values, then $w_0 = \cdots =w_{j_0} < w_{j_0+1} = \cdots = w_{m-1}$, so that $j_1 = j_0+1$.  The homogeneous weight is an example of this, with $j_0 = m-2$, $j_1 = m-1$.
\item  If $w$ has three or more values, define $j_2$ so that
\[  \cdots w_{j_2-1} < w_{j_2} = \cdots = w_{j_1-1} < w_{j_1} = \cdots = w_{m-1} . \]
Said another way, $j_1$ is the largest index less than $m$ with $\epsilon_{j_1} > 0$, and $j_2$ is the second-largest index less than $m$ with $\epsilon_{j_2} > 0$.
\end{itemize}

The weight $w$ has a least two nonzero values if and only if $j_1>0$.  In that case, $j_0 < j_1 \leq m-1$, with
\begin{equation}  \label{eq:vanishingEpsilons}
\epsilon_{j_1} >0 \text{ and } \epsilon_{j_1+1} = \cdots = \epsilon_{m-1} = 0 .
\end{equation}
If $w$ has at least three values, then $j_0 < j_2 < j_1 \leq m-1$ and, in addition to \eqref{eq:vanishingEpsilons}, we have
\begin{equation}  \label{eq:vanishingEpsilonsAgain}
\epsilon_{j_2} >0 \text{ and } \epsilon_{j_2+1} = \cdots = \epsilon_{j_1-1} = 0 .
\end{equation}

The key to our analysis is a simplified expression for $\delta_k$; cf., \eqref{eq:DefnOfpl}.

\begin{lem}  \label{lem:deltak-in-terms-of-epsilons}
Suppose $\{0,1, \ldots, j_0\} = \mathring{I}$, with $j_0 \geq 1$.  If $k$ is an integer, $2 \leq k \leq j_0+1$, then
\[  \delta_k = - q^m (q-1) \sum_{j=1}^{k-1} q^{j-1} p_{k-j} \epsilon_{m-j} . \]
If $k$ is an integer, $j_0+2 \leq k \leq m$, then
\begin{align*}
\delta_k = - ((k-j_0-1)& q^{m+k-j_0-2} -q^{m+k-j_0-3} - \cdots -q^{m-1}) \epsilon_m  \\
- q^m &(q-1) \sum_{j=1}^{k-1} q^{j-1} p_{k-j} \epsilon_{m-j}  .
\end{align*}
\end{lem}

\begin{proof}  For any integer $k$, $2 \leq k \leq m$, Lemma~\ref{lem:formulaForLittledelta} and  \eqref{eq:OrbitSizes} imply 
\begin{align*}
\delta_k &= -\size{\orb(1)} \Delta_k - \sum_{i=1}^{\min\{j_0,k-1\}} (k-i) q^{k-1}w_{m-1} \size{\orb(\theta^i)} \\
&= - q^{m-1}(q-1) \Delta_k - \sum_{i=1}^{\min\{j_0,k-1\}} (k -i) q^{m+k -i-2} (q-1) w_{m-1} .
\end{align*}
In the formula for $\Delta_k$ given in Proposition~\ref{prop:FormulaForDelta}, use $\epsilon_{m-1}' = q \epsilon_{m-1} + \epsilon_m = q \epsilon_{m-1} - w_{m-1}$.  Because of telescoping sums, the $w_{m-1}$-terms cancel completely when $2 \leq k \leq j_0+1$ or cancel partially when $j_0+2 \leq k \leq m$.  The terms that remain are as stated.
\end{proof}

In the case of $j_0+2 \leq k \leq m$, note that the $\epsilon_m$-term is positive, as $-\epsilon_m = w_{m-1}>0$ and the numerical sum is positive.  The other terms are nonpositive.  The balance between the terms appears to be problematic.  When $2 \leq k \leq j_0+1$, $\delta_k \leq 0$.

We will also need information about how $\delta_k$ changes when $k \geq j_0+2$.
\begin{lem}  \label{lem:changeIndelta}
If $k \geq j_0+2$ and $k \geq m-j_1$, then
\begin{align*}
\delta_{k+1} - \delta_k = -&q^{m+k-j_0-2}(q-1) \\
& \times \left\{ (k-j_0) \epsilon_m + q^{j_0+1} \sum_{i=m-j_1}^k (k-i+1) \epsilon_{m-i} \right\} .
\end{align*}
\end{lem}

\begin{proof}
Because $k \geq j_0+2$, the second formula in Lemma~\ref{lem:deltak-in-terms-of-epsilons} applies to both $\delta$'s.  All but the highest order terms cancel, leaving
\begin{align*}
\delta_{k+1} - \delta_k = -&(k-j_0) q^{m+k-j_0-2}(q-1) \epsilon_m \\
& -q^{m+k-1}(q-1) \sum_{i=1}^k (k-i+1) \epsilon_{m-i} .
\end{align*}
Because $j_1 \geq m-k$, \eqref{eq:vanishingEpsilons} implies that terms vanish in the summation for $i=1, 2, \ldots, m-j_1-1$.
\end{proof}

\begin{cor}  \label{cor:specialCaseChangedelta}
Suppose $j_0 + j_1 = m-1$.  If $k \geq j_0+2$, then
\begin{align*}
\delta_{k+1} - \delta_k = -&q^{m+k-j_0-2}(q-1) (k-j_0) \left( \epsilon_m +q^{j_0+1} \epsilon_{j_1} \right) \\
& -q^{m+k-1}(q-1) \sum_{i=m-j_1+1}^k (k-i+1) \epsilon_{m-i} .
\end{align*}
\end{cor}

\begin{proof}
Note first that $j_0 + j_1 = m-1$ implies $j_0 + 2 = m-j_1+1$.  Thus, if $k \geq j_0+2$, then $k \geq m-j_1$ is automatic.  Apply Lemma~\ref{lem:changeIndelta} and notice that the term inside the summation with $i=m-j_1$ is $(k-m+j_1 +1) \epsilon_{j_1} = (k-j_0) \epsilon_{j_1}$.
\end{proof}

The next theorem gives a more detailed description of the claims in Theorem~\ref{thm:MainThmWeaklyMonotone}.

\begin{thm} \label{thm:detailedMainThm}
Let $R$ be a finite chain ring with a weakly monotone weight $w$.  Let $j_0$ and $j_1$ be as defined in \eqref{eq:Firstj} and \eqref{eq:Secondj}.  Then the following statements hold.
\begin{enumerate}
\item  If $j_0=m-1$, then $w$ is a multiple of the Hamming weight.
\item If $j_0=0$, then Corollary~\textup{\ref{cor:Only0}} applies.
\end{enumerate}
In the following statements, assume $1 \leq j_0 < m-1$.
\begin{enumerate} \addtocounter{enumi}{2}
\item  If $j_0 + j_1 \geq m$, then $\delta_k < 0$ for $m- j_1 + 1 \leq k \leq j_0 +1$.
\item  If $j_0 + j_1 \leq m -2$, then $\delta_k > 0$ for $j_0 + 2 \leq k \leq m-j_1$.
\item  If $j_0 +j_1 = m-1$ and $w_{m-1} \neq q^{j_0+1} \epsilon_{j_1}$, then $\delta_{j_0+2} \neq 0$.
\item  Suppose $j_0 + j_1 = m-1$ and $w_{m-1} = q^{j_0+1} \epsilon_{j_1}$.  If $w$ has at least three nonzero values, then $\delta_k <0$ for  $k \geq m-j_2+1$.
\item  Suppose $j_0 + j_1 = m-1$, $w_{m-1} = q^{j_0+1} \epsilon_{j_1}$, and $w$ has two nonzero values.  Then $j_1 = j_0 + 1$, $m = 2j_0 + 2$ is even, and 
\[  A_{w_{m-1}}(C_k^\perp) - A_{w_{m-1}}(D_k^\perp) < 0  \]
for $k > j_1$.
\end{enumerate}
\end{thm}

\begin{rem}
The length of the `run' $w_{j_1} = \cdots = w_{m-1}$ is $m-j_1$.  The length of the `run' $w_0 = \cdots = w_{j_0}$ is $j_0+ 1$.  Their difference is
\[  (m-j_1) - (j_0+1) = \begin{cases}
+ , & j_0 + j_1 \leq m-2, \\
0 , & j_0 + j_1 = m-1, \\
- , & j_0 + j_1 \geq m.
\end{cases}  \]
\end{rem}

\begin{proof}[Proof of Theorem~\textup{\ref{thm:detailedMainThm}}]
The first two claims were explained after \eqref{eq:Secondj}.
If $1 \leq j_0 < m-1$ and $j_0 + j_1 \geq m$, then any $k$ satisfying $m- j_1 + 1 \leq k \leq j_0 +1$ has $\delta_k <0$ by Lemma~\ref{lem:deltak-in-terms-of-epsilons} and \eqref{eq:vanishingEpsilons}.
Namely, the first formula in Lemma~\ref{lem:deltak-in-terms-of-epsilons} applies, and $\epsilon_{j_1} >0$ appears in the formula because $j_1 \geq m-k+1$.

If $1 \leq j_0 < m-1$ and $j_0 + j_1 \leq m -2$, then any $k$ satisfying $j_0 + 2 \leq k \leq m-j_1$ has $\delta_k >0$.  Now the second formula in Lemma~\ref{lem:deltak-in-terms-of-epsilons} applies, and only the (positive) $\epsilon_m$-term survives, because $j_1 < m-k+1$.

If $1 \leq j_0 < m-1$ and $j_0 + j_1 = m-1$, then $j_0 + 2 = m - j_1 +1$.  Setting $k = j_0 + 2$, we see from the second formula of Lemma~\ref{lem:deltak-in-terms-of-epsilons} that 
\begin{align*}
\delta_{j_0 + 2} =& -(q^m - q^{m-1}) \epsilon_m - q^m (q-1) q^{j_0} \epsilon_{j_1} \\
&= -q^{m-1} (q-1) \left( \epsilon_m + q^{j_0+1} \epsilon_{j_1} \right) .
\end{align*} 
As $\epsilon_m = -w_{m-1}$, we see that $w_{m-1} \neq q^{j_0+1} \epsilon_{j_1}$ implies $\delta_{j_0+2} \neq 0$.

Next, suppose $1 \leq j_0 < m-1$, $j_0 + j_1 = m-1$, and $w_{m-1} = q^{j_0+1} \epsilon_{j_1}$.  We just saw that this implies $\delta_{j_0+2} = 0$.  Applying Corollary~\ref{cor:specialCaseChangedelta}, first with $k=j_0+2$, and then recursively, we see, for any $k \geq j_0+2$, that $\delta_k$ has the form:
\[  \delta_k = - \sum_{i=m-j_1+1}^{k-1} c(k)_i \epsilon_{m-i} , \]
where each $c(k)_i$ is a positive integer depending on $k$.  Because each $\epsilon_{m-i} \geq 0$, each $\delta_k \leq 0$, and $\delta_k < 0$ if at least one $\epsilon_{m-i} > 0$ in the interval of summation.  Remember \eqref{eq:vanishingEpsilonsAgain}.  If $k \geq m-j_2+1$, then $m-j_2 +1 > m-j_1+1 = j_0+2$, $m-(k-1) \leq j_2$, and $\epsilon_{j_2} > 0$ appears in the expression for $\delta_k$.  Thus $\delta_k < 0$. 

Finally, suppose $1 \leq j_0 < m-1$, $j_0 + j_1 = m-1$, $w_{m-1} = q^{j_0+1} \epsilon_{j_1}$, and $w$ has exactly two nonzero values.  That means that $j_1 = j_0 +1$, so that $m = 2j_0 + 2$.  In addition, we must have $\epsilon_{j_1} = w_{m-1} - w_0$.  Using this in the equation $w_{m-1} = q^{j_0+1} \epsilon_{j_1}$ yields $q^{j_0+1} w_0 = (q^{j_0+1} -1) w_{m-1}$.

Because the coefficients in this last equation are relatively prime, there exists a positive integer $s$ such that $w_0 = (q^{j_0+1}-1) s$ and $w_{m-1} = q^{j_0+1} s$.  Calculate: $2 w_0 - w_{m-1} = (q^{j_0+1} -2) s > 0$, because $j_0 \geq 1$ and $q \geq 2$.  (But see Remark~\ref{rem:curse}.)  Thus $\mathring{w}=w_0 < w_{m-1} < 2 \mathring{w}$.  Using the second formula of Lemma~\ref{lem:singletonContribution}, with $d = w_{m-1}$ and $0 \not\in I_d$, we see that $A_{w_{m-1}}(C_k^\perp) - A_{w_{m-1}}(D_k^\perp) < 0$ for any $k > j_1 = \min I_d$.  
\end{proof}

\begin{rem} \label{rem:curse}
In general, for positive integers $q \geq 2$, $q^{j_0+1} -2 \geq 0$, with equality if and only if $j_0=0$ and $q=2$.  Equality, again, points to Theorem~\ref{thm:chainMacWids}.
\end{rem}

\begin{ex}  \label{ex:m3exception}
When $m=3$, the only situation not covered by Corollaries~\ref{cor:w0big}, \ref{cor:Only0}, or Theorem~\ref{thm:MainThmWeaklyMonotone} is when $w_0=w_2 < w_1$.  In this case, 
$\epsilon_1' =  w_1-w_0 > 0$ and $\epsilon_2' = (q-1)w_2 - q w_1<0$.
Then Proposition~\ref{prop:FormulaForDelta} and Lemma~\ref{lem:formulaForLittledelta} imply that
\[ \delta_2 = -\size{\orb(1)} \Delta_2 = -q^2 (q-1) ((q-1)w_2 - q w_1) >0 . \]
Thus, this $w$ does not respect duality.  We conclude that, for $m=3$, the only weights that respect duality are multiples of the Hamming weight.
\end{ex}

\section{Symmetrized enumerators and examples}  \label{sec:examples}
In this section there are details of the MacWilliams identities for the symmetrized enumerator for the action of the full group $\U$ of units on a finite chain ring $R$, followed by several examples.
The details supplement the general outline provided in Appendix~\ref{sec:appendix}.

Suppose $R$ is a finite chain ring with $\m=(\theta)$, $R/\m \cong \F_q$, and $\theta^m=0$.  As we have seen earlier, the group of units $\U$ acts on $R$ on the left, with orbits $\orb(1), \orb(\theta), \ldots, \orb(\theta^{m-1}), \orb(\theta^m) = \{0\}$.  The dual action of $\U$ on $R^\sharp \cong R_R$ has the same orbit structure.  

For an element $r \in R$, define $\nu(r)$ via $r \in \orb(\theta^{\nu(r)})$; i.e., $\nu(r)$ is the exponent $i$ of $\theta$ such that $r = u \theta^i$ for some unit $u \in \U$.   
Define $f: R \ra \C[Z_0, \ldots, Z_m]$ by $f(r) = Z_{\nu(r)}$, $r \in R$.  Then define $F: R^n \ra \C[Z_0, \ldots, Z_m]$ by
\begin{equation}
F(x) = \prod_{\ell=1}^n f(x_\ell) =  \prod_{\ell=1}^n Z_{\nu(x_\ell)}, \quad x = (x_1, x_2, \ldots, x_n) \in R^n . 
\end{equation}
The \emph{symmetrized enumerator} of an additive code $C \subseteq R^n$ is the following element $\se_C \in \C[Z_0, \ldots, Z_m]$:
\[  \se_C = \se_C(Z_0, \ldots, Z_m) = \sum_{x \in C} F(x) = \sum_{x \in C} \prod_{\ell=1}^n Z_{\nu(x_\ell)} . \]

Any finite chain ring $R$ is a Frobenius ring, with a generating character $\chi$.  We will use the following properties of $\chi$ \cite{wood:duality}:
\begin{enumerate}
\item  for any nonzero ideal $I$ of $R$, $\sum_{r \in I} \chi(r) = 0$; \label{enum:vanishing}
\item  every $\pi \in \hr$ has the form $\pi = \chi^r$ for some unique $r \in R$;
\item $\chi^r(0)=1$ for all $r \in R$.
\end{enumerate}

We first calculate the sum of a character over the orbits $\orb(\theta^i)$.

\begin{lem}  \label{lem:chi-OrbitSums}
Suppose $i,j=0, 1, \ldots, m$ and $r \in \orb(\theta^j)$.  Then
\[  \sum_{s \in \orb(\theta^i)} \chi^r(s) = \begin{cases}
\hphantom{-}0, & i+j \leq m-2, \\
-q^{m-i-1} , & i+j=m-1, \\  
\hphantom{-}q^{m-i-1}(q-1) , & i+j \geq m, i<m, \\  
\hphantom{-} 1, & i+j \geq m, i=m.
\end{cases}  \]
\end{lem} 

\begin{proof}
When $i=m$, $\theta^m=0$. Then $\sum_{s \in \orb(\theta^m)} \chi^r(s) = 1$, as $\chi^r(0)=1$.  For $i=0, 1, \ldots, m-1$, $\orb(\theta^i) = (\theta^i) - (\theta^{i+1})$, so, for $\chi$ itself,
\begin{align}
\sum_{s \in \orb(\theta^i)} \chi(s) &= \sum_{s \in (\theta^i)} \chi(s) - \sum_{s \in (\theta^{i+1})} \chi(s) \notag \\
&= \begin{cases}
\hphantom{-}0, & i=0, 1, \ldots, m-2, \\
-1, & i=m-1,
\end{cases} \label{eq:char-sum}
\end{align} 
using property \eqref{enum:vanishing} of $\chi$.

Now suppose $r \in \orb(\theta^j)$, so that $r = u \theta^j$, $u \in \U$.  Left multiplication by $r$ maps $\orb(\theta^i)$ onto $\orb(\theta^{i+j})$, with each element in $\orb(\theta^{i+j})$ being hit $\size{\orb(\theta^i)}/\size{\orb(\theta^{i+j})}$ times.
(Because $\U \theta^i = \theta^i \U$, units can be moved across powers of $\theta$; for any unit $u$, $u \theta^i = \theta^i u'$ for some unit $u'$.)
This implies that
\[  \sum_{s \in \orb(\theta^i)} \chi^r(s) = \sum_{s \in \orb(\theta^i)} \chi(rs) = \frac{\size{\orb(\theta^i)}}{\size{\orb(\theta^{i+j})}} \sum_{t \in \orb(\theta^{i+j})} \chi(t) .\]
Using \eqref{eq:OrbitSizes} and \eqref{eq:char-sum}, we get the stated result.
\end{proof}

\begin{rem}
The formulas in Lemma~\ref{lem:chi-OrbitSums} depend only on the orbit of $r$, not $r$ itself.  This is a general feature of character sums over the blocks of a partition coming from a group action \cite[Theorem~2.6]{MR3336966}.
\end{rem}

Define the \emph{generalized Kravchuk matrix} $K$  by
\[  K_{ij} = \sum_{s \in \orb(\theta^i)} \chi^r(s), \quad r \in \orb(\theta^j) . \]

We calculate the Fourier transform of $f$ and $F$ as in \eqref{eqn:FTdefn}.

\begin{lem}  \label{lem:FTforChain}
For any $r \in R$, the Fourier transform of $f$ is
\[  \ft{f}(r) = \sum_{i=0}^m Z_i K_{ij}, \quad r \in \orb(\theta^j) . \]
The Fourier transform of $F$ is
\[  \ft{F}(r_1, r_2, \ldots, r_n) = \prod_{\ell=1}^n \ft{f}(r_\ell) = \prod_{\ell=1}^n \left(\sum_{i=0}^m Z_i K_{i,\nu(r_\ell)} \right). \]
\end{lem}

\begin{proof}
The identity $\ft{F}(r_1, r_2, \ldots, r_n) = \prod_{\ell=1}^n \ft{f}(r_\ell)$ is well-known, \cite[Proposition~A.5]{wood:duality}, so it is enough to
calculate $\ft{f}(r)$:
\begin{align*}
\ft{f}(r) &= \sum_{s \in R} \chi(rs) f(s) = \sum_{s \in R} \chi^r(s) Z_{\nu(s)} \\
&= \sum_{i=0}^m \sum_{s \in \orb(\theta^i)} \chi^r(s) Z_i = \sum_{i=0}^m Z_i \sum_{s \in \orb(\theta^i)} \chi^r(s) .  \qedhere
\end{align*}
\end{proof}

We now have the MacWilliams identities for the symmetrized enumerator over a finite chain ring; cf., \cite[Theorem~8.4]{wood:duality}, \cite[Theorem~3.5]{MR3336966}.

\begin{thm}  \label{thm:MWforSE}
Suppose $R$ is a finite chain ring and $C \subseteq R^n$ is a left $R$-linear code.  Then, using $C^\perp = \mathcal{R}(C)$, 
\[  \se_{C^\perp}(Z_0, \ldots, Z_m) = \left. \frac{1}{\size{C}}  \se_{C}(\mathcal{Z}_0, \ldots, \mathcal{Z}_m) \right|_{\mathcal{Z}_j = \sum_{i=0}^m Z_i K_{ij}} . \]
\end{thm}

\begin{proof}
Follow the outline in Appendix~\ref{sec:appendix}, apply Lemma~\ref{lem:FTforChain}, and note that $\mathcal{R}(C) = \mathfrak{R}(C)$ for left $R$-linear codes.
 \end{proof}
 
A version of this theorem, valid for the partition determined by the homogeneous weight, appears in \cite[Theorem~2.1]{MR4319450}

\begin{ex}
Let $R = \Z/8\Z$.  Then $\U = \{1,3,5,7\}$.  The $\U$-orbits are $\orb(1)= \U$, $\orb(2) = \{2,6\}$, $\orb(4) = \{4\}$, and $\orb(0)=\{0\}$.  The generalized Kravchuk matrix is
\[  K = \left[ \begin{array}{rrrr}
0 & 0 & -4 & 4 \\
0 & -2 & 2 & 2 \\
-1 & 1 & 1 & 1 \\
1 & 1 & 1 & 1
\end{array} \right] .  \]

Let $\wg$ be the homogeneous weight, so that $\wg_0=\wg_1=1$, $\wg_2=2$, and $\wg_3=0$.   The linear codes $C_3$ and $D_3$ of Theorem~\ref{thm:SameWWE} have the following codewords,
with multiplicities listed above the horizontal line, telling how many times the given entry is repeated.
\[  \begin{array}{rrr}
8 & 8 & 8 \\ \hline
0 & 0 & 0 \\
1 & 2 & 4 \\
2 & 4 & 0 \\
3 & 6 & 4 \\
4 & 0 & 0 \\
5 & 2 & 4 \\
6 & 4 & 0 \\
7 & 6 & 4 
\end{array}  \quad\quad
\begin{array}{rrrrrrr}
2 & 2 & 2 & 2 & 4 & 4 & 8 \\ \hline
0 & 0 & 0 & 0 & 0 & 0 & 0 \\
4 & 4 & 4 & 4 & 0 & 0 & 0 \\
0 & 0 & 4 & 4 & 4 & 4 & 0 \\
4 & 4 & 0 & 0 & 4 & 4 & 0 \\
0 & 4 & 0 & 4 & 0 & 4 & 4 \\
4 & 0 & 4 & 0 & 0 & 4 & 4 \\
0 & 4 & 4 & 0 & 4 & 0 & 4 \\
4 & 0 & 0 & 4 & 4 & 0 & 4
\end{array}  \]
Then the symmetrized enumerators of the codes are:
\begin{align*}
\se_{C_3} &= 4 Z_0^8 Z_1^8 Z_2^8 + 2 Z_1^8 Z_2^8 Z_3^8 + Z_2^8 Z_3^{16} + Z_3^{24}, \\
\se_{D_3} &= 4 Z_2^{16} Z_3^8 + 2 Z_2^{12} Z_3^{12} + Z_2^8 Z_3^{16} + Z_3^{24}.
\end{align*}
One can compute the symmetrized enumerators of the dual codes (via Theorem~\ref{thm:MWforSE} and  SageMath \cite{sagemath}, say), but the results have too many terms to include here.

Specializing $Z_i \leadsto t^{\wg_i}$ (and taking a Taylor expansion for the dual codes) yields the homogeneous weight enumerators:
\begin{align*}
\howe_{C_3} &= 1 + t^{16} + 2 t^{24} + 4 t^{32} , \\
\howe_{D_3} &= 1 + t^{16} + 2 t^{24} + 4 t^{32} , \\
\howe_{C_3^\perp} &= 1 + 16 t + 1848 t^2 + 60400 t^3 + \cdots , \\
\howe_{D_3^\perp} &= 1 + 48 t + 1832 t^2 + 64656 t^3 + \cdots .
\end{align*}
The computed value $\delta_3 =A_1(C_3^\perp) - A_1(D_3^\perp) = -32$ matches the value given in the second formula of Lemma~\ref{lem:deltak-in-terms-of-epsilons}.
\end{ex}

\begin{ex}
Still use $R=\Z/8\Z$, but change the weight to $w_0=1, w_1=w_2=2, w_3=0$.  The multiplicities of both codes (call them $C_{3'}$ and $D_{3'}$) change accordingly:
\[  \begin{array}{rrrr}
8 & 8 & 8 & 6 \\ \hline
0 & 0 & 0 & 0 \\
1 & 2 & 4 & 0 \\
2 & 4 & 0 & 0 \\
3 & 6 & 4 & 0 \\
4 & 0 & 0 & 0 \\
5 & 2 & 4 & 0 \\
6 & 4 & 0 & 0 \\
7 & 6 & 4 & 0 
\end{array}  \quad\quad
\begin{array}{rrrrrrr}
2 & 2 & 2 & 2 & 6 & 6 & 10 \\ \hline
0 & 0 & 0 & 0 & 0 & 0 & 0 \\
4 & 4 & 4 & 4 & 0 & 0 & 0 \\
0 & 0 & 4 & 4 & 4 & 4 & 0 \\
4 & 4 & 0 & 0 & 4 & 4 & 0 \\
0 & 4 & 0 & 4 & 0 & 4 & 4 \\
4 & 0 & 4 & 0 & 0 & 4 & 4 \\
0 & 4 & 4 & 0 & 4 & 0 & 4 \\
4 & 0 & 0 & 4 & 4 & 0 & 4
\end{array}  \]
Now the symmetrized enumerators are:
\begin{align*}
\se_{C_{3'}} &= 4 Z_0^8 Z_1^8 Z_2^8 Z_3^6 + 2 Z_1^8 Z_2^8 Z_3^{14} + Z_2^8 Z_3^{22} + Z_3^{30}, \\
\se_{D_{3'}} &= 4 Z_2^{20} Z_3^{10} + 2 Z_2^{16} Z_3^{14} + Z_2^8 Z_3^{22} + Z_3^{30}.
\end{align*}
The $w$-weight enumerators are
\begin{align*}
\wwe_{C_{3'}} &= 1 + t^{16} + 2 t^{32} + 4 t^{40} , \\
\wwe_{D_{3'}} &= 1 + t^{16} + 2 t^{32} + 4 t^{40} , \\
\wwe_{C_{3'}^\perp} &= 1 + 24 t + 1074 t^2 + 36584 t^3 + \cdots , \\
\wwe_{D_{3'}^\perp} &= 1  \phantom{ {}+ 48 t} + 1354 t^2 + 34304 t^3 + \cdots .
\end{align*}
For this $w$, $\mathring{I} = \{0\}$.  Lemma~\ref{lem:formulaForLittledelta} implies that $\delta_{3'} = -4 \Delta_{3'} = 24$, which matches the computed value.
\end{ex}

\begin{ex}  \label{ex:outlier}
For a final example, still use $R=\Z/8\Z$, but change the weight to $w_0=1, w_1=2, w_2=1, w_3=0$, which lies outside the scope of the main results given in previous sections; cf., Example~\ref{ex:m3exception}.  The multiplicities of both codes (call them $C_{3''}$ and $D_{3''}$) change accordingly:
\[  \begin{array}{rrrr}
4 & 4 & 4 & 11 \\ \hline
0 & 0 & 0 & 0 \\
1 & 2 & 4 & 0 \\
2 & 4 & 0 & 0 \\
3 & 6 & 4 & 0 \\
4 & 0 & 0 & 0 \\
5 & 2 & 4 & 0 \\
6 & 4 & 0 & 0 \\
7 & 6 & 4 & 0 
\end{array}  \quad\quad
\begin{array}{rrrrrrr}
1 & 1 & 1 & 1 & 5 & 5 & 9 \\ \hline
0 & 0 & 0 & 0 & 0 & 0 & 0 \\
4 & 4 & 4 & 4 & 0 & 0 & 0 \\
0 & 0 & 4 & 4 & 4 & 4 & 0 \\
4 & 4 & 0 & 0 & 4 & 4 & 0 \\
0 & 4 & 0 & 4 & 0 & 4 & 4 \\
4 & 0 & 4 & 0 & 0 & 4 & 4 \\
0 & 4 & 4 & 0 & 4 & 0 & 4 \\
4 & 0 & 0 & 4 & 4 & 0 & 4
\end{array}  \]
Now the symmetrized enumerators are:
\begin{align*}
\se_{C_{3''}} &= 4 Z_0^4 Z_1^4 Z_2^4 Z_3^{11} + 2 Z_1^4 Z_2^4 Z_3^{15} + Z_2^4 Z_3^{19} + Z_3^{23}, \\
\se_{D_{3''}} &= 4 Z_2^{16} Z_3^7 + 2 Z_2^{12} Z_3^{11} + Z_2^4 Z_3^{19} + Z_3^{23}.
\end{align*}
The $w$-weight enumerators are
\begin{align*}
\wwe_{C_{3''}} &= 1 + t^4 + 2 t^{12} + 4 t^{16} , \\
\wwe_{D_{3''}} &= 1 + t^4 + 2 t^{12} + 4 t^{16} , \\
\wwe_{C_{3''}^\perp} &= 1 + 63 t + 2111 t^2 + 51635 t^3 + \cdots , \\
\wwe_{D_{3''}^\perp} &= 1 + 23 t + 1195 t^2 + 38431 t^3 + \cdots .
\end{align*}
In this example $0 \in \mathring{I}$, and Lemma~\ref{lem:formulaForLittledelta} implies that $\delta_{3''} = 44-4=40$, which matches the computed value.
The weight $w$ does not satisfy the hypotheses of Corollaries~\ref{cor:w0big}, \ref{cor:Only0} or Theorem~\ref{thm:MainThmWeaklyMonotone}; nonetheless, we see that $w$ does not respect duality.  Example~\ref{ex:m3exception} uses $k=2$ to reach the same conclusion.
\end{ex}

\part{Matrix Rings over Finite Fields}

\section{Matrix modules, their orbits, and a positive result}
We begin our study of matrix rings by describing certain matrix modules, the orbits of the group of units, the homogeneous weight, and the MacWilliams identities for $M_{2 \times 2}(\F_2)$.

Fix integers $k, m$ with $2 \leq k \leq m$.  Let $R = M_{k \times k}(\F_q)$ be the ring of $k \times k$ matrices over a finite field $\F_q$, and let $M = M_{k \times m}(\F_q)$ be the left $R$-module of $k \times m$ matrices over $\F_q$.  Both $R$ and $M$ are vector spaces over $\F_q$.  The scalar multiplication of $R$ on $M$ is the multiplication of matrices.  
The group $\GL(k, \F_q)$ of invertible $k \times k$ matrices is the group of units $\U=\U(R)$ of $R$; $\U$ acts on $M$ on the left via matrix multiplication.  The size of $\GL(k, \F_q)$ is:
\begin{equation}  \label{eq:size-of-GLk}
\size{\GL(k, \F_q)} = (q^k-1) (q^k-q) \cdots (q^k-q^{k-1}) .
\end{equation}
Our first objective is to understand the cyclic left $R$-submodules of $M$ and the $\U$-orbits in $M$.

Given a $k \times m$ matrix $x \in M$, denote by $\rowsp(x)$ the row space of $x$, i.e., the $\F_q$-linear subspace of $\F_q^m$ spanned by the rows of $x$.
A left $R$-module is \emph{cyclic} if it is generated by one element; i.e., it has the form $Rx$ for some element $x$ in the module.  Denote the orbit of $x \in M$ under the action of $\U$ by $\orb(x)$ or $[x]$; denote the rank of a matrix $x$ by $\rk x$.

\begin{lem}  \label{lem:Rx-rowspace}
Let $x \in M = M_{k \times m}(\F_q)$.  Then
\begin{itemize}
\item for $y \in M$, $y \in Rx$ if and only if $\rowsp(y) \subseteq \rowsp(x)$;
\item  for $y \in M$, $Ry = Rx$ if and only if $\rowsp(y) = \rowsp(x)$ if and only if $\orb(y)=\orb(x)$;
\item  if $y \in \orb(x)$, then $\rk y = \rk x$.
\end{itemize}
\end{lem}

\begin{proof}
If $y = rx$, then the rows of $y$ are linear combinations of the rows of $x$.  This implies $\rowsp(y) \subseteq \rowsp(x)$.  Conversely, if $\rowsp(y) \subseteq \rowsp(x)$, then each row of $y$ is a linear combination of the rows of $x$, say $y_i = \sum_{j=1}^k r_{ij} x_j$, for some $r_{ij} \in \F_q$, where the rows of $x$ and $y$ are denoted with subscripts.  Define $r \in R$ by $r=(r_{ij})$; then $y=rx$.

For the second item, apply the first item twice, symmetrically in $y$ and $x$.  When $\rowsp(y) = \rowsp(x)$, both $x$ and $y$ row reduce to the same row-reduced echelon form, which means they are in the same $\U$-orbit.  
\end{proof}

Let $\mathcal{P}_M$ be the partially ordered set (poset) of all cyclic left $R$-submodules of $M = M_{k \times m}(\F_q)$, and let $\mathcal{P}_{k,m}$ be the poset of all linear subspaces of dimension at most $k$ in $\F_q^m$.  Define $\rho: \mathcal{P}_M \ra \mathcal{P}_{k,m}$ by $\rho(Rx) = \rowsp(x)$; $\rho$ is well-defined by Lemma~\ref{lem:Rx-rowspace}.  Conversely, given a linear subspace $V \subseteq \F_q^m$, define
\[  \psi(V) = \{ x \in M: \rowsp(x) \subseteq V \} . \]

\begin{prop}  \label{prop:poset-isom}
When $\dim V \leq k$, $\psi(V)$ is a cyclic left $R$-submodule of $M$.  The map $\rho: \mathcal{P}_M \ra \mathcal{P}_{k,m}$  is an isomorphism of posets, with inverse given by $\psi$.
\end{prop}

\begin{proof}
Suppose $\dim V \leq k$.  Choose a basis of $V$, and define $x \in M$ to have the chosen basis of $V$ as its first $\dim V$ rows, followed by rows of zeros.  Then $\rowsp(x) = V$.  By Lemma~\ref{lem:Rx-rowspace}, $\psi(V) = Rx$ is a cyclic module.  The argument also shows that $\rho$ is surjective.  
Lemma~\ref{lem:Rx-rowspace} implies that $\rho$ is injective and preserves inclusion.
\end{proof}

\begin{cor}  \label{cor:orbit-correspondence}
The orbits $\orb(x)$, $x \in M$, of the left action of $\U = \GL(k, \F_q)$ on $M= M_{k \times m}(\F_q)$ are in one-to-one correspondence with the linear subspaces of dimension at most $k$ contained in $\F_q^m$.  More precisely, the orbits of matrices of rank $j$ are in one-to-one correspondence with linear subspaces of dimension $j$ in $\F_q^m$.  The linear subspace corresponding to $\orb(x)$ is $\rowsp(x)$.
\end{cor}

There are similar results for linear functionals on $M$.  A \emph{linear functional} on $M$ is a homomorphism $\la: M \ra R$ of left $R$-modules; inputs will be written on the left, so that $(rx)\la = r(x \la)$ for $r \in R$ and $x \in M$.  The collection of all linear functionals is denoted $M^\sharp = \Hom_R(M,R)$; $M^\sharp$ is a right $R$-module, with addition defined point-wise and $\la r$ defined by $x(\la r) = (x\la)r$, where $\la \in M^\sharp$, $r \in R$, $x \in M$.  When $M = M_{k \times m}(\F_q)$, then $M^\sharp = M_{m \times k}(\F_q)$, with the evaluation $x\la \in R$, $x \in M$, $\la \in M^\sharp$, being matrix multiplication.  

The orbits $\orb(\la)$, $\la \in M^\sharp$, of the right action of $\U = \GL(k, \F_q)$ on $M^\sharp= M_{m \times k}(\F_q)$ are in one-to-one correspondence with the linear subspaces of dimension at most $k$ contained in $\F_q^m$.  More precisely, the orbits of linear functionals of rank $j$ are in one-to-one correspondence with linear subspaces of dimension $j$ in $\F_q^m$.  The linear subspace corresponding to $\orb(\la)$ is the column space $\colsp(\la)$.

\begin{rem}   \label{rem:functionalEvaluation}
Given a linear functional $\la \in M^\sharp$, i.e., $\la: M \ra R$, its kernel consists of $x \in M$ such that $\rowsp(x) \subseteq \colsp(\la)^\perp$.  Here, for a linear subspace $Y \subseteq \F_q^m$, denote by $Y^\perp$ its orthogonal with respect to the standard dot product on $\F_q^m$.  Given $\la_i \in M^\sharp$, $i=1, 2, \ldots, n$, define $\La: M \ra R^n$ by $x \La = (x \la_1, \ldots, x \la_n) \in R^n$.  Then $\ker \La$ consists of $x \in M$ with $\rowsp(x) \subseteq \cap \colsp(\la_i)^\perp = (\colsp(\la_1) + \cdots + \colsp(\la_n))^\perp$.  This latter uses the fact that $(X + Y)^\perp = X^\perp \cap Y^\perp$, for linear subspaces $X,Y \subseteq \F_q^m$.  In particular, if $\colsp(\la_1) , \ldots , \colsp(\la_n)$ span $\F_q^m$, then $\La$ is injective.
\end{rem}

We record the number of $\U$-orbits and their sizes, depending on their rank.  The \emph{$q$-binomial coefficient} $\smqbinom{m}{j}$ for $1 \leq j \leq m$ is defined by
\[  \qbinom{m}{j} = \frac{(q^m-1) (q^{m-1}-1) \cdots (q^{m-j+1}-1)}{(q^j-1) (q^{j-1}-1) \cdots (q-1)} . \]
For $m \geq 0$, $\smqbinom{m}{0} = 1$.  If $j<0$ or $j>m$, then $\smqbinom{m}{j}=0$.

\begin{prop}  \label{prop:orbitNumberSize}
There is one orbit of rank $0$, of size $\mathscr{S}_0=1$, in $M$. For any integer $j$, $1 \leq j \leq k$, all orbits of rank $j$ matrices in $M$ have the same size.  The number and size $\mathscr{S}_j$ of orbits of rank $j$ matrices in $M$ are:
\[  \begin{array}{c|c}
\text{\rm number} & \text{\rm size} \\ \hline
\qbinom{m}{j} & \mathscr{S}_j = \displaystyle\prod_{i=0}^{j-1} (q^k - q^i)
\end{array}  \]
The ratio $\mathscr{S}_{j+1}/\mathscr{S}_j$ satisfies $\mathscr{S}_{j+1}/\mathscr{S}_j = q^k - q^j$, for $0 \leq j \leq k-1$.

The size of a cyclic submodule $R r$ depends only on $\rk r$: $\size{R r} = q^{k \rk r}$.
\end{prop}

We note that the sizes $\mathscr{S}_j$ and $\size{R r}$ do not depend on $m$.

\begin{proof}
It is well-known (e.g., \cite[Theorem~3.2.6]{MR4249619}) that the $q$-binomial coefficient $\smqbinom{m}{j}$ counts the number of $j$-dimensional linear subspaces in $\F_q^m$, so the number of orbits follows from Corollary~\ref{cor:orbit-correspondence}.  

In addition to the left action of $\U=\GL(k, \F_q)$ on $M=M_{k \times m}(\F_q)$, there is also a right action of $\GL(m,\F_q)$ on $M$ via matrix multiplication; both actions preserve rank.  As matrix multiplication is associative, these two actions commute.  Thus right multiplication by $P \in \GL(m, \F_q)$ maps the $\U$-orbit $\orb(x)$, $x \in M$, to $\orb(xP)$, and the two orbits have the same size.

Suppose the integer $j$ satisfies $1 \leq j \leq k$.  Choose $x_0 \in M$ to have the first $j$ standard basis vectors (i.e., $(1, 0, \ldots, 0)$, etc.) as its first $j$ rows, with the remaining rows being all zeros.    Pick any $y \in M$ with $\rk y = j$.  Choose a basis of $\rowsp(y)$, and extend it to a basis of $\F_q^m$.  Use this basis of $\F_q^m$ as the rows of a matrix $P \in \GL(m,\F_q)$.  Then the rows of $x_0 P$ consist of the chosen basis of $\rowsp(y)$, followed by $k-j$ zero-rows.  Thus we have $\rowsp(x_0 P) = \rowsp(y)$, so that $\orb(y) = \orb(x_0 P)$, by Lemma~\ref{lem:Rx-rowspace}.  We conclude that $\size{\orb(y)} = \size{\orb(x_0 P)} = \size{\orb(x_0)}$, so that all orbits of rank $j$ matrices have the same size.

As for the size $\mathscr{S}_j$ of an orbit of rank $j$ matrices, it is enough to calculate $\size{\orb(x_0)}$ using 
$\size{\orb(x_0)} = \size{\U}/\size{\stab(x_0)}$, where $\stab(x_0) = \{ u \in \U: u x_0 = x_0 \}$ is the stabilizer subgroup of $x_0$.  Then $u \in \stab(x_0)$ has the form
\[  u =  \left[ \begin{array}{c|c}
I_j & B \\ \hline
0 & D
\end{array} \right]  ,\]
with $I_j$ the $j \times j$ identity matrix, $B$ arbitrary, and $D$ invertible.  Then 
\[  \size{\orb(x_0)} = \size{\GL(k, \F_q)}/ (q^{j(k-j)} \size{\GL(k-j, \F_q)}) , \]
which simplifies as claimed.

The same argument using the right action of $\GL(m, \F_q)$ implies that the size of a cyclic submodule $R r \subseteq M$ depends only on $\rk r$.  Indeed, suppose $\rk r_1 = \rk r_2$.  By row and column operations, there are units $u_1 \in \GL(k, \F_q)$ and $u_2 \in \GL(m, \F_q)$ such that $r_2 = u_1 r_1 u_2$.  By  Lemma~\ref{lem:Rx-rowspace}, $R r_2 = R r_1 u_2$.  Right multiplication by $u_2$ maps $R r_1$ isomorphically to $R r_1 u_2$.  We conclude that $\size{R r_2} = \size{R r_1 u_2} = \size{R r_1}$.  For $j = 1, 2, \ldots, k$, let $r \in M$ be the following matrix of rank $j$:
\[  r =  \left[ \begin{array}{c|c}
I_j & 0 \\ \hline
0 & 0
\end{array} \right]  .\]
Then $R r$ consists of all matrices in $M$ whose last $m-j$ columns are zero.  The cyclic submodule $R r$ has size $\size{R r} = q^{kj} = q^{k \rk r}$.
\end{proof}

We now turn our attention to the homogeneous weight on $R = M_{k \times k}(\F_q)$.
Because of Proposition~\ref{prop:poset-isom}, the M\"obius function for the poset of principal left ideals of $R$ equals the M\"obius function for the poset of linear subspaces of $\F_q^k$, which, following \cite[(2.7)]{phall:eulerian-functions}, is 
\begin{equation}  \label{eq:matrixMoebius}
\mu(V_1, V_2) = (-1)^c q^{\binom{c}{2}} , \quad V_1 \subseteq V_2 \subseteq \F_q^k, 
\end{equation}
where $c = \dim V_2 - \dim V_1$ is the codimension of $V_1$ in $V_2$.

Equation \eqref{eq:GS-homogwt} yields the following formula for the homogeneous weight $\wg$ on $R$; this formula also appears in \cite[Proposition~7]{Honold-Nechaev:weighted-modules}.  We write $\rho=\rk(r)$, $r \in R$:
\begin{equation}  \label{eq:MatrixHomogWeight}
\wg(r) = \begin{cases}
0, & \rho=0, \\
\zeta \left( 1 -  \frac{(-1)^{\rho}}{(q^k -1)(q^{k-1} -1) \cdots (q^{k-\rho +1}-1)} \right), & \rho >0.
\end{cases}
\end{equation}
Note that $\wg(r)$ depends only on $\rho$, which is consistent with $\wg$ being constant on left $\U$-orbits.  Write $\wg_{\rho}$ for the common value of $\wg(r)$ where $\rk r = \rho$.
By choosing $\zeta$ appropriately, namely
\[  \zeta = (q^k-1)(q^{k-1}-1) \cdots (q-1)/q,  \]
the homogeneous weight will have integer values.  

\begin{ex}  \label{ex:2by2case}
For $M_{2 \times 2}(\F_q)$, the homogeneous weight is
\[  \begin{array}{c|ccc}
 & \wg_0 & \wg_1 & \wg_2 \\ \hline
\text{general $q$} & 0 & q^2 -q & q^2 - q -1 \\
q=2 & 0 & 2 & 1 \\
q=3 & 0 & 6 & 5
\end{array} , \]
with average weight $\zeta= (q^2 -1)(q-1)/q$ in general, so that $\zeta = 3/2$ for $q=2$, and $\zeta = 16/3$ for $q=3$.
\end{ex}

\begin{ex}
For $M_{3 \times 3}(\F_q)$, the homogeneous weight is:
\[  \begin{array}{c|c|c|c|c}
 & \wg_0 & \wg_1 & \wg_2 & \wg_3 \\ \hline
q & 0 & q^5 -q^4 -q^3 +q^2 & q^5 -q^4 -q^3 +q & q^5 -q^4 -q^3 +q +1 \\
q=2 & 0 & 12 & 10 & 11 \\
q=3 & 0 & 144 & 138 & 139 
\end{array} , \]
with average weight $\zeta= (q^3 - 1)(q^2 -1)(q-1)/q$ in general, so that $\zeta = 21/2$ for $q=2$, and $\zeta = 416/3$ for $q=3$.
\end{ex}

\begin{lem}  \label{lem:order-of-homog}
The homogeneous weight on $M_{k \times k}(\F_q)$ satisfies 
\[  0 = \wg_0 < \wg_2 < \wg_4 < \cdots < \zeta < \cdots <\wg_3 < \wg_1 . \]
Moreover, $2\wg_2 - \wg_1 >0$ for all $k \geq 2$, $q\geq 2$, \emph{except} for $k=q=2$, where $2\wg_2 - \wg_1 =0$.
\end{lem}

\begin{proof}
The denominator in \eqref{eq:MatrixHomogWeight} is an increasing function of $\rho = \rk(r)$.   This, together with the alternating sign of $(-1)^\rho$, yields the inequalities among the $\wg_{\rho}$.  By moving $\zeta$ to the left side of \eqref{eq:MatrixHomogWeight}, one sees that $\wg_{\rho} - \zeta = -(-1)^{\rho} \zeta/((q^k-1)\cdots (q^{k-\rho+1}-1))$.  This implies that $\wg_{\rho}-\zeta$ is positive when $\rho$ is odd and negative when $\rho$ is even.

For $k \geq 2$, one calculates that 
\[  2 \wg_2 - \wg_1 = \zeta \frac{(q^k-2)(q^{k-1}-1) -2}{(q^k-1)(q^{k-1}-1)} . \]
Using $k \geq 2$ and $q \geq 2$, the numerator satisfies
\[  (q^k-2)(q^{k-1}-1) -2 \geq (q^2 -2)(q-1)-2 = q(q+1)(q-2). \]
This last expression is positive when $q>2$ and vanishes when $q=2$.  Even for $q=2$, the earlier inequality is strict when $k>2$.  Thus, $2 \wg_2 - \wg_1 >0$, except for $k=q=2$, where $2 \wg_2 - \wg_1 =0$.
\end{proof}

Using Lemma~\ref{lem:singletons}, we see that any nonzero vector $v$ with $\wg(v) < \wg_1$ must be a singleton.  Any nonzero vector with $\wg(v) = \wg_1$ must be a singleton (or a doubleton, i.e., two nonzero entries, only for $M_{2 \times 2}(\F_2)$).

\begin{thm}  \label{thm:MWIds-homog-2by2-F2}
The MacWilliams identities hold for the homogeneous weight over $R=M_{2 \times 2}(\F_2)$.  For a linear code $C \subseteq R^n$, 
\[  \howe_{C^\perp}(X,Y) = \frac{1}{\size{C}} \howe_C(X+3Y, X-Y) .  \]
\end{thm}

\begin{proof}
As in the proof of Theorem~\ref{thm:chainMacWids}, we provide details to be used in the argument outlined in Appendix~\ref{sec:appendix}.

From Example~\ref{ex:2by2case}, we have that $\wg_0=0$, $\wg_1=2$, and $\wg_2=1$.  A generating character for $R$ is $\chi(r) = (-1)^{\tr r}$, $r \in R$, where $\tr$ is the matrix trace.
Define $f: R \ra \C[X,Y]$ by $f(r) = X^{2-\wg(r)} Y^{\wg(r)}$.  The value of $f(r)$ depends only on $\rk r$:
\[  \begin{array}{c|ccc}
\rk r & 0 & 1 & 2 \\ \hline
f(r) & X^2 & Y^2 & XY
\end{array}. \]
A calculation shows the Fourier transform \eqref{eqn:FTdefn} depends only on $\rk r$:
\[  \ft{f}(r) = \begin{cases}
X^2 + 9 Y^2 + 6 XY = (X+3Y)^2, & \rk r = 0, \\
X^2 +  Y^2 - 2 XY = (X-Y)^2, & \rk r = 1, \\
X^2 - 3 Y^2 + 2 XY = (X+3Y)(X-Y), & \rk r = 2. \\
\end{cases} \]
Note that $\ft{f}(r)$ has the form of $f(r)$ with a linear substitution: $X \leftarrow X+3Y$, $Y \leftarrow X-Y$.
Applying these details, the rest of the argument in Appendix~\ref{sec:appendix} carries through.
\end{proof}

There is a Gray map from $M_{2 \times 2}(\F_2)$ equipped with the homogeneous weight to $\F_4^2$ equipped with the Hamming weight (of $\F_4$) \cite{MR3146951}.

\section{\texorpdfstring{$W$}{W}-matrix}  \label{sec:W-matrix}
In this section we determine the $W$-matrix of \eqref{eq:DefnOfWmatrix} for a weight $w$ on $R=M_{k \times k}(\F_q)$ that has maximal symmetry. 

As usual, let $R=M_{k \times k}(\F_q)$, and suppose $w$ is a weight on $R$ with maximal symmetry.  Suppose $r \in R$.  By row and column reduction there exist units $u_1, u_2 \in \U = \GL(k,\F_q)$ such that 
\[  u_1 r u_2 = \left[ \begin{array}{c|c}
I_{\rho} & 0 \\ \hline
0 & 0
\end{array} \right]  ,\]
where $\rho = \rk r$.  Thus $w(r) = w(u_1 r u_2) = w(\left[ \begin{smallmatrix}
I_{\rho} & 0 \\ 0 & 0
\end{smallmatrix} \right])$, which says that the value of $w(r)$ depends only on the rank of $r$.  Write $w_0, w_1, \ldots, w_k$ for the value of $w$ on matrices of rank $0, 1, \ldots, k$, respectively.

\begin{rem}  \label{rem:w0AsIndeterminate}
While $w(0)=0$ is part of the definition of a weight, some of the results of this section will be more natural to state if we allow $w_0$ to be viewed as an indeterminate.  We will proceed with $w_0$ as an indeterminate, and later show, in Theorem~\ref{thm:W0-block-diagonal-form} and Corollary~\ref{cor:whenisW0invertible}, how the general results are affected when we set $w_0=0$.
\end{rem}

The information module $M$ will be $M = M_{k \times m}(\F_q)$ with $m \geq k$.  Then $\Hom_R(M,R) = M_{m \times k}(\F_q)$, 
achieved by right multiplication against $M$; i.e., the evaluation pairing $M \times \Hom_R(M,R) \ra R$ sends $x \in M$ and $\la \in \Hom_R(M,R)$ to $x \la \in R$.

Because of maximal symmetry, the symmetry groups of $w$ are $\glt=\grt=\GL(k, \F_q)$.  The orbit space $\OO = \glt \backslash M$ is represented by row-reduced echelon matrices of size $k \times m$, 
and $\OS = \Hom_R(M,R)/\grt$ is represented by column-reduced echelon matrices of size $m \times k$.  
The matrix transpose maps $\OO \leftrightarrow \OS$ bijectively.
The sets $\OO$ and $\OS$ are partitioned by rank, and elements of $\OO$ and $\OS$ correspond to linear subspaces of $\F_q^m$ of dimension at most $k$, by Corollary~\ref{cor:orbit-correspondence}. (Left orbits in $\OO$ are viewed in terms of row spaces; their transposes, right orbits in $\OS$, are viewed in terms of column spaces.) 

The rows of $W$ are indexed by elements $[x] \in \OO$, ordered so that ranks go from $0$ to $k$.  Similarly, the columns of $W$ are indexed by elements $[\la] \in \OS$, ordered to match $\OO$ under the bijection $\OO \leftrightarrow \OS$.  The $[x], [\la]$-entry of $W$ is simply $w(x \la)$, i.e., the value of $w$ at the evaluation $x \la \in R$.  The value $w(x \la)$ is well-defined by the definition of the symmetry groups.  By maximal symmetry, the value $w(x \la)$ depends only on the rank $\rk(x \la)$.  The matrix $W$ is square of size $N_{k,m}$, the number of linear subspaces of dimension at most $k$ in $\F_q^m$.

Suppose $[x] \in \OO$ corresponds to the linear subspace $X \subseteq \F_q^m$ and $[\la] \in \OS$ corresponds to $Y$.  We seek to express $\rk(x \la)$, and hence $w(x \la)$, in terms of $X$ and $Y$.  

For a linear subspace $X \subseteq \F_q^m$, denote by $X^\perp$ its orthogonal with respect to the standard dot product on $\F_q^m$.  Then $\dim X^\perp = m - \dim X$ and $(X^\perp)^\perp = X$, for all linear subspaces $X \subseteq \F_q^m$.  Also, $(X \cap Y)^\perp = X^\perp + Y^\perp$, for linear subspaces $X,Y \subseteq \F_q^m$.

\begin{lem}  \label{lem:dim-formulas}
For linear subspaces $X,Y \subseteq \F_q^m$ representing $[x] \in \OO$ and $[\la] \in \OS$, respectively:
\begin{enumerate}
\item $\dim X - \dim (X \cap Y^\perp) = \dim Y - \dim(Y \cap X^\perp)$, and
\item  $\rk(x \la) = \dim X - \dim (X \cap Y^\perp) = \dim Y - \dim(Y \cap X^\perp)$.
\end{enumerate}
\end{lem}

\begin{proof}
Consider $(X \cap Y^\perp)^\perp = X^\perp + Y$, and compare dimensions:
\begin{align*}
m - \dim(X \cap Y^\perp) &= \dim X^\perp + \dim Y - \dim(X^\perp \cap Y) \\
&= m -\dim X + \dim Y - \dim(X^\perp \cap Y).
\end{align*}
We conclude that $\dim X - \dim(X \cap Y^\perp) =  \dim Y - \dim(X^\perp \cap Y)$.

Choose as a representative $x$ a $k \times m$ matrix whose first $\dim(X \cap Y^\perp)$ rows form a basis for $X \cap Y^\perp$, whose next $\dim X - \dim(X \cap Y^\perp)$ rows complete to a basis of $X$, and whose remaining rows are zeros.  Choose a representative $\la$ by reversing the roles of $X$ and $Y$: its first $\dim(Y \cap X^\perp)$ columns form a basis for $Y \cap X^\perp$, its next $\dim Y - \dim(Y \cap X^\perp)$ columns complete to a basis of $Y$, and its remaining columns are zeros.  Then $x \la$ has the form
\[  x \la = \left[ \begin{array}{c|c|c}
0& 0 & 0 \\ \hline
0 & Z & 0 \\ \hline
0 & 0 & 0
\end{array} \right], \]
where $Z$ is a square matrix of size $(\dim X - \dim(X \cap Y^\perp)) \times ( \dim Y - \dim(X^\perp \cap Y)$.  The matrix $Z$ is nonsingular.  If not, there exists a nonzero vector $v$ with $vZ=0$.  Then $[0|v|0] x \in Y^\perp$, which violates the choice of basis of $X \cap Y^\perp$ in the construction of $x$.  The formula for $\rk(x \la)$ now follows.
\end{proof}

The matrix $W$ is symmetric when we use the bijection $\OO \leftrightarrow \OS$ to align the indexing.

Because some of our later results depend upon inverting $W$, we need to understand when the matrix $W$ is invertible.  We will be able to transform $W$ into a block diagonal format by making use of the M\"obius function of the poset $\mathcal{P}_{k,m}$ of linear subspaces of dimension at most $k$ in $\F_q^m$.  Versions of this block diagonal format can be found in \cite[\S 4]{MR3207471}, \cite[Theorem~9.6]{wood:turkey}, and \cite[\S 6]{MR3866771}.

Recall that we index the rows and columns of $W$ by linear subspaces of dimension at most $k$ in $\F_q^m$, with ranks increasing from $0$ to $k$.

Define a matrix $P$, with rows and columns indexed by linear subspaces of dimension $\leq k$ in $\F_q^m$, using the same ordering as for $W$.  The entry $P_{\al,\be}$ is given by
\begin{equation}  \label{eq:DefnOfPmatrix}
P_{\al,\be} = \begin{cases}
\mu(0,\be) = (-1)^{\dim\be} q^{\binom{\dim\be}{2}}, & \text{ if $\be \subseteq \al$}, \\
0, & \text{ if $\be \not\subseteq \al$}.
\end{cases}
\end{equation}
Because we are ordering rows and columns so that ranks increase, we see that $P$ is lower triangular.  Its diagonal entries are $P_{\al,\al} = (-1)^{\dim \al} q^{\binom{\dim\al}{2}} \neq 0$.  Thus $P$ is invertible over $\Q$.

For $j=0,1, \ldots, k$, define an incidence matrix $\mathscr{I}_j$ over $\Q$, square of size $\smqbinom{m}{j}$,  with rows and columns indexed by linear subspaces of dimension $j$ in $\F_q^m$, using the dimension $j$ portion of the ordering used for $W$ and $P$.  The $\al, \delta$-entry of $\mathscr{I}_j$ is given by
\[  (\mathscr{I}_j)_{\al,\delta} = \begin{cases}
1, & \al \cap \delta^\perp = 0, \\
0, & \al \cap \delta^\perp \neq 0.
\end{cases}  \]
The incidence matrices $\mathscr{I}_j$ are invertible by \cite[Proposition~6.7]{MR3866771}.

Our main objective in this section is to prove the next theorem.
\begin{thm}  \label{thm:W-block-diagonal-form}
For positive integers $2 \leq k \leq m$ and a weight $w$ on $M_{k \times k}(\F_q)$ having maximal symmetry and $w_0$ indeterminate, we have 
\[ PWP^\top = \left[ \begin{array}{rrrr}
c_0 \mathscr{I}_0 & 0 & \cdots & 0 \\
0 & c_1 \mathscr{I}_1 & \cdots & 0 \\
\vdots & & \ddots & \vdots \\
0 & 0 & \cdots & c_k \mathscr{I}_k 
\end{array} \right], \]
where, for $j=0, 1, \ldots, k$,
\begin{equation}  \label{eqn:diag-coefficents}
c_j = (-1)^j q^{\binom{j}{2}} \sum_{\ell =0}^j (-1)^\ell q^{\binom{\ell}{2}} \qbinom{j}{\ell} w_{\ell} .
\end{equation}
\end{thm}

Before we prove Theorem~\ref{thm:W-block-diagonal-form}, we prove some preliminary lemmas and propositions that will be used in the proof.  We begin by quoting the well-known Cauchy Binomial Theorem, e.g., \cite[Theorem~3.2.4]{MR4249619}.

\begin{thm}[Cauchy Binomial Theorem] \label{thm:CauchyBT}
For a positive integer $k$,
\[  \prod_{i=0}^{k-1} (1+x q^i) = \sum_{j=0}^k q^{\binom{j}{2}} \qbinom{k}{j} x^j .    \]
In particular, using $x=-1$, for $k$ positive,
\begin{equation}  \label{eq:Cauchy-vanishing}
  \sum_{j=0}^k (-1)^j q^{\binom{j}{2}} \qbinom{k}{j} =0 .  
\end{equation}
\end{thm}

In a vector space, a \emph{frame} is an ordered set of linearly independent vectors; if there are $b$ such vectors, we call the frame a $b$-frame.

\begin{lem}  \label{lem:count-hit-trivially}
Let $V$ be a vector space over $\F_q$ with $\dim V = v$, and let $D$ be a linear subspace of $V$ with $\dim D = d$.  Then
\[ \size{\{B \subseteq V: \dim B = b, B \cap D = 0\}} = q^{bd} \qbinom{v-d}{b} . \]
\end{lem}

\begin{proof}
We count the number of $b$-frames outside of $D$, and divide by $\size{\GL(b, \F_q)}$, \eqref{eq:size-of-GLk}.  Then, factoring out $b$ factors of $q^d$ from
\[ \frac{(q^v-q^d) (q^v - q^{d+1}) \cdots (q^v-q^{d+b-1})}{(q^b-1)(q^b-q) \cdots (q^b-q^{b-1})} \]
yields the stated result.
\end{proof}

\begin{lem}  \label{lem:count-hit}
Let $V$ be a vector space over $\F_q$ with $\dim V = v$, and let $D$ be a linear subspace of $V$ with $\dim D = d$.  For any $j=1,2, \ldots, d$, 
\[ \size{\{B \subseteq V: \dim B = b, \dim(B \cap D) = j \}} = q^{(b-j)(d-j)} \qbinom{d}{j} \qbinom{v-d}{b-j} . \]
\end{lem}

\begin{proof}
The count equals the number of $j$-dimensional subspaces $J \subseteq D$ times the number of $B$'s of dimension $b$ with $B \cap D = J$.  The number of $j$-dimensional subspaces of $D$ is $\smqbinom{d}{j}$.  The set of $b$-dimensional subspaces $B \subseteq V$ with $B \cap D = J$ is in one-to-one correspondence with the set of $(b-j)$-dimensional subspaces of $V/J$ that intersect $D/J$ trivially.  By Lemma~\ref{lem:count-hit-trivially}, the number of such subspaces is $q^{(b-j)(d-j)} \smqbinom{v-d}{b-j}$.
\end{proof}
Lemma~\ref{lem:count-hit} sharpens \cite[Lemma~6.9]{MR3866771}; the latter's $C_1(b)$ is $q^{b(m-b)}$.

\begin{lem}  \label{lem:counting-intersections}
Let $V$ be a vector space over $\F_q$ with $\dim V = v$, and let $A$ and $D$ be linear subspaces of $V$ with $\dim A =a$ and $\dim D = d$.  If $\dim(A \cap D) = i$, then, for any $j=0,1, \ldots, i$, 
\[ \size{\{B \subseteq A: \dim B = b, \dim(B \cap D) = j \}} = q^{(b-j)(i-j)} \qbinom{i}{j} \qbinom{a-i}{b-j} . \]
\end{lem}

\begin{proof}
Note that $B \cap D = B \cap (A \cap D)$.  Apply Lemmas~\ref{lem:count-hit-trivially} and~\ref{lem:count-hit} with ambient space $A$ and subspace $A \cap D$.
\end{proof}

\begin{lem}  \label{lem:dimOfIntersection}
Suppose $\al, \delta \subseteq \F_q^m$ are linear subspaces.  Then $\dim(\al \cap \delta^\perp) \geq \dim \al - \dim \delta$.  If $\dim \al > \dim \delta$, then $\dim(\al \cap \delta^\perp) >0$.
\end{lem}

\begin{proof}
Using $\al + \delta^\perp \subseteq \F_q^m$, compare dimensions:
\begin{align*}
m &\geq \dim(\al + \delta^\perp) 
= \dim \al +m-\dim \delta - \dim(\al \cap \delta^\perp),
\end{align*}  
from which the result follows.
\end{proof}

\begin{prop}  \label{prop:PW-format}
If $\dim(\al \cap \delta^\perp) > 0$, then $(P W)_{\al, \delta} = 0$.  If $\al \cap \delta^\perp = 0$, then, writing $a = \dim \al$,
\[  (P W)_{\al, \delta} = \sum_{\be \subseteq \al} \mu(0,\be) w_{\dim \be} 
= \sum_{r=0}^a (-1)^r q^{\binom{r}{2}} \qbinom{a}{r} w_r . \]
In particular, the matrix $P W$ is block upper triangular.

Likewise, if $\dim(\ga \cap \be^\perp) > 0$, then $(W P^\top )_{\be, \ga} = 0$.  If $\ga \cap \be^\perp = 0$, then, writing $c = \dim \ga$,
\[  (W P^\top )_{\be, \ga}  = \sum_{\epsilon \subseteq \ga} \mu(0,\epsilon) w_{\dim \epsilon} 
= \sum_{s=0}^c (-1)^s q^{\binom{s}{2}} \qbinom{c}{s} w_s . \]
In particular, the matrix $W P^\top$ is block lower triangular.
\end{prop}

\begin{proof}
From the definition of the matrix $P$ and Lemma~\ref{lem:dim-formulas}, 
\[  (P W)_{\al, \delta} = \sum_{\be \subseteq \al} \mu(0,\be) w_{\dim \be - \dim(\be \cap \delta^\perp)}. \] 
In this formula, the subscript $\dim \be - \dim(\be \cap \delta^\perp) = \dim \delta - \dim(\delta \cap \be^\perp) \leq \dim \delta - \dim(\delta \cap \al^\perp) = \dim \al - \dim(\al \cap \delta^\perp)$, as $\be \subseteq \al$, so that $\al^\perp \subseteq \be^\perp$.  Writing $i = \dim(\al \cap \delta^\perp)$, we see that the subscript $\dim \be - \dim(\be \cap \delta^\perp) \leq a-i$.  Thus
\begin{align*}
(P W)_{\al, \delta} &= \sum_{r=0}^{a-i} \sum_{\substack{\be \subseteq \al \\ \dim \be - \dim(\be\cap \delta^\perp)=r}} \mu(0,\be) w_r \\
&= \sum_{r=0}^{a-i} \sum_{j=0}^i \sum_{\substack{\be \subseteq \al \\ \dim \be = r + j \\ \dim(\be\cap \delta^\perp)=j}} \mu(0,\be) w_r .
\end{align*}  
By Lemma~\ref{lem:counting-intersections} and \eqref{eq:matrixMoebius}, the coefficient $C_r$ of $w_r$ is
\begin{align*}
C_r &= \sum_{j=0}^i q^{r(i-j)} \qbinom{i}{j} \qbinom{a-i}{r} (-1)^{r+j} q^{\binom{r+j}{2}} \\
&= (-1)^r q^{\binom{r}{2}} q^{i r} \qbinom{a-i}{r} \left(\sum_{j=0}^i (-1)^j q^{\binom{j}{2}} \qbinom{i}{j}\right) \\
&= \begin{cases}
(-1)^r q^{\binom{r}{2}} \smqbinom{a}{r}, & i=0, \\
0 , & i > 0.
\end{cases}
\end{align*}
To simplify, we used the identity $\binom{r+j}{2} = \binom{r}{2} + j r + \binom{j}{2}$ and \eqref{eq:Cauchy-vanishing}.

If $\dim \delta < \dim \al$, Lemma~\ref{lem:dimOfIntersection} implies $i=\dim(\al \cap \delta^\perp) > 0$.  Using dimension to create blocks, we see that $P W$ is block upper triangular.

Essentially the same arguments yield the results about $W P^\top$. 
\end{proof}

Recall that the matrix $P$ is lower triangular, so that $P^\top$ is upper triangular.  Thus $P W P^\top$ will be both lower and upper block triangular, hence block diagonal.  The exact form of $P W P^\top$ is the next result.

\begin{prop}  \label{prop:W-block-diagonal}
If $\dim \al \neq \dim \delta$, then $(P W P^\top)_{\al, \delta} = 0$.  If $\dim \al = \dim\delta= a$, then
\[  (P W P^\top)_{\al, \delta} = \begin{cases}
(-1)^a q^{\binom{a}{2}} \sum_{j=0}^a (-1)^j q^{\binom{j}{2}} \smqbinom{a}{j} w_j, & \al \cap \delta^\perp = 0, \\
0, & \al \cap \delta^\perp \neq 0.
\end{cases} \]
\end{prop}

\begin{proof}
We first show that $\al \cap \delta^\perp \neq 0$ implies $(P W P^\top)_{\al, \delta} =0$.  Assume $\al \cap \delta^\perp \neq 0$.  For any $\ga \subseteq \delta$, we have $\delta^\perp \subseteq \ga^\perp$, so that $\al \cap \delta^\perp \subseteq \al \cap \ga^\perp$.  Thus $\al \cap \ga^\perp \neq 0$ for all $\ga \subseteq \delta$.

Using the definition of $P$, we see that
\begin{align}  \label{eq:PWP-formula}
(P W P^\top)_{\al, \delta} &= \sum_{\ga \subseteq \delta} (P W)_{\al,\ga} \mu(0,\ga) .
\end{align}
By Proposition~\ref{prop:PW-format}, all the $(P W)_{\al,\ga}$-terms in \eqref{eq:PWP-formula} vanish, so that $(P W P^\top)_{\al, \delta} =0$, as claimed.

Essentially the same argument using 
\[ (P W P^\top)_{\al, \delta} = \sum_{\be \subseteq \al} \mu(0,\be) (W P^\top)_{\be,\delta} \]
shows that $\delta \cap \al^\perp \neq 0$ implies $(P W P^\top)_{\al, \delta} =0$.

From Lemma~\ref{lem:dim-formulas} we have the equation $\dim \al - \dim \delta = \dim(\al \cap \delta^\perp) - \dim(\delta \cap \al^\perp)$.  If $\dim \al \neq \dim \delta$, at least one of $\al \cap \delta^\perp$ or $\delta \cap \al^\perp$ is nonzero.  Thus $(P W P^\top)_{\al, \delta} =0$.

At last, suppose $\dim \al = \dim \delta$.  In this situation, note that $\al \cap \delta^\perp = 0$ if and only if $\delta \cap \al^\perp =0$.  As above, if $\al \cap \delta^\perp \neq 0$, then $(P W P^\top)_{\al, \delta} =0$.  If $\al \cap \delta^\perp = 0$, then the equality of dimensions yields
$\F_q^m = \al \oplus \delta^\perp$.  For any $\ga \subsetneq \delta$, we have $\delta^\perp \subsetneq \ga^\perp$, so that $\al \cap \ga^\perp \neq 0$.  By Proposition~\ref{prop:PW-format}, $(P W)_{\al,\ga} =0$ for $\ga \subsetneq \delta$.  Thus, \eqref{eq:PWP-formula} implies that 
\[ (P W P^\top)_{\al, \delta} =  (P W)_{\al,\delta}  \mu(0,\delta),  \]
which gives the stated formula.
\end{proof}

\begin{proof}[Proof of Theorem~{\textup{\ref{thm:W-block-diagonal-form}}}]
The nonzero entries in Proposition~\ref{prop:W-block-diagonal} depend only on $\dim\al$ and equal the $c_j$ in the statement of Theorem~\ref{thm:W-block-diagonal-form}.  Whether $\al \cap \delta^\perp =0$ is marked by the incidence matrix $\mathscr{I}_j$.
\end{proof}

\begin{cor}  \label{cor:whenisWinvertible}
The matrix $W$ is invertible over $\Q$ if and only if
\[  c_j = (-1)^j q^{\binom{j}{2}} \sum_{\ell =0}^j (-1)^\ell q^{\binom{\ell}{2}} \qbinom{j}{\ell} w_{\ell} \neq 0, \]
for all $j=0, 1, \ldots, k$.
\end{cor}

\begin{proof}
The matrices $P$, $P^\top$, and $\mathscr{I}_j$ are all invertible over $\Q$.
\end{proof}

Recall from Remark~\ref{rem:w0AsIndeterminate} that we have been treating $w_0$ as an indeterminate.  When $w_0$ is set equal to $0$, the first row and first column of the matrix $W$ consist of $0$'s.  Then $W$ cannot be invertible.  Equivalently, $c_0=0$ in Corollary~\ref{cor:whenisWinvertible}.

To get around this lack of invertibility, we make the following adjustments, as in Remark~\ref{rem:W0-map}.  Define a matrix $W_0$ with rows and columns indexed by the \emph{nonzero} elements of $\OO$ and $\OS$, respectively, ordered so that ranks go from $1$ to $k$.  The $[x],[\la]$-entry is again $w(x \la)$, with $w(0)=0$ now.  Similarly, define a matrix $P_0$ with rows and columns indexed by the nonzero elements of $\OO$ and $(\al,\be)$-entry given by \eqref{eq:DefnOfPmatrix}.  Then the counterparts of Theorem~\ref{thm:W-block-diagonal-form} and Corollary~\ref{cor:whenisWinvertible} are the following.

\begin{thm}  \label{thm:W0-block-diagonal-form}
For positive integers $2 \leq k \leq m$ and a weight $w$ on $M_{k \times k}(\F_q)$ having maximal symmetry and $w(0)=0$, we have 
\[ P_0 W_0 P_0^\top = \left[ \begin{array}{rrrr}
c_1 \mathscr{I}_1 & 0 & \cdots & 0 \\
0 & c_2 \mathscr{I}_2 & \cdots & 0 \\
\vdots & & \ddots & \vdots \\
0 & 0 & \cdots & c_k \mathscr{I}_k 
\end{array} \right], \]
where, for $j=1, \ldots, k$,
\begin{equation}  \label{eqn:0-diag-coefficents}
c_j = (-1)^j q^{\binom{j}{2}} \sum_{\ell =1}^j (-1)^\ell q^{\binom{\ell}{2}} \qbinom{j}{\ell} w_{\ell} .
\end{equation}
\end{thm}

\begin{proof}
The relationships between the matrices $P$ and $P_0$ and between $W$ and $W_0$, when $w_0=0$, are given below.  The notations $\row(a)$ and $\col(a)$ mean a row, resp., column, vector, all of whose entries are $a$.
\[ P = \left[ \begin{array}{c|c}
1 & \row(0) \\ \hline
\col(1) & P_0
\end{array} \right], \quad
W|_{w_0=0} = \left[ \begin{array}{c|c}
0 & \row(0) \\ \hline
\col(0) & W_0
\end{array} \right] .
 \]
Then
\[  (P W P^\top)|_{w_0=0} = \left[ \begin{array}{c|c}
0 & \row(0) \\ \hline
\col(0) & \vphantom{\hat{P}}P_0 W_0 P_0^\top
\end{array} \right] ,
 \]
 and the result follows from Theorem~\ref{thm:W-block-diagonal-form}, with $w_0=0$.
\end{proof}

\begin{cor}  \label{cor:whenisW0invertible}
The matrix $W_0$ is invertible over $\Q$ if and only if
\[  c_j = (-1)^j q^{\binom{j}{2}} \sum_{\ell =1}^j (-1)^\ell q^{\binom{\ell}{2}} \qbinom{j}{\ell} w_{\ell} \neq 0, \]
for all $j=1, \ldots, k$.
\end{cor}

\begin{rem}  \label{rem:EPand W0}
The extension property (EP) for $w$ holds when the $W_0$ map is injective (zero right null space) for all information modules.  We see that EP holds if and only if all $c_j \neq 0$ for $j=1,2, \ldots, k$.  See \cite[Theorem~9.5]{wood:turkey}.  In particular, when $W_0$ is invertible, then $W_0: F_0(\OS, \Q) \ra F_0(\OO, \Q)$ is an isomorphism, Remark~\ref{rem:W0-map}.
\end{rem}

\begin{rem}
It is possible to generalize Propositions~\ref{prop:PW-format} and~\ref{prop:W-block-diagonal} to the context of an alphabet $A=M_{k \times \ell}(\F_q)$ and information module $M=M_{k \times m}(\F_q)$, with $m \geq \ell \geq k$.  This paper does not require such a level of generality.
\end{rem}

\section{Locally constant functions}

In this section we examine a type of multiplicity function that will feature prominently in the construction of linear codes in Section~\ref{sec:construction}.

As in previous sections, $R = M_{k \times k}(\F_q)$ and $M=M_{k \times m}(\F_q)$, with $k \leq m$.  The ring $R$ is equipped with a weight $w$ having maximal symmetry and positive integer values.  As in Remark~\ref{rem:w0AsIndeterminate}, we will treat $w_0$ as an indeterminate.
Recall that the orbit spaces $\OO$ and $\OS$ are ordered by rank, from $0$ to $k$.  Write $\OO_i$ and $\OS_i$ for the collections of orbits of rank equal to $i$.  Recall that orbits are denoted by $\orb(x)$ or $[x]$.  Given sets $B \subseteq A$, we denote the indicator function of $B$ by $1_B: A \ra \Z$, with
\[ 1_B (x) = \begin{cases}
1, & x \in B, \\
0, & x \not\in B.
\end{cases}  \]
The type of multiplicity function to be considered has the form $1_{\OS_i}$ or a linear combination of such indicator functions.

Given a multiplicity function $\eta: \OS \ra \N$ (or, more generally, $\eta: \OS \ra \Q$), we refer to $\omega = W \eta$ as the \emph{list of orbit weights} of $\eta$.  The list $\omega$ is a function $\omega: \OO \ra \Q$.  We say that $\eta$, resp., $\omega$, is \emph{locally constant} if $\eta = \sum_{j = 0}^k a_j 1_{\OS_j}$, resp., $\omega = \sum_{j=0}^k b_j 1_{\OO_j}$, for rational constants $a_j, b_j$.  Said another way, $\eta$ is locally constant if $\rk \la_1 = \rk \la_2$ implies $\eta([\la_1]) = \eta([\la_2])$; i.e., $\eta$ is constant on each $\OS_j$.   Similar comments apply to $\omega$.

Let $\eta_j = 1_{\OS_j}$, and $\omega_j = W \eta_j$.  We show that $\omega_j$ is locally constant.

\begin{prop}  \label{prop:loc-constant-eta}
Let $\eta_j = 1_{\OS_j}$.  Then $\omega_j = W \eta_j$ is locally constant.  If $i = \rk x$, then 
\[  \omega_j([x]) =  
\sum_{d=0}^{\min\{j,m-i\}}  q^{(m-i-d)(j-d)} \qbinom{m-i}{d} \qbinom{i}{j-d} w_{j-d}. \]
In particular, when $j=1$, $i = \rk x = 1, 2, \ldots, k$, and $w_0=0$,
\[ \omega_1([x]) = q^{m-i} \qbinom{i}{1} w_1 = (q^{m-i} + q^{m-i+1} + \cdots + q^{m-1}) w_1, \]
is an increasing function of $i$.
\end{prop}

\begin{proof}
Because $\omega_j = W \eta_j$, we have
\[ \omega_j([x]) = \sum_{[\la] \in \OS} w(x \la) \eta_j([\la]) = \sum_{[\la] \in \OS_j} w(x \la) . \]

Using Lemma~\ref{lem:dim-formulas}, represent $\orb(x)=[x]$ by a linear subspace $X \subseteq \F_q^m$, $\dim X = i$, and $\orb(\la)=[\la]$ by $Y$, $\dim Y = j$.  Then $w(x\la) = w_{i - \dim(X \cap Y^\perp)} = w_{j - \dim(Y \cap X^\perp)}$.  Lemma~\ref{lem:count-hit} now implies
\begin{align*}
\omega_j([x]) 
&= \sum_{d=0}^{\min\{j,m-i\}} \size{\{Y: \dim Y = j, \dim(Y \cap X^\perp) = d\}} w_{j-d} \\
&= \sum_{d=0}^{\min\{j,m-i\}} q^{(m-i-d)(j-d)} \qbinom{m-i}{d} \qbinom{i}{j-d} w_{j-d} .
\end{align*}
When $j=1$, the formula simplifies as stated.
\end{proof}

\begin{cor}  \label{cor:LocallyConstant}
If $\eta$ is locally constant, then so is $\omega = W \eta$.
\end{cor}

\section{\texorpdfstring{$\sW$}{Wbar}-matrix}

There is an `averaged' version $\sW$ of the $W$-matrix that will be useful in our analysis of dual codewords in Section~\ref{sec:DualCodewords}.  We describe $\sW$ in this section.  We continue to assume that the information module is $M = M_{k \times m}(\F_q)$ and that $w_0$ is indeterminate.   

\begin{lem}  \label{lem:sumDependsOnlyOnj}
For $i=0, 1, \ldots, k$ and $[\la] \in \OS_j$, the value of
\[  \sum_{[x] \in \OO_i} W_{[x],[\la]} 
= \sum_{d=0}^i q^{(i-d)(m-j-d)} \qbinom{m-j}{d} \qbinom{j}{i-d} w_{i-d} \]
depends only on $i$ and $j=\rk \la$.
\end{lem}

\begin{proof}
The proof of Proposition~\ref{prop:loc-constant-eta} applies, interchanging the roles of $i=\rk x$ and $j = \rk \la$.
\end{proof}

\begin{defn}
Define $\sW$ to be the integer $(k+1)\times (k+1)$ matrix 
\begin{equation}  \label{eq:Defn-sW}
\sW_{ij} = \sum_{[x] \in \OO_i} W_{[x],[\la]}, \quad [\la] \in \OS_j,
\end{equation}
for $i,j=0, 1, \ldots, k$.  The definition is well-defined by Lemma~\ref{lem:sumDependsOnlyOnj}.
\end{defn}

\begin{ex}  \label{ex:sWfor2}
When $k=2$ and $m=3$, we see that
\[ \sW = \begin{bmatrix}
w_0 & w_0 & w_0 \\
(1+q+q^2)w_0 & (1+q) w_0 + q^2 w_1 & w_0 + (q + q^2) w_1 \\
(1+q+q^2)w_0 & w_0 + (q + q^2) w_1 & (1+q) w_1 + q^2 w_2
\end{bmatrix} . \]
\end{ex}

We next formalize the relationship between $\sW$ and $W$, in order to determine when $\sW$ is invertible.

Recall that $W$ is square of size $N_{k,m}$, the number of linear subspaces of dimension at most $k$ in $\F_q^m$.
Define $B$ to be a $(k+1) \times N_{k,m}$ matrix.
For $i=0,1, \ldots, k$, row $i$ of $B$ is the indicator function $1_{\OO_i}$ for the collection of linear subspaces of dimension $i$ in $\F_q^m$.  Similarly, define an $N_{k,m} \times (k+1)$ matrix $E$.  For $j=0, 1, \ldots, k$, column $j$ of $E$ is $(1/\smqbinom{m}{j}) 1_{\OS_j}$.  Notice that $BE = I_{k+1}$, and $EB$ is block diagonal, with the block indexed by rank $j$ being $1/\smqbinom{m}{j}$ times the square all-one matrix of size $\smqbinom{m}{j}$.

\begin{lem}
We have: $\sW = B W E$.
\end{lem}

\begin{proof}
Left multiplication by $B$ givens the sums of Lemma~\ref{lem:sumDependsOnlyOnj}.  Right multiplication by $E$ averages those (equal!) sums over $[\la] \in \OS_j$.
\end{proof}

Recall the $P$-matrix of \eqref{eq:DefnOfPmatrix}.

\begin{lem}  \label{lem:BE-simplification}
We have: $BPEB = BP$ and $EBP^\top E = P^\top E$.
\end{lem}

\begin{proof}
The matrix $BP$ has rows that are locally constant, and $EB$ acts as the identity when it right multiplies matrices with rows that are locally constant.  Similarly, $P^\top E$ has columns that are locally constant, and $EB$ acts as the identity when it left multiplies matrices with columns that are locally constant.
\end{proof}

Our next result is the counterpart of Theorem~\ref{thm:W-block-diagonal-form}.
Set $Q_1 = BPE$ and $Q_2 = B P^\top E$.  One verifies that $Q_1$ is lower triangular and $Q_2$ is upper triangular.  Their 
$i,j$-entries are
\begin{equation}  \label{eq:Q-matrices}
(Q_1)_{i,j} = (-1)^j q^{\binom{j}{2}} \qbinom{m-j}{i-j} , \quad (Q_2)_{i,j} = (-1)^i q^{\binom{i}{2}} \qbinom{j}{i}  .
\end{equation}
In particular, the diagonal entries of both are
$(-1)^j q^{\binom{j}{2}}$, $j=0, 1, \ldots, k$, so that both $Q_1$ and $Q_2$ are invertible.  

\begin{thm}
We have $Q_1 \sW Q_2 = B P W P^\top E$ and 
\[  Q_1 \sW Q_2  = \left [ 
\begin{array}{cccccc}
c_0 & 0 & \dots & 0 & \dots & 0 \\
0 & q^{m-1} c_1 & \dots & 0 & \dots & 0  \\
\vdots & \vdots & \ddots & \vdots & & \vdots\\
0 & 0 & \dots & q^{j(m-j)} c_j & \dots & 0 \\
\vdots & \vdots & & \vdots & \ddots & \vdots\\
0 & 0 & \dots & 0 & \dots & q^{k(m-k)} c_k
\end{array} \right ] , \]
where $c_0, c_1, \ldots, c_k$ are in \eqref{eqn:diag-coefficents}.  
\end{thm}

\begin{proof}
The first equation follows from Lemma~\ref{lem:BE-simplification}.  The second equation follows from Proposition~\ref{prop:W-block-diagonal}.  The factor $q^{j(m-j)}$ arises from counting the number of $\al$ of dimension $j$ in $\F_q^m$ that satisfy $\al \cap \delta^\perp = 0$ for a fixed $\delta$ of dimension $j$.  That count uses Lemma~\ref{lem:count-hit-trivially}.
\end{proof}

\begin{cor}
The matrix $W$ is invertible if and only if $\sW$ is invertible.
\end{cor}

When $w_0=0$ one can define a $k \times k$ matrix $\sW_0$ using \eqref{eq:Defn-sW}, but only for $i,j=1,2, \ldots, k$.  By defining smaller versions of $B, E, Q_1, Q_2$, one can prove the following results using a proof similar to that of Theorem~\ref{thm:W0-block-diagonal-form}.

\begin{thm}  \label{thm:Diagonalize-sW}
Suppose $w_0=0$.  Then $Q_{0,1} \sW_0 Q_{0,2} = B_0 P_0 W_0 P_0^\top E_0$ and
\[  Q_{0,1} \sW_0 Q_{0,2} = \left [ 
\begin{array}{ccccc}
q^{m-1} c_1 & \dots & 0 & \dots & 0  \\
\vdots & \ddots & \vdots & & \vdots\\
0 & \dots & q^{j(m-j)} c_j & \dots & 0 \\
\vdots & & \vdots & \ddots & \vdots\\
0 & \dots & 0 & \dots & q^{k(m-k)} c_k
\end{array} \right ] , \]
where $c_1, \ldots, c_k$ are in \eqref{eqn:diag-coefficents}, but with $w_0=0$.  
\end{thm}

\begin{cor}
When $w_0=0$, the matrix $W_0$ is invertible if and only if $\sW_0$ is invertible.
\end{cor}

\section{Constructions}  \label{sec:construction}

In this section we construct two linear codes $C$ and $D$ over $R=M_{k \times k}(\F_q)$ with $\wwe_C = \wwe_D$, assuming that $w$ has maximal symmetry, that $w_0=0$, and that the associated $W_0$ matrix is invertible, Corollary~\ref{cor:whenisW0invertible}.  Both linear codes will have information module $M=M_{k \times m}(\F_q)$, with $m>k$.

Here is a sketch of the main idea behind the construction.   Suppose $C$ 
is the image of $\La: M \ra R^n$.  Recall that the maximal symmetry hypothesis means that $W_{\La}$ is constant on any left $\U$-orbit in the information module $M$. 
Suppose $\La$ has the property that $N$ chosen orbits of rank $s$ have the same value $v_1$ of $W_{\La}$.  Pick one orbit of rank $s+1$, and denote by $v_2$ its value of $W_{\La}$.  Now try to swap values on those orbits: try to find a linear code $D$, the image of $\Gamma: M \ra R^n$, so that $W_{\Gamma}$ has value $v_2$ on the $N$ chosen orbits of rank $s$, value $v_1$ on the chosen orbit of rank $s+1$, and $W_{\Gamma}= W_{\La}$ on all other orbits.  If $N$ times the size of a rank $s$ orbit equals the size of a rank $s+1$ orbit, then $\wwe_C = \wwe_D$, provided such a $D$ exists
and provided that $C$ and $D$ have the same length; cf., Remark~\ref{rem:appendZeroFuncl}.  If the codes have different lengths, we append enough zero-functionals to the shorter code so that the lengths become equal.  Because $w_0=0$, the additional zero-functionals have no effect on the weights.

We first need a few facts about $q$-binomial coefficients.

\begin{lem}  \label{lem:q-binomialCounts}
For integers $0 \leq s \leq m$, 
\begin{align*}
\qbinom{m}{s} &= \qbinom{m}{m-s} ; \quad
\qbinom{m}{s+1} = \frac{q^{m-s}-1}{q^{s+1}-1} \qbinom{m}{s} ; \\
\qbinom{m}{s} &\leq \qbinom{m}{s+1}, \quad 0 \leq s < m/2 ; \\
\qbinom{m}{s} &= q^{m-s} \qbinom{m-1}{s-1} + \qbinom{m-1}{s}, \quad 0 \leq s <m ; \\
\qbinom{m}{s} &= \qbinom{m-1}{s-1} + q^s \qbinom{m-1}{s}, \quad 0 \leq s <m ;
\end{align*}
\end{lem}

\begin{proof}
Most of the identities are in \cite[Chapter~3]{MR4249619} or its exercises.  The inequality follows from the preceding identity, because the multiplying fraction is at least $1$ for $s<m/2$.
\end{proof}

\begin{lem}  \label{lem:enoughSubspaces}
Suppose integers $k,m,s$ satisfy $0<s<k<m$.  Then
\[  q^k \leq \qbinom{m}{s} - \qbinom{m-1}{s} . \]
\end{lem}

\begin{proof}
The result is true for $s=1$: $\smqbinom{m}{1} - \smqbinom{m-1}{1} = q^{m-1} \geq q^k$.  Now suppose $s \geq 2$, so that $2 \leq s \leq m-2$.  By symmetry and monotonicity in Lemma~\ref{lem:q-binomialCounts}, we see that $\smqbinom{m-1}{s-1} \geq \smqbinom{m-1}{1} = q^{m-2} + q^{m-3} +\cdots +q + 1 \geq q^{m-2}$.  Use Lemma~\ref{lem:q-binomialCounts} again to see that 
\[  \qbinom{m}{s} - \qbinom{m-1}{s} = q^{m-s} \qbinom{m-1}{s-1} \geq q^2 q^{m-2} = q^m > q^k .  \qedhere \]
\end{proof}

For the constructions, assume $k \geq 2$, $m > k$, and set $w_0=0$.  We assume a weight $w$ on $R$ has maximal symmetry, and we assume the associated $W_0$-matrix is invertible.  There will be a construction for each value of $s = 1, 2, \ldots, k-1$.  The integer $s$ determines the ranks of the orbits that will be swapped.

Fix an integer $s$ with $1 \leq s < k<m$, and
choose an orbit $\orb(\la_0) \in \OS$ with $\rk \la_0=1$.  As in Corollary~\ref{cor:orbit-correspondence}, this orbit corresponds to the linear subspace $L_0 = \colsp(\la_0) \subseteq \F_q^m$, with $\dim L_0 = 1$.  As $\dim L_0^\perp = m-1$, there are $\smqbinom{m-1}{s}$ linear subspaces of dimension $s$ contained in $L_0^\perp$, leaving $\smqbinom{m}{s} - \smqbinom{m-1}{s}$ linear subspaces of dimension $s$ in $\F_q^m$ that are not contained in $L_0^\perp$.  
By Lemma~\ref{lem:enoughSubspaces}, there are at least $q^k - q^s < q^k$ linear subspaces of dimension $s$ in $\F_q^m$ that are not contained in $L_0^\perp$.  Similarly, using $m-s \geq 2$, there are $\smqbinom{m-1}{s+1} \geq 1$ linear subspaces of dimension $s+1$ contained in $L_0^\perp$.

Choose distinct $s$-dimensional linear subspaces $X_1, X_2, \ldots, X_{q^k-q^s}$ of $\F_q^m$ that are not contained in $L_0^\perp$.  Also choose one $(s+1)$-dimensional linear subspace $Y \subseteq L_0^\perp$.  By Corollary~\ref{cor:orbit-correspondence}, these linear subspaces correspond to distinct orbits $[x_1]=\orb(x_1), \ldots, [x_{q^k-q^s}] =\orb(x_{q^k-q^s})$ and $[y]=\orb(y)$ in $\OO$, with $\rk x_i = s$ and $\rk y=s+1$.  
Then $\size{\orb(y)} = \sum_{i=1}^{q^k-q^s} \size{\orb(x_i)}$, by Proposition~\ref{prop:orbitNumberSize}.

We now consider several indicator functions on $\OS$ and find the weights they determine at the orbits $[x_1], \ldots, [x_{q^k-q^s}]$ and $[y]$.

Let $1_{[\la_0]}: \OS \ra \Z$ be the indicator function of $[\la_0] = \orb(\la_0) \in \OS$: 
\[  1_{[\la_0]}([\la]) = \begin{cases}
1, & [\la] = [\la_0], \\
0, &  [\la] \neq [\la_0]. 
\end{cases} \]
From Lemma~\ref{lem:dim-formulas}, the weights of $1_{[\la_0]}$ at the orbits $[x_i]$, $1 \leq i \leq q^k-q^s$, and $[y]$ are $(W_0 1_{[\la_0]})([x]) 
= w(x \la_0) = w_{\dim X - \dim(X \cap L_0^\perp)}$: 
\[  (W_0 1_{[\la_0]})([x]) = \begin{cases}
w_1, & [x] = [x_i], \\
0, & [x] = [y] .
\end{cases} \]
The exact values of $W_0 1_{[\la_0]}$ at other inputs will not be relevant.
What is crucial is that the value at $[y]$ is $0$ and that the values at the $[x_i]$ are equal and positive.

Recall that $\OS_1 = \{ \orb(\la) \in \OS: \rk \la =1\}$; let $1_{\OS_1}$ be its indicator function.  The orbit weights $\omega_1 = W_0 1_{\OS_1}$ are found in Proposition~\ref{prop:loc-constant-eta}.
The indicator function $1_{\OS_+}$ of the set of all nonzero orbits $\OS_+$ is locally constant, so $W_0 1_{\OS_+}$ is also locally constant, Corollary~\ref{cor:LocallyConstant}.  Set $\al_1 = (W_0 1_{\OS_+})(\orb(x))$, when $\rk x =1$, and $\al_2 = (W_0 1_{\OS_+})(\orb(x))$ when $\rk x =2$.  Both $\al_1, \al_2$ are positive integers.

Define a function $\varsigma^{(s)} \in F_0(\OO, \Z)$ by 
\[  \varsigma^{(s)}(\orb(x)) = \begin{cases}
-1, & \orb(x) = \orb(x_i), \\
1, & \orb(x) = \orb(y), \\
0, & \text{otherwise}.
\end{cases} \]
Because we are assuming $W_0$ is invertible, Remark~\ref{rem:EPand W0} says that $W_0: F_0(\OS, \Q) \ra F_0(\OO, \Q)$ is an isomorphism.  Thus $W_0^{-1} \varsigma^{(s)} \in F_0( \OS , \Q)$ exists, but it has rational values of both signs.  To clear denominators, choose a positive integer $c$ sufficiently large so that $\sigma^{(s)}=c w_1 W_0^{-1} \varsigma^{(s)}$ has integer values.  Then $\sigma^{(s)}(\orb(0)) = 0$, as $\sigma^{(s)} \in F_0(\OS, \Z)$, and $W_0 \sigma^{(s)} = c w_1 \varsigma^{(s)}$, so that $(W_0 \sigma^{(s)})(\orb(x))$ equals $0$ or $\pm c w_1$.

We collect the values of $W_0 \eta$ at the $[x_i]$ and $[y]$, for various $\eta$.
\begin{equation}
\begin{array}{c|c|c}  \label{eq:weightsForVariousEtas}
\eta & (W_0 \eta)([x_i]) & (W_0 \eta)([y]) \\ \hline
1_{[\la_0]} & w_1 & 0 \\
1_{\OS_+} & \al_1 & \al_2 \\
1_{\OS_1} & (q^{m-s} + \cdots + q^{m-1}) w_1 & (q^{m-s-1} + \cdots + q^{m-1}) w_1 \\
\sigma^{(s)} & -c w_1 & c w_1 
\end{array}
\end{equation}
For later use, write $B_{m,s} = q^{m-s-1} + \cdots + q^{m-1}$; $B_{m,s}>0$.

Some of the values of $\sigma^{(s)}$ will be negative.  Choose a positive integer $a$ sufficiently large so that all the values of $a w_1 1_{\OS_+} + \sigma^{(s)}$ are nonnegative.  
Choose an integer $b \geq 1$ large enough that $a(\al_1 - \al_2) < b q^{m-s-1}$.  (If $\al_1 < \al_2$, then $b=1$ suffices.)   Set
\begin{align} \label{eq:DefnDelta}
 \Delta &= \sum_{\orb(\la) \in \OS_+} \sigma^{(s)}(\orb(\la)), \\
 \varepsilon &= a w_1 1_{\OS_+} + b 1_{\OS_1} + (c + a(\al_2 - \al_1) + b q^{m-s-1}) 1_{[\la_0]}. \notag
\end{align} 
Define two multiplicity functions $\eta_C, \eta_D \in F(\OS,\Q)$ by setting
\begin{align}
\eta_C &= \varepsilon + \max(\Delta,0) 1_{[0]} , \notag \\
\eta_D &=  \varepsilon + \sigma^{(s)} -\min(\Delta,0) 1_{[0]} .  \label{eq:defnCD}
\end{align}

\begin{thm}  \label{thm:MainConstruction}
Let $R = M_{k \times k}(\F_q)$, $k \geq 2$, and $M= M_{k \times m}(\F_q)$, $k<m$.  Let $w$ be a weight on $R$ with maximal symmetry, positive integer values, and $w(0)=0$.  Assume the associated $W_0$-matrix is invertible.

Then, for any integer $s$, $1 \leq s <k$, the multiplicity functions $\eta_C$ and $\eta_D$ of \eqref{eq:defnCD} have nonnegative integer values, 
i.e., $\eta_C, \eta_D \in F(\OS,\N)$,  and they define left $R$-linear codes $C$ and $D$, respectively.  
The two codes have the same length.  Their weights at the orbits $\orb(x_i)$, $i=1, 2, \ldots, q^k-q^s$, and $\orb(y)$ are 
\[  \begin{array}{c|c|c} 
\eta & (W \eta)(\orb(x_i)) & (W \eta)(\orb(y)) \\ \hline
\eta_C & (c + a \al_2 + b B_{m,s}) w_1  & (a \al_2 + b B_{m,s})w_1 \\
\eta_D & (a \al_2 + b B_{m,s})w_1 &  (c + a \al_2 + b B_{m,s}) w_1
\end{array}.  \]
At all other orbits, their weights agree.  In particular, $\wwe_C = \wwe_D$.
\end{thm}

\begin{proof}
All of $a, w_1, b, c, \al_1, \al_2, \max(\Delta,0)$, and $- \min(\Delta,0)$ are nonnegative integers.  The integer $b \geq 1$ was chosen so that $a(\al_2-\al_1) + b q^{m-s-1}>0$, so the values of $\eta_C$ are nonnegative integers.  
The positive integer $a$ was chosen so that $a w_1 1_{\OS_+} + \sigma^{(s)}$ has nonnegative values, which implies that the values of $\eta_D$ are nonnegative integers.  

Both $\eta_C$ and $\eta_D$ contain the term $b 1_{\OS_1}$.
Because the linear subspace spanned by $\{ \colsp(\la): \la \in \OS_1 \}$ is all of $\F_q^m$, Remark~\ref{rem:functionalEvaluation} implies that $\La_C$ and $\La_D$ are both injective.  The value of $\Delta$ was chosen to satisfy $\Delta = \efflength(D)-\efflength(C)$.  If $\Delta >0$, then $D$ is longer, and we add $\Delta$ zero-functionals to $C$; if $\Delta <0$, then $C$ is longer, and we add $-\Delta$ zero-functionals to $D$.  Thus the codes have the same length.

The weights $(W \eta)(\orb(x_i))$ and $(W \eta)(\orb(y))$ follow from \eqref{eq:weightsForVariousEtas} and \eqref{eq:defnCD}.  Because $W \eta_D - W \eta_C = W \sigma^{(s)} = c w_1 \varsigma^{(s)}$, the form of $\varsigma^{(s)}$ implies that 
the weights for $C$ and $D$ are equal at all other orbits.  As noted earlier, $\size{\orb(y)} = \sum_{i=1}^{q^k-q^s} \size{\orb(x_i)}$, by Proposition~\ref{prop:orbitNumberSize}.  We conclude that $\wwe_C = \wwe_D$.
\end{proof}

\begin{rem}  \label{rem:reflectionOnConstruction}
Here is a summary of the motivations for various parts of the construction.
The integer $c$ was chosen to clear denominators so that $\sigma^{(s)}$ would have integer values.  However, $\sigma^{(s)}$ has both positive and negative values, because $W$ has nonnegative entries and $W \sigma^{(s)} = c w_1 \varsigma^{(s)}$ has mixed signs.  
So, $a$ was chosen so that $a w_1 1_{\OS_+} + \sigma^{(s)}$ has nonnegative integer values.  The function $1_{\OS_+}$ (resp., $1_{\OS_1}$) is used because its weights are the same at all the $[x_i]$'s.  The integer $b \geq 1$ was chosen so that $(W_0(a w_1 1_{\OS_+} + b 1_{\OS_1}))(\orb(x_i)) < (W_0(a w_1 1_{\OS_+} + b 1_{\OS_1}))(\orb(y))$ and also to guarantee that $\La_C$ and $\La_D$ are injective.  Then the coefficient of $1_{[\la_0]}$ was chosen so that the weights interchange when $\sigma^{(s)}$ is added to $\eta_C$.  The function $1_{[\la_0]}$ is used because it changes the weights at the $[x_i]$'s in the same way but not the weight at $[y]$.
\end{rem}

\section{Degeneracies}  \label{sec:degeneracy}
The constructions of Theorem~\ref{thm:MainConstruction} made use of the invertibility of the matrix $W_0$.  In this section, we describe a construction in the situation where $W_0$ is not invertible.  We continue to assume information module $M = M_{k \times m}(\F_q)$ and weight $w$ with maximal symmetry, positive integer values, and $w_0=0$. 

Corollary~\ref{cor:whenisW0invertible} says that $W_0$ is singular when at least one
\[  c_j = (-1)^j q^{\binom{j}{2}} \sum_{s=1}^j (-1)^s q^{\binom{s}{2}} \qbinom{j}{s} w_s , \]
$j=1, 2, \ldots, k$, vanishes.  By hypothesis on $w$, $w_1>0$, so that $c_1 \neq 0$.  

\begin{lem}
Suppose $c_j=0$ for some $j=2, 3, \ldots, k$.  Then any row of the matrix $P_0$ indexed by a linear subspace $\ga \subseteq \F_q^m$ with $\dim \ga = j$ belongs to $\ker W_0$.
\end{lem}

\begin{proof}
The second part of Proposition~\ref{prop:PW-format} shows that, if $\dim \ga = j$, then $(W P^\top )_{\be, \ga} =0$.  This holds for any $\be$ and any $\ga$ with $\dim\ga = j$.  This means the columns of $P^\top$ indexed by $\ga$ with $\dim \ga = j$ belong to $\ker W$.  Those columns are the same as the rows of $P$ indexed by $\ga$ with $\dim\ga = j$.  Because $w_0=0$, the same analysis applies with $W_0$ and $P_0$.
\end{proof}

Suppose $c_j=0$ for some $j=2, 3, \ldots, k$.  Pick a row $v_{\ga}$ of $P_0$ indexed by a linear subspace $\ga \subseteq \F_q^m$ with $\dim\ga = j$.  Define two multiplicity functions $\eta_{\pm} \in F_0(\OS, \N)$ based on the positive, resp., negative, parts of $v_{\ga}$:
\begin{align*}
\eta_+ &= \hphantom{-}(v_{\ga} + \abs{v_{\ga}})/2 + 1_{\OS_1}, \\
\eta_- &= -(v_{\ga} - \abs{v_{\ga}})/2 + 1_{\OS_1},
\end{align*}
where $\abs{v_{\ga}}$ means the vector obtained from $v_{\ga}$ by taking the absolute value of each entry.  The terms $1_{\OS_1}$ are included so that the associated homomorphisms $\La_{\eta_\pm}$ are injective; see  Remark~\ref{rem:functionalEvaluation}.  Note that $\eta_+ - \eta_- = v_{\ga}$, so that $W_0 \eta_+ = W_0 \eta_-$.  Modify $\eta_+$ by setting $\eta_+([0]) = 1$; set $\eta_-([0])=0$.

\begin{prop}  \label{prop:DegenerateCase}
Let $C_{\pm}$ be the linear codes determined by $\eta_{\pm}$.  Then $C_+$ and $C_-$ have the same length, and $\wwe_{C_+} = \wwe_{C_-}$.
\end{prop}

\begin{proof}
In \eqref{eq:DefnOfPmatrix}, there are $\smqbinom{j}{r}$ nonzero terms in rank $r$ positions of the row $v_{\ga}$.  Then the difference in lengths of $C_+$ and $C_-$ is 
\begin{align*}
\length(C_+)  - \length(C_-) &=
\sum_{\substack{r=0 \\ \text{$r$ even}}}^j q^{\binom{r}{2}} \qbinom{j}{r} - \sum_{\substack{r=0 \\ \text{$r$ odd}}}^j q^{\binom{r}{2}} \qbinom{j}{r} \\
&= \sum_{r=0}^j (-1)^r q^{\binom{r}{2}} \qbinom{j}{r} =0,
\end{align*} 
by 
\eqref{eq:Cauchy-vanishing}.  Zero-functionals don't change the value of $W \eta_+$.  We then have $W \eta_+ = W \eta_-$, so that the $w$-weight enumerators are equal.
\end{proof}

This same construction was used (with the Hamming weight) in \cite[p.\ 703]{wood:code-equiv} and \cite[p.\ 145]{wood:turkey}.

\begin{rem}  \label{rem:NoSwapDegenerate}
The swapping idea in Theorem~\ref{thm:MainConstruction} does not always work in the degenerate case because $\varsigma^{(s)}$ is not always in $\im W_0$.  See Remark~\ref{rem:DegenerateImageW0} for more details.
\end{rem}

\section{Analysis of singleton dual codewords}  \label{sec:DualCodewords}
In this section we analyze dual codewords that are singleton vectors.  We will then apply this analysis 
to the codes constructed using \eqref{eq:defnCD}.  The key result, Theorem~\ref{thm:SingletonContributions}, shows how the contributions of singletons of rank $i$ to $\asing_{w_i}(D^\perp) - \asing_{w_i}(C^\perp)$ depend on the parameter $s$ used in \eqref{eq:defnCD}.

Given a left $R$-linear code $C \subseteq R^n$, recall the right dual code $\mathcal{R}(C)$ from \eqref{eq:RightDualCode}.  We will often denote $\mathcal{R}(C)$ by $C^\perp$.

When an $R$-linear code $C$ is given by a multiplicity function $\eta$, a singleton vector $v$ belongs to $C^\perp$ when the nonzero entry $r$ of $v$ right-annihilates the coordinate functional $\la$ in that position: $\la r = 0$.  Remember that $\la \in \Hom_R(M,R) = M_{m \times k}(\F_q)$ and $r \in R = M_{k \times k}(\F_q)$.  

Given a functional $\la$ with $\rk \la = j$, we will determine how many elements $r \in R$ with $\rk r = i$ satisfy $\la r = 0$.  
For $\la \in \Hom_R(M,R)$, define 
\[  \ann(i,\la) = \{ r \in R: \rk r = i \text{ and } \la r =0 \}. \]
Recall that the sizes $\mathscr{S}_j$ of $\U$-orbits were given in Proposition~\ref{prop:orbitNumberSize}.

\begin{lem}  \label{lem:CountAnnihilatorsRanki}
Suppose $\la \in \Hom_R(M,R)$.   Then the size of $ \ann(i,\la)$ depends only on $\rk \la$.  If $\rk \la = j$, then
\[  \size{\ann(i, \la)} = \begin{cases}
\mathscr{S}_i \smqbinom{k-j}{i} , & i \leq k-j, \\
0, & i > k-j .
\end{cases}  \]
\end{lem}

\begin{proof}
View $M_{m \times k}(\F_q)$ and $M_{k \times k}(\F_q)$ as spaces of $\F_q$-linear transformations $\F_q^m \ra \F_q^k$ and $\F_q^k \ra \F_q^k$, respectively, with inputs written on the left.  Then $\la r=0$ means that $\im \la \subseteq \ker r$.  Given that $\rk \la = j$ and $\rk r = i$, we see that a necessary condition for $\la r =0$ is that $j = \dim\im\la \leq \dim\ker r = k-i$.

Given $\im \la \subseteq \F_q^k$, the number of linear subspaces $K$ (candidates for $\ker r$) satisfying $\im \la \subseteq K \subseteq \F_q^k$ and $\dim K = k-i$ equals the number of linear subspaces of dimension $k-i-j$ in $\F_q^k/\im\la$, a vector space of dimension $k-j$.  That number is $\smqbinom{k-j}{k-i-j}=\smqbinom{k-j}{i}$.  For a given $K$, there is a $\U$-orbit's worth of matrices $r$ with $\ker r = K$, and hence $\rk r =i$.  The size of that orbit is $\mathscr{S}_i$.
\end{proof}

We will abuse notation and write $\size{\ann(i,j)}$ for the common value $\size{\ann(i,\la)}$ when $\rk \la = j$.

\begin{ex}  \label{ex:AnnMatrixk2}
Let $k=2$.  We display the values of $\size{\ann(i,j)}$ in a matrix with row index $i$ and column index $j$, $i,j=0,1,2$:
\[ \begin{bmatrix}
1 & 1 & 1 \\
(q^2-1)(q+1) & q^2-1 & 0 \\
(q^2-1)(q^2-q) & 0 & 0
\end{bmatrix} . \]
\end{ex}

Suppose a left $R$-linear code $C$ is given by a multiplicity function $\eta: \OS \ra \N$.  For $j=0,1, \ldots, k$, define 
\[   \eb_j =  \sum_{[\la] \in \OS_j} \eta([\la])  . \]
Then $\eb_j$ counts the number of coordinate functionals having rank $j$.  Set $\eb = \langle \eb_0, \eb_1, \ldots, \eb_k \rangle \in \N^{k+1}$.
We call the $\eb_j$ `rank-sums'. 

\begin{prop}  \label{prop:SingletonContribs}
Suppose a left $R$-linear code $C$ is given by a multiplicity function $\eta$.  Then the contribution of singletons of rank $i$ to $\asing_{w_i}(C^\perp)$ is
\[  \sum_{j=0}^{k-i} \size{\ann(i,j)}  \eb_j . \]
In particular, the contribution of singletons of rank $k$ to $\asing_{w_k}(C^\perp)$ is $\size{\GL(k,\F_q)} \eb_0$.
\end{prop}

\begin{proof}
There are a total of $\eb_j$ coordinate functionals of rank $j$.  For each one, apply Lemma~\ref{lem:CountAnnihilatorsRanki}.
\end{proof}

The larger $i$ is, the fewer terms there are in the summation.

\begin{cor}  \label{cor:contribsBySingletons}
For a linear code $C$ given by multiplicity function $\eta$,
\[  \asing_d(C^\perp) = \sum_{i: w_i=d} \sum_{j=0}^{k-i} \size{\ann(i,j)} \eb_j . \]
\end{cor}

Having seen the importance of the rank-sums $\eb_j$, we next see how they behave with respect to the $\sW$-matrix.  Given a function $\omega: \OO \ra \Q$, define
\[  \ob_i = \sum_{[x] \in \OO_i} \omega([x]), \quad i=0,1, \ldots, k. \]
Set $\ob = \langle \ob_0, \ob_1, \ldots, \ob_k \rangle$.

\begin{prop}
Suppose $\eta: \OS \ra \Q$.  If $\omega = W \eta$, then $\ob = \sW \eb$.
\end{prop}

\begin{proof}
Sum the rows indexed by rank $i$ elements of $\OO$ in $\omega = W \eta$, change the order of summation, and use Lemma~\ref{lem:sumDependsOnlyOnj} and \eqref{eq:Defn-sW}.
\end{proof}

Because our ultimate objective is to find linear codes $C$ and $D$ of the same length with $\wwe_C = \wwe_D$ and $\wwe_{C^\perp} \neq \wwe_{D^\perp}$, we now apply Proposition~\ref{prop:SingletonContribs} and Corollary~\ref{cor:contribsBySingletons} to two linear codes $C$ and $D$ of the same length.  Write $\delta \eta = \eta_D - \eta_C$ and $\delta \asing_d = \asing_d(D^\perp) - \asing_d(C^\perp)$.  The net contribution of rank $i$ singletons to $\delta \asing_{w_i}$ then simplifies to 
\begin{equation}  \label{eq:netContribOfRankiSingletons}
  \sum_{j=0}^{k-i} \size{\ann(i,j)}  \delta\eb_j ,
\end{equation}
and $\delta \asing_d = \sum_{i: w_i=d} \sum_{j=0}^{k-i} \size{\ann(i,j)} \delta\eb_j$.

\begin{rem}
The last two formulas will be the main tools for showing that two dual codes have different $w$-weight enumerators.  But the fomulas cut both ways.  If all the $\delta\eb_j$ vanish, then singletons cannot detect differences between $\wwe_C$ and $\wwe_D$.  See Example~\ref{ex:linIndepDep}.
\end{rem}

Now consider specifically the linear codes $C$ and $D$ constructed by \eqref{eq:defnCD}. In this case, we see that 
\begin{align*}
\delta\eta &= \sigma^{(s)} - \Delta 1_{[0]} , \\
\eb_D - \eb_C &= \langle -\Delta, \os_1^{(s)}, \ldots, \os_k^{(s)} \rangle , \\
\ob_D - \ob_C &= \sW \os^{(s)} =  \langle 0, \ldots, 0, -cw_1 (q^k-q^s), c w_1, 0, \ldots, 0 \rangle ,
\end{align*}
where the nonzero entries of $\ob_D - \ob_C$ are in positions $s$ and $s+1$.

Recall from Theorem~\ref{thm:MainConstruction} that $C$ and $D$ have the same length.  This
is reflected in the fact that the sum of the entries of $\eb_D - \eb_C$ is $-\Delta + \sum_j \os_j^{(s)} = 0$, from \eqref{eq:DefnDelta}.  This allows us to re-write 
the net contribution \eqref{eq:netContribOfRankiSingletons} of rank $i$ singletons to $\delta \asing_{w_i}$ as $\sum_{j=1}^k \left( \size{\ann(i,j)} - \size{\ann(i,0)} \right) \os_j^{(s)}$.  To write this equation in matrix form, define a $k \times k$ matrix $\sAnn$ by
\[  \sAnn_{i,j} = \size{\ann(i,j)} - \size{\ann(i,0)} \quad i,j=1,2, \ldots, k. \]
In summary, the net contributions of singletons of rank $i$ to $\delta \asing_{w_i}$ are given by the entries of
\begin{equation}  \label{eq:SingletonContribs}
\sAnn \os^{(s)} .
\end{equation}

\begin{lem}
The matrix $\sAnn$ is invertible over $\Q$.
\end{lem}

\begin{proof}
Define a $(k+1) \times (k+1)$ matrix $\A$ by
\[  \A_{i,j} = \size{\ann(i,j)}, \quad i,j=0,1,2, \ldots, k. \]
Then $\mathcal{A}$ is upper `anti-triangular' by Lemma~\ref{lem:CountAnnihilatorsRanki}, i.e., $\A_{i,j} =0$ when $i+j>k$.  Then $\det \A = \pm \prod_{i=0}^k \mathscr{S}_i$ is nonzero.  Thus $\A$ is invertible.

Subtract the zeroth column of $\A$ from every other column.  The resulting matrix is still invertible and has the form
\[  \left[ \begin{array}{c|c}
1 & \row(0) \\ \hline
\col(*) & \vphantom{\ah} \sAnn
\end{array} \right]  ,\]
where $\col(*)$ is a column of nonzero entries.  Thus $\sAnn$ is invertible.
\end{proof}

The matrices $\sW_0$ and $\sAnn$ are both $k \times k$, so they both define $\Q$-linear transformations $\Q^k \ra \Q^k$.  We next explore how these linear transformations behave with respect to certain filtrations of $\Q^k$.

Define a $k \times k$ matrix $T$ over $\Q$:
\[  T_{i,j} = (-1)^i q^{\binom{i}{2}} \smqbinom{j}{i} , \quad i,j=1, 2, \ldots, k . \]
By standard conventions, $\smqbinom{j}{i} = 0$ when $i>j$.  This implies that $T$ is upper-triangular and invertible.  (By \eqref{eq:Q-matrices}, $T$ is just $Q_{0,2}$.)

\begin{lem}  \label{lem:AnnT-lowerTriangular}
The matrices $\sW_0 T$ and $\sAnn T$ are lower triangular.
\end{lem}

The proof of Lemma~\ref{lem:AnnT-lowerTriangular} will utilize the next lemma.

\begin{lem}  \label{lem:NewSumVanishing}
For integers $0 \leq i < j \leq k$, 
\[  \sum_{\ell = 0}^j (-1)^{\ell} q^{\binom{\ell}{2}} \qbinom{j}{\ell} \qbinom{k-\ell}{i} = 0 . \]
\end{lem}

\begin{proof}
We first prove the edge cases $i=0$ or $j=k$, and then prove the remaining cases by induction on $k$.  

Suppose $i=0$.  Then the term $\smqbinom{k-\ell}{0} = 1$ for all $\ell = 0 , \ldots j$.  The sum reduces to \eqref{eq:Cauchy-vanishing}, which vanishes, as $j>0$.

Suppose $j=k$.  We observe that 
\[  \qbinom{k}{\ell} \qbinom{k-\ell}{i} = \qbinom{k}{i} \qbinom{k-i}{\ell} , \]
and that $\smqbinom{k-\ell}{i} = 0$ for $\ell>k-i$.  Then
\begin{align*}
\sum_{\ell = 0}^k (-1)^{\ell} q^{\binom{\ell}{2}} \qbinom{k}{\ell} \qbinom{k-\ell}{i} &= \sum_{\ell = 0}^{k-i} (-1)^{\ell} q^{\binom{\ell}{2}} \qbinom{k}{i} \qbinom{k-i}{\ell} \\
&= \qbinom{k}{i}  \sum_{\ell = 0}^{k-i} (-1)^{\ell} q^{\binom{\ell}{2}} \qbinom{k-i}{\ell} =0,
\end{align*}
by \eqref{eq:Cauchy-vanishing}, as $k-i>0$.

We prove the remaining cases, $0<i<j<k$, by induction on $k$.  The first case is when $k=3$, with $i=1<j=2$.  By direct calculation,
\begin{align*}
\sum_{\ell = 0}^2 (-1)^{\ell} q^{\binom{\ell}{2}} \qbinom{2}{\ell} \qbinom{3-\ell}{1} 
&= (q^2+q+1) - (q+1)^2 +q = 0.
\end{align*}  

For the induction step, suppose $0<i<j<k+1$.  Apply an identity from Lemma~\ref{lem:q-binomialCounts}, so that
\begin{align*}
\sum_{\ell = 0}^j (-1)^{\ell} q^{\binom{\ell}{2}} \qbinom{j}{\ell} \qbinom{k+1-\ell}{i}
&= \sum_{\ell = 0}^j (-1)^{\ell} q^{\binom{\ell}{2}} \qbinom{j}{\ell} \qbinom{k-\ell}{i-1}  \\
&\quad + q^i \sum_{\ell = 0}^j (-1)^{\ell} q^{\binom{\ell}{2}} \qbinom{j}{\ell} \qbinom{k-\ell}{i} = 0,
\end{align*}
because both sums on the right side vanish by the induction hypothesis or the edge cases.
\end{proof}

\begin{proof}[Proof of {Lemma~\textup{\ref{lem:AnnT-lowerTriangular}}}]
We need to show that the $i,j$-entry of each product vanishes when $i<j$.  For $\sW_0 T$, the result is the $w_0=0$ analog of Proposition~\ref{prop:PW-format}, as $T$ is  
$Q_{0,2}$.  For $\sAnn T$, 
\[
(\sAnn T)_{i,j}
= \sum_{\ell = 1}^j \mathscr{S}_i \left(\qbinom{k-\ell}{i} - \qbinom{k}{i} \right) (-1)^{\ell} q^{\binom{\ell}{2}} \qbinom{j}{\ell} .
\]
Note that the sum does not change if we include $\ell=0$.  Applying the distributive law and  
\eqref{eq:Cauchy-vanishing}, as $j>0$, we see from Lemma~\ref{lem:NewSumVanishing} that
\[
(\sAnn T)_{i,j} = \mathscr{S}_i  \sum_{\ell = 0}^j  (-1)^{\ell} q^{\binom{\ell}{2}} \qbinom{j}{\ell} \qbinom{k-\ell}{i} =0.  \qedhere
\]
\end{proof}

We now define two filtrations of $\Q^k$.  Let $e_1 = (1, 0, \ldots, 0), \ldots, e_k=(0, \ldots, 0,1)$ be the standard basis vectors of $\Q^k$.  For $i = 1, 2, \ldots, k$, define linear subspaces of $\Q^k$:
\begin{align*}
\mathcal{V}_i &= \spn\{e_i, e_{i+1}, \ldots, e_k\}, \\
\mathcal{T}_i &= \spn\{ \text{columns $i, i+1, \ldots, k$ of the matrix $T$} \} .
\end{align*}
Then $\dim \mathcal{V}_i = \dim \mathcal{T}_i = k-i+1$, and
\begin{align*}
\Q^k &= \mathcal{V}_1 \supset \mathcal{V}_2 \supset \cdots \supset \mathcal{V}_k \supset \{0\}, \\
\Q^k &= \mathcal{T}_1 \supset \mathcal{T}_2 \supset \cdots \supset \mathcal{T}_k \supset \{0\}.
\end{align*}
Lemma~\ref{lem:AnnT-lowerTriangular} shows that the matrices $\sAnn$ and $\sW_0$ (when invertible) map the $\mathcal{T}$-filtration isomorphically to the $\mathcal{V}$-filtration.  

\begin{thm}  \label{thm:SingletonContributions}
Let $C$ and $D$ be the linear codes constructed using \eqref{eq:defnCD} and $\varsigma^{(s)}$ for some $s$, $1 \leq s < k$, with $\sW_0$ invertible.  Then the contribution by rank $i$ singletons to $\delta \asing_{w_i} = \asing_{w_i}(D^\perp) - \asing_{w_i}(C^\perp)$ is zero if $i<s$.  The contribution of rank $s$ singletons to $\delta \asing_{w_s}$ is nonzero.
\end{thm}

\begin{proof}
The vector $\bar\varsigma^{(s)}$ belongs to $\mathcal{V}_s - \mathcal{V}_{s+1}$.  Then $\os^{(s)} = \sW_0^{-1} \bar\varsigma^{(s)} \in \mathcal{T}_s - \mathcal{T}_{s+1}$.  This, in turn, implies that $\sAnn \os^{(s)} \in \mathcal{V}_s - \mathcal{V}_{s+1}$.  By \eqref{eq:SingletonContribs}, singletons of rank $i$, $i < s$, make zero contribution to $\delta \asing_{w_i}$, while the singletons of rank $s$ make a nonzero contribution to $\delta \asing_{w_s}$.
\end{proof}

We point out that 
Theorem~\ref{thm:SingletonContributions} makes no claims about the contribution of singletons of rank $k$ to $\delta \asing_{w_k}$.

\begin{rem}  \label{rem:DegenerateImageW0}
If $W_0$ is degenerate because $c_j = 0$, then $\sW_0$ is also degenerate, and $\sW_0$ maps $\mathcal{T}_j$ into $\mathcal{V}_{j+1}$, and $\mathcal{T}_i$, $i < j$, will map to a proper linear subspace of $\mathcal{V}_i$.  If $j$ is the largest index such that $c_j=0$, then $\os^{(s)}$ will be in the image of $\sW_0$ provided $j<s$.  Recalling that $s < k$ and that $c_1 \neq 0$, we see that $\os^{(s)} \in \im\sW_0$ when $2 \leq j <s<k$.  For example, when $k=3$, no such $s$ can exist. 
\end{rem}

\section{Main results}  \label{sec:MainResults}

We are now in a position to prove that a large number of weights with maximal symmetry, including the homogeneous weight, do not respect duality.  

As usual in this part of the paper, let $R=M_{k \times k}(\F_q)$ with $k \geq 2$.  Suppose $w$ is a weight on $R$ that has maximal symmetry, positive integer values, and $w(0)=0$.  The value $w(r)$, $r \in R$, depends only on the rank $\rk r$ of $r$.  Write $w_i$ for the common value $w(r)$ when $\rk r = i$. 

We will say that a weight $w$ is \emph{nondegenerate} if the expressions $c_i$ of \eqref{eqn:diag-coefficents} are nonzero for all $i=1, 2, \ldots, k$.  (Note that $c_0=0$ because $w(0)=0$, and $c_1 = w_1 >0$.)  If at least one of $c_2, \ldots, c_k$ vanishes, we say that $w$ is \emph{degenerate}.

Let $\mathring{w} = \min\{w_1, w_2, \ldots, w_k \}$; $\mathring{w}$ is a positive integer.  Write $\mathring{I} = \{ i: w_i = \mathring{w} \}$ for set of indices $i$ where $w_i$ achieves the minimum positive value $\mathring{w}$.  In general, for an integer $d \geq \mathring{w}$, set $I_d = \{ i: w_i=d \}$; depending on $d$, $I_d$ may be empty.  Of course, $\mathring{I} = I_{\mathring{w}}$ is nonempty.

\begin{thm}  \label{thm:generalMainResult}
Assume $w$ is a nondegenerate weight on $R = M_{k \times k}(\F_q)$, with maximal symmetry, positive integer values, and $w(0)=0$.  Suppose there is an integer $d$ such that $\mathring{w} \leq d < 2 \mathring{w}$, $I_d$ is nonempty, and $k \not\in I_d$.  Then $w$ does not respect duality: there exist linear codes $C$ and $D$ over $R$ of the same length such that $\wwe_C = \wwe_D$ and $A_d(C^\perp) \neq A_d(D^\perp)$.
In particular, if $k \not\in \mathring{I}$, then $w$ does not respect duality.
\end{thm}

\begin{proof}
Suppose $d$ satisfies the stated hypotheses, and let $s = \max I_d$; then $s < k$.  Construct linear codes $C$ and $D$ over $R$ using \eqref{eq:defnCD} using $s = \max I_d$.  Because $s<k$, Theorem~\ref{thm:MainConstruction} implies that the codes have the same length and that $\wwe_C = \wwe_D$.  

As for the dual codes, Corollary~\ref{cor:only-singletons} says that $A_d(D^\perp) - A_d(C^\perp) = \asing_d(D^\perp) - \asing_d(C^\perp) = \delta \asing_d$.  The only singletons that can contribute to $\delta\asing_d$ are those of rank $i$ with $i \in I_d$, Corollary~\ref{cor:contribsBySingletons}.  Theorem~\ref{thm:SingletonContributions} says that singletons of rank $i<s$ make zero contributions to $\delta\asing_d$, while singletons of rank $s$ make a nonzero contribution to $\delta\asing_d$.  We conclude that $A_d(D^\perp) - A_d(C^\perp) = \delta\asing_d \neq 0$.
\end{proof}

\begin{rem}
If $k \in I_d$, the arguments given above are not conclusive.  In Theorem~\ref{thm:SingletonContributions} there is always the possibility that the contributions of singletons of rank $k$ could cancel the contributions of singletons of lower rank $i \in I_d$.  Even if $I_d = \{k\}$, it is possible that singletons of rank $k$ make zero contributions.  This happens for the homogeneous weight, for example; see Corollary~\ref{cor:sameEssentialLength} and Proposition~\ref{prop:SingletonContribs}.
\end{rem}

\begin{cor}   \label{cor:twoSmallValues}
Suppose a nondegenerate $w$ satisfies $\mathring{w} < w_i < 2 \mathring{w}$ for some $i$, $1 \leq i \leq k$.  Then $w$ does not respect duality.
\end{cor}

\begin{proof}
By the definition of $\mathring{w}$, there is some index $j$ so that $w_j = \mathring{w}$.  As $w_j < w_i$, we have $i \neq j$.  The index $k$ can equal at most one of $i$ or $j$.  Apply 
Theorem~\ref{thm:generalMainResult} to the other one.
\end{proof}

We can, at long last, prove that the homogeneous weight on $M_{k \times k}(\F_q)$ does not respect duality, provided $k>2$ or $q>2$. 
The homogeneous weight on any finite Frobenius ring has the Extension Property \cite[Theorem~2.5]{greferath-schmidt:combinatorics}.  Since the matrix ring $M_{k \times k}(\F_q)$ is Frobenius, it follows that $W_0$ and $\sW_0$ are invertible for the homogeneous weight and any information module $M$, Remark~\ref{rem:EPand W0}.  This says that the homogeneous weight is nondegenerate.

\begin{thm}  \label{thm:homog-failures}
Let $R = M_{k \times k}(\F_q)$, $k \geq 2$, with the homogeneous weight $\wg$.  Then $\wg$ respects duality if and only if $k=2$ and $q=2$.
\end{thm}

\begin{proof}
The `if' portion is Theorem~\ref{thm:MWIds-homog-2by2-F2}.  For the `only if' portion,  Lemma~\ref{lem:order-of-homog} says $\wg_2 = \mathring{\wg}$ is the smallest nonzero value of $\wg$, while $\wg_1$ is the largest value, with $\wg_2 < \wg_1 \leq 2\wg_2$.  There is equality $\wg_1 = 2\wg_2$ if and only if $k=2$ and $q=2$.  For all other values of $k$ and $q$, there is strict inequality: $\wg_2 < \wg_1 < 2 \wg_2$.  Now apply Theorem~\ref{thm:generalMainResult} with $d=\wg_1$.
\end{proof}

The last result of this section determines, for $k=2$, all the weights with maximal symmetry (nondegenerate or not) that respect duality.  The list is a short one: the Hamming weight (any $q$) and the homogeneous weight (only when $q=2$).

\begin{thm}  \label{main-theorem-for-2by2}
Let $w$ be a weight on $R=M_{2 \times 2}(\F_q)$ having maximal symmetry, positive integer values, and $w(0)=0$.  Assume $w$ is neither a multiple of the Hamming weight \textup{(}$w_1=w_2$\textup{)} nor, when $q=2$, a multiple of the homogeneous weight \textup{(}$w_1=2w_2$ when $q=2$\textup{)}.  Then $w$ does not respect duality.
\end{thm}

\begin{proof}
Using $m=k+1=3$ and $s=1$ (the only possible value of $s$), Example~\ref{ex:sWfor2} gives the $\sW_0$-matrix:
\[  \sW_0 = \begin{bmatrix}
q^2 w_1 & (q^2+q) w_1 \\
(q^2+q) w_1 & (q+1) w_1 + q^2 w_2
\end{bmatrix} .  \]
Then $\det \sW_0 = q^3 w_1 (-(q+1) w_1 + q w_2)$.  As $w_1 >0$ by hypothesis, $\det \sW_0$ vanishes only when $-(q+1) w_1 + q w_2=0$.

First consider the nondegenerate case, where $-(q+1) w_1 + q w_2 \neq 0$.  Use the construction of \eqref{eq:defnCD} with $m=3$ and $s=1$ to produce linear codes $C$ and $D$ with $\wwe_C = \wwe_D$.  The net counts of singleton vectors in the dual codes depends only on $\os$.  A calculation shows that 
\begin{align*}
\sW_0 \langle (q+1) w_1 &+ (q^2-q) w_2, -q^2 w_1 \rangle \\
&= (q(-(q+1) w_1 + q w_2)) \langle -(q^2-q), 1 \rangle \\
&= (q(-(q+1) w_1 + q w_2)) \bar\varsigma^{(1)} .
\end{align*} 
Thus, up to scaling, we take $\os^{(1)} = \langle (q+1) w_1 + (q^2-q) w_2, -q^2 w_1 \rangle$.

Using Example~\ref{ex:AnnMatrixk2}, we see that the matrix $\sAnn$ is
\[  \sAnn = \begin{bmatrix}
-q(q^2-1) & -(q^2-1)(q+1) \\
-(q^2-1)(q^2-q) & -(q^2-1)(q^2-q)
\end{bmatrix} .  \]
By \eqref{eq:SingletonContribs}, the net contributions to $\delta \asing_{w_j} = \asing_{w_j}(D^\perp) - \asing_{w_j}(C^\perp)$ by rank $j$ singletons are given by the entries of $\sAnn \os$:
\begin{equation}  \label{eq:k2singletons}
 \sAnn \os = q (q-1) (q^2-1) \begin{bmatrix}
 (q+1) w_1 - q w_2 \\
(q^2-q-1) w_1 - (q^2-q) w_2
\end{bmatrix} .
\end{equation}
The contribution for $j=1$ is nonzero because $w$ is nondegenerate.  
The contribution for $j=2$ is nonzero, provided $w$ is not a multiple of the homogeneous weight; see Example~\ref{ex:2by2case}.

Because $k=2$, there are just a few (nondegenerate) cases:
\begin{itemize}
\item If $w_1 < w_2$, then $\delta A_{w_1} = \delta \asing_{w_1} \neq 0$, by \eqref{eq:k2singletons}.
\item  If $w_2 < w_1$ and $w$ is not homogeneous, then $\delta A_{w_2} = \delta \asing_{w_2} \neq 0$, by \eqref{eq:k2singletons}.
\item  If $w$ is homogeneous, apply Theorem~\ref{thm:homog-failures} (except when $q=2$).
\item  If $w_1=w_2$, $w$ is a multiple of the Hamming weight.
\end{itemize}

Now suppose $w$ is degenerate, so that
$q w_2 = (q+1) w_1$.  Because $q$ and $q+1$ are relatively prime, there exists a positive integer $\tau$ such that $w_1 = q \tau$ and $w_2 = (q+1) \tau$.  Then $w_1 < w_2 < 2 w_1$, as $q \geq 2$.  The degenerate matrix $\sW_0$ becomes
\[  \sW_{\rm degen} = \begin{bmatrix}
q^2 w_1 & q(q+1) w_1 \\
q(q+1) w_1 & (q+1)^2 w_1
\end{bmatrix} .  \]
A basis for $\ker \sW_{\rm degen}$ is $\langle -(q+1), q \rangle$.

Use the linear codes $C_{\pm}$ of Proposition~\ref{prop:DegenerateCase}.  They have the same length and $w$-weight enumerators.  Their net rank-sums $\eb_i(C_+) - \eb_i(C_-)$ are $1, -(q+1), q$, respectively.  Then the contributions of singletons are:
\[  \sAnn \begin{bmatrix}
-(q+1) \\
q
\end{bmatrix}
=  \begin{bmatrix}
0 \\
(q^2-1)(q^2-q)
\end{bmatrix} . \]
Because $w_1 < w_2 < 2 w_1$, Corollary~\ref{cor:only-singletons} applies to $w_2$.  Then $A_{w_2}(C_+^\perp) - A_{w_2}(C_-^\perp) = \asing_{w_2}(C_+^\perp) - \asing_{w_2}(C_-^\perp) = (q^2-1)(q^2-q) >0$.
\end{proof}

In Section~\ref{sec:3by3case}, the case of $R = M_{3 \times 3}(\F_2)$ is discussed in detail.

\section{Rank partition enumerators}  \label{sec:RankPartitionEnum}

Section~\ref{sec:MWIdentities} described various enumerators including the complete enumerator and symmetrized enumerators.  In this section we focus on a particular enumerator, the rank partition enumerator, over the matrix ring $M_{k \times k}(\F_q)$.  The rank partition enumerator is a partition enumerator associated to rank, and it is coarser than the symmetrized enumerator associated to the group action of $\GL(k, \F_q)$ acting on $M_{k \times k}(\F_q)$.

On $R = M_{k \times k}(\F_q)$, define the \emph{rank partition} $\mathcal{RK} = \{ R_i \}_{i=0}^k$, with 
\[  R_i = \{ s \in R: \rk s = i \} . \]
As in Section~\ref{sec:MWIdentities}, define counting functions $n_i: R^n \ra \N$, $i=0, 1, \ldots, k$, by $n_i(x) = \size{\{ j : x_j \in R_i\} }$, for $x = (x_1, x_2, \ldots, x_n) \in R^n$.  For a linear code $C \subseteq R^n$, define the \emph{rank partition enumerator} $\re_C$ associated to $C$ to be the homogeneous polynomial of degree $n$ in the variables $Z_0, Z_1, \ldots, Z_k$ given by
\[  \re_C(Z) = 
\sum_{x \in C} \prod_{j=1}^n Z_{\rk x_j} = \sum_{x \in C} \prod_{i=0}^k Z_i^{n_i(x)} . \]

If $w$ is a weight on $R$ with maximal symmetry and positive integer values, then the value of $w(r)$, $r \in R$, depends only on the rank $\rk r$ of $r$.  Write $w_i$ for $w(r)$ when $\rk r = i$, and denote by $w_{\max}$ the largest value of $w$.  Then the specialization of variables $Z_i \leadsto X^{w_{\max}-w_i} Y^{w_i}$ allows one to write $\wwe_C$ in terms of $\re_C$:
\begin{equation}  \label{eq:wweViare}
\wwe_C(X,Y) = \re_C(Z)|_{Z_i \leadsto X^{w_{\max}-w_i} Y^{w_i}} .
\end{equation}

As an example of some of the results of \cite{MR3336966}, we will show that the rank partition enumerator satisfies the MacWilliams identities.  Then \eqref{eq:wweViare} will allow us to calculate $\wwe_{C^\perp}$ for many examples.  This will be one way to illustrate the main results of Section~\ref{sec:MainResults} (and to prove additional results).

In order to show that the MacWilliams identities hold for the rank partition enumerator, 
we refer to the argument outlined in Appendix~\ref{sec:appendix} and provide details on the relevant Fourier transforms.

It is well-known (\cite[Example~4.4]{wood:duality}) that $R=M_{k \times k}(\F_q)$ is a Frobenius ring with a generating character $\chi$.  To describe the standard generating character $\chi$, we first recall the standard generating character $\theta_q$ of $\F_q$:  $\theta_q(a) = \zeta^{\Tr_{q \ra p}(a)}$, $a \in \F_q$, where $q = p^e$, $p$ prime, $\zeta=\exp(2 \pi i/p) \in \C^\times$, and $\Tr_{q \ra p}$ is the absolute trace from $\F_q \ra \F_p$.  Then define 
$\chi(s) = \theta_q(\tr s)$, $s \in R$, where $\tr$ is the matrix trace.  
Because the matrix trace over $\F_q$ satisfies $\tr(rs)=\tr(sr)$, we see that $\chi(rs) = \chi(sr)$ for all $r,s \in R$.  
As is the case for all generating characters, $\chi$ has the property that $\ker\chi$ contains no nonzero left or right ideal of $R$, \cite[Lemma~4.1]{wood:duality}.

\begin{lem}  \label{lem:block-char-sums}
Each partition block $R_i$ of $\mathcal{RK}$ is invariant under left or right multiplication by units. 
If $\rk r_1= \rk r_2$, then $\sum_{s \in R_i} \chi(s r_1) = \sum_{s \in R_i} \chi(s r_2)$ for all $i$.
\end{lem}

\begin{proof}
The rank of a matrix is invariant under multiplication by units.  If $\rk r_1= \rk r_2$, then, using row and column operations, there are units $u_1, u_2$ such that $r_2 = u_1 r_1 u_2$.  Thus,
\begin{align*}
\sum_{s \in R_i} \chi(s r_2) = \sum_{s \in R_i} \chi(s u_1 r_1 u_2) = \sum_{s \in R_i} \chi(u_2 s u_1 r_1) = \sum_{s \in R_i} \chi(s r_1) ,
\end{align*}  
using the property $\chi(rs) = \chi(sr)$ and the bi-invariance of $R_i$.
\end{proof}

Define the \emph{Kravchuk matrix} $K$ of size $(k+1) \times (k+1)$ for the rank partition $\mathcal{RK}$ by
\[  K_{i,j} = \sum_{s \in R_i} \chi(s r) , \quad i,j=0, 1, \ldots, k, \]
where $r \in R$ has $\rk r = j$.  This sum is well-defined by Lemma~\ref{lem:block-char-sums}.  

In order to develop an explicit formula for $K_{i,j}$, we 
first remark that $R_i$, being invariant under left multiplication by units, equals the disjoint union of the left $\U$-orbits it contains.  We already know the number and sizes of the $\U$-orbits in $R$,  Proposition~\ref{prop:orbitNumberSize}.  So, we turn our attention to sums of the form $\sum_{t \in \orb(s)} \chi(tr)$, for $r \in R$.

Let $\PP_R$ be the poset of all principal left ideals of $R$ under set containment; all the left ideals of $R$ are principal.  (This is the type of poset used in Proposition~\ref{prop:poset-isom}, with $M=R$ and $m=k$.)  Fix an element $r \in R$, and define two functions $f_r, g_r: \PP_R \ra \C$ by 
\begin{align*}
f_r(Rs) = \sum_{t \in Rs} \chi(tr) ,  \quad
g_r(Rs) = \sum_{t \in \orb(s)} \chi(tr) .
\end{align*}
The definition of $g_r$ is well-defined by Lemma~\ref{lem:Rx-rowspace}.
We collect some facts about $f_r$ and $g_r$ in the next lemma.  For $r \in R$, its left \emph{annihilator} is $\ann_{\lt}(r) = \{ s \in R: s r =0 \}$; the left annihilator is a left ideal of $R$.

\begin{lem}  \label{lem:FactAboutfr}
For any $r,s \in R$, we have
\[  f_r(Rs) = \sum_{Rt \subseteq Rs} g_r(Rt) . \]
The values of $f_r$ are
\[  f_r(Rs) = \begin{cases}
\size{Rs}, & \text{if $Rs \subseteq \ann_{\lt}(r)$,} \\
0, & \text{otherwise}.
\end{cases} \]
\end{lem}

\begin{proof}
The first equality reflects the fact that the left ideal $Rs$ is invariant under left multiplication by units, and hence $Rs$ equals the disjoint union of the left $\U$-orbits it contains. 

As in \eqref{eq:FrakDualCode}, denote $\mathfrak{R}(Rs) = \{r \in R: \chi(t r) = 1, \text{for all $t \in Rs$}\}$.  We claim that $r \in \mathfrak{R}(Rs)$ if and only if $Rs \subseteq \ann_{\lt}(r)$.  The `if' direction is clear: if $Rs \subseteq \ann_{\lt}(r)$, then $tr=0$, and hence $\chi(tr)=1$, for all $t \in Rs$.  Conversely, suppose $r \in  \mathfrak{R}(Rs)$.  Then $Rsr \subseteq \ker\chi$.  We conclude that $Rsr=0$, because any left ideal in $\ker\chi$ must be zero.  Thus $Rs \subseteq \ann_{\lt}(r)$.  The second formula now follows from Lemma~\ref{lem:charSumAnnih}.
\end{proof}

\begin{prop}
Let $\mu$ be the M\"obius function of the poset $\PP_R$.  Then
\[ g_r(Rs) = \sum_{Rt \subseteq Rs} \mu(Rt, Rs) f_r(Rt), \quad Rs \in \PP_R . \]
Simplifying, we have
\[  g_r(Rs) = \sum_{j=0}^{\rk s - \rk(sr)} (-1)^{\rk s - j} q^{\binom{\rk s - j}{2}} q^{k j} \qbinom{\rk s - \rk(sr)}{j} . \]
\end{prop}

\begin{proof}
The first equation comes from M\"obius inversion, \cite[Theorem~5.5.5]{MR4249619}, because of the first equation in Lemma~\ref{lem:FactAboutfr}.

By Proposition~\ref{prop:poset-isom}, the poset $\PP_R$ is isomorphic to the poset $\PP_{k,k}$ of linear subspaces of $\F_q^k$.  This allows us to translate the equation for $g_r(Rs)$ into geometric language.  Let $r \in R$ correspond to a linear subspace $Y \subseteq \F_q^k$, with $\dim Y = \rk r$.  Then $\size{Rr} = q^{k \dim Y} = q^{k \rk r}$, by Proposition~\ref{prop:orbitNumberSize}.  Similarly, let $s, t \in R$ correspond to $X, T$, with $\dim X = \rk s$ and $\dim T = \rk t$.  The M\"obius function of $\PP_{k,k}$ is in \eqref{eq:matrixMoebius}.  The condition $Rs \subseteq \ann_{\lt}(r)$ becomes $X \subseteq Y^\perp$.

Using Lemma~\ref{lem:FactAboutfr}, the formula for $g_r(Rs)$ simplifies:
\begin{align*}
g_r(Rs) 
&= \sum_{Rt \subseteq Rs \cap \ann_{\lt}(r)} \mu(Rt,Rs) \size{Rt} \\
&= \sum_{T \subseteq X \cap Y^\perp} (-1)^{\dim X - \dim T} q^{\binom{\dim X - \dim T}{2}} q^{k \dim T} \\
&= \sum_{j=0}^{\dim(X \cap Y^\perp)} (-1)^{\rk s - j} q^{\binom{\rk s - j}{2}} q^{k j} \qbinom{\dim(X \cap Y^\perp)}{j} . 
\end{align*}
Finally, Lemma~\ref{lem:dim-formulas} implies $\dim(X \cap Y^\perp) = \rk s - \rk(sr)$.
\end{proof}

To simplify notation slightly, define 
\begin{equation}
B(i,\ell) = \sum_{j=0}^{\ell} (-1)^{i-j} q^{\binom{i-j}{2}} q^{kj} \qbinom{\ell}{j} ,
\end{equation}
for $i = 0, 1, \ldots, k$ and $0 \leq \ell \leq i$.  Then $\sum_{t \in \orb(s)} \chi(tr) = g_r(Rs) = B(\rk s, \rk s - \rk(sr))$.
In addition, suppose there are linear subspaces $C \subseteq A \subseteq \F_q^k$ with $\dim A = a$ and $\dim C = c$.  Then
define $I(a,b,c,d)= \size{\{ B \subseteq A: \dim B = b, \text{ and } \dim(B \cap C)= d\}}$.  By Lemma~\ref{lem:counting-intersections} (with $D=A \cap C$), 
\[  I(a,b,c,d) = q^{(b-d)(c-d)} \qbinom{c}{d} \qbinom{a-c}{b-d} . \]

\begin{prop}
The Kravchuk matrix $K$ has entries
\[  K_{i,j} = \sum_{\ell = 0}^i I(k,i,k-j,\ell) B(i,\ell) , \]
for $i,j = 0, 1, \ldots, k$.
\end{prop}

\begin{proof}
As mentioned earlier, if $j = \rk r$, the sum in $K_{i,j} = \sum_{s \in R_i} \chi(sr)$ can be split up into sums over the left $\U$-orbits contained in $R_i$.  The individual sums over orbits depend upon $\rk s$ and $\rk(sr)$, so we need to count the number of orbits $\orb(s)$ with $\rk s = i$ for various values of $\rk(sr)$.

In terms of linear subspaces, we need to count the number of linear subspaces $X$ of $\F_q^k$ with $\dim X = i$ and $\dim(X \cap Y^\perp) = \ell$.  Because $\dim Y^\perp = k - \rk r = k-j$, this count is $I(k,i,k-j,\ell)$.  
\end{proof}

\begin{ex}  \label{ex:Kravk2}
For $k=2$, the Kravchuk matrix is:
\[  K = \left[\begin{array}{rrr}
1 & 1 & 1 \\
{\left(q^{2} - 1\right)} {\left(q + 1\right)} & q^{2} - q - 1 & -q - 1 \\
(q^2-q)(q^2-1) & -q^{2} + q & q
\end{array}\right] . \]
\end{ex}

Suppose $C \subseteq R^n$ is an additive code.  The annihilators $\mathfrak{L}(C)$ and $\mathfrak{R}(C)$ were defined in \eqref{eq:FrakDualCode}.  The MacWilliams identities for the rank partition enumerator are next.
\begin{thm}  \label{thm:MWids-rankPartition}
Let $R = M_{k \times k}(\F_q)$ with Kravchuk matrix $K$.  If $C$ is an additive code in $R^n$,  then
\begin{align*}
\re_{\mathfrak{L}(C)}(Z_i) = \frac{1}{\size{C}} \re_C(\mathcal{Z}_j)|_{\mathcal{Z}_j = \sum_i K_{i,j} Z_i}  , \\
\re_{\mathfrak{R}(C)}(Z_i) = \frac{1}{\size{C}} \re_C(\mathcal{Z}_j)|_{\mathcal{Z}_j = \sum_i K_{i,j} Z_i} .
\end{align*}
The formulas are reversible in $C$ and $\mathfrak{L}(C)$, resp., $C$ and $\mathfrak{R}(C)$.

If $C \subseteq R^n$ is a left, resp., right, $R$-linear code, then $\mathcal{R}(C) =\mathfrak{R}(C)$, resp., $\mathcal{L}(C) = \mathfrak{L}(C)$.
\end{thm}

\begin{proof}
We add details to the outline provided in Appendix~\ref{sec:appendix}.  Let $V = \C[Z_0, Z_1, \ldots, Z_k]$, and define $f: R \ra V$ by $f(s) = Z_{\rk s}$, $s \in R$.  We calculate the Fourier transform of $f$, as in \eqref{eqn:FTdefn}.  Write $j = \rk r$.
\begin{align*}
\ft{f}(r) &= \sum_{s \in R} \chi(rs) Z_{\rk s} = \sum_{s \in R} \chi(sr) Z_{\rk s} \\
&= \sum_{i=0}^k \sum_{s \in R_i} \chi(sr) Z_i = \sum_{i=0}^k K_{i,j} Z_i . 
\end{align*}
Note that $\ft{f}(r)$ depends only on $\rk r$, so that $\ft{f}(r)$ equals $f(r)$ after applying the linear substitution $Z_j \leftarrow  \sum_{i=0}^k K_{i,j} Z_i$.

To reverse roles, use Lemma~\ref{lem:DualityFacts} and apply the formulas to the pair $\mathfrak{L}(C)$ and $\mathfrak{R}(\mathfrak{L}(C))=C$ and the pair $\mathfrak{R}(C)$ and $\mathfrak{L}(\mathfrak{R}(C))=C$.  For the case of linear codes, see Remark~\ref{rem:LinearCaseAnnihilators}.
\end{proof}

\section{Examples}  \label{sec:MatrixExamples}
In this section we calculate a number of examples over $R = M_{2 \times 2}(\F_2)$.
Set $m=3$, so that the information module is $M=M_{2 \times 3}(\F_2)$.  The orbit spaces $\OO_1$ and $\OO_2$ have representatives in reduced row-echelon form and are ordered as follows:  
\begin{align}
\OO_1 &= \left[\left[\begin{smallmatrix}
1 & 0 & 1 \\
0 & 0 & 0
\end{smallmatrix}\right], \left[\begin{smallmatrix}
0 & 1 & 0 \\
0 & 0 & 0
\end{smallmatrix}\right], \left[\begin{smallmatrix}
0 & 0 & 1 \\
0 & 0 & 0
\end{smallmatrix}\right], \left[\begin{smallmatrix}
1 & 1 & 1 \\
0 & 0 & 0
\end{smallmatrix}\right], \left[\begin{smallmatrix}
0 & 1 & 1 \\
0 & 0 & 0
\end{smallmatrix}\right], \left[\begin{smallmatrix}
1 & 0 & 0 \\
0 & 0 & 0
\end{smallmatrix}\right], \left[\begin{smallmatrix}
1 & 1 & 0 \\
0 & 0 & 0
\end{smallmatrix}\right]\right] , \notag \\
\OO_2 &= \left[\left[\begin{smallmatrix}
1 & 0 & 0 \\
0 & 1 & 0
\end{smallmatrix}\right], \left[\begin{smallmatrix}
1 & 0 & 1 \\
0 & 1 & 0
\end{smallmatrix}\right], \left[\begin{smallmatrix}
1 & 0 & 0 \\
0 & 0 & 1
\end{smallmatrix}\right], \left[\begin{smallmatrix}
1 & 1 & 0 \\
0 & 0 & 1
\end{smallmatrix}\right], \left[\begin{smallmatrix}
1 & 0 & 0 \\
0 & 1 & 1
\end{smallmatrix}\right], \left[\begin{smallmatrix}
0 & 1 & 0 \\
0 & 0 & 1
\end{smallmatrix}\right], \left[\begin{smallmatrix}
1 & 0 & 1 \\
0 & 1 & 1
\end{smallmatrix}\right]\right] . \label{ex:OrderingOO}
\end{align}
The rank $1$ orbits have size $3$, and the rank $2$ orbits have size $6$, as in Proposition~\ref{prop:orbitNumberSize}.
The representatives of $\OS_1$ and $\OS_2$ are the transposes of the representatives of $\OO_1$ and $\OO_2$, using the same orderings,
with rank $1$ coming before rank $2$.
The $W_0$-matrix has size $14 \times 14$, while $\sW_0$ is $2 \times 2$:
\begin{align*}
W_0 &= \left[\begin{smallmatrix}
0 & 0 & w_{1} & 0 & w_{1} & w_{1} & w_{1} & w_{1} & 0 & w_{1} & w_{1} & w_{1} & w_{1} & w_{1} \\
0 & w_{1} & 0 & w_{1} & w_{1} & 0 & w_{1} & w_{1} & w_{1} & 0 & w_{1} & w_{1} & w_{1} & w_{1} \\
w_{1} & 0 & w_{1} & w_{1} & w_{1} & 0 & 0 & 0 & w_{1} & w_{1} & w_{1} & w_{1} & w_{1} & w_{1} \\
0 & w_{1} & w_{1} & w_{1} & 0 & w_{1} & 0 & w_{1} & w_{1} & w_{1} & w_{1} & w_{1} & w_{1} & 0 \\
w_{1} & w_{1} & w_{1} & 0 & 0 & 0 & w_{1} & w_{1} & w_{1} & w_{1} & w_{1} & 0 & w_{1} & w_{1} \\
w_{1} & 0 & 0 & w_{1} & 0 & w_{1} & w_{1} & w_{1} & w_{1} & w_{1} & w_{1} & w_{1} & 0 & w_{1} \\
w_{1} & w_{1} & 0 & 0 & w_{1} & w_{1} & 0 & w_{1} & w_{1} & w_{1} & 0 & w_{1} & w_{1} & w_{1} \\
w_{1} & w_{1} & 0 & w_{1} & w_{1} & w_{1} & w_{1} & w_{2} & w_{2} & w_{1} & w_{1} & w_{2} & w_{1} & w_{2} \\
0 & w_{1} & w_{1} & w_{1} & w_{1} & w_{1} & w_{1} & w_{2} & w_{1} & w_{1} & w_{2} & w_{2} & w_{2} & w_{1} \\
w_{1} & 0 & w_{1} & w_{1} & w_{1} & w_{1} & w_{1} & w_{1} & w_{1} & w_{2} & w_{2} & w_{2} & w_{1} & w_{2} \\
w_{1} & w_{1} & w_{1} & w_{1} & w_{1} & w_{1} & 0 & w_{1} & w_{2} & w_{2} & w_{1} & w_{2} & w_{2} & w_{1} \\
w_{1} & w_{1} & w_{1} & w_{1} & 0 & w_{1} & w_{1} & w_{2} & w_{2} & w_{2} & w_{2} & w_{1} & w_{1} & w_{1} \\
w_{1} & w_{1} & w_{1} & w_{1} & w_{1} & 0 & w_{1} & w_{1} & w_{2} & w_{1} & w_{2} & w_{1} & w_{2} & w_{2} \\
w_{1} & w_{1} & w_{1} & 0 & w_{1} & w_{1} & w_{1} & w_{2} & w_{1} & w_{2} & w_{1} & w_{1} & w_{2} & w_{2}
\end{smallmatrix}\right]  , \\
\sW_0 &= \begin{bmatrix}
4 w_1 & 6 w_1 \\
6 w_1 & 3 w_1 + 4 w_2
\end{bmatrix} . 
\end{align*}

Akin to Figure~\ref{fig:comm-diagram} on page~\pageref{fig:comm-diagram}, the rank partition enumerator of a linear code specializes to the $w$-weight enumerator of the code under the specialization $Z_0 \leadsto 1$, $Z_1 \leadsto t^{w_1}$ and $Z_2 \leadsto t^{w_2}$.  Using the Kravchuk matrix $K$ from Example~\ref{ex:Kravk2}, with $q=2$: 
\[  K= \left[ \begin{array}{rrr}
1 & 1 & 1 \\
9 & 1 & -3 \\
6 & -2 & 2
\end{array} \right] , \]
the MacWilliams identities for the rank partition enumerator yield the rank partition enumerator for the dual code, Theorem~\ref{thm:MWids-rankPartition}.  Using the same specialization, the $w$-weight enumerator of the dual code is obtained.  

In all of the examples that follow, multiplicity functions $\eta$ and lists $\omega$ of orbit weights are written in terms of the ordering of $\OO$ given in \eqref{ex:OrderingOO}, with ranks separated by vertical lines.  All calculations were performed in SageMath \cite{sagemath}.  Rank partition enumerators of dual codes are not listed because they would use too much space.  Only the lowest order terms in the $w$-weight enumerators of dual codes are displayed.

It is differences such as $A_j(D^\perp) - A_j(C^\perp)$ that ultimately matter, so we will write $\delta\eb = \eb(D) - \eb(C)$.  By \eqref{eq:netContribOfRankiSingletons} and Example~\ref{ex:AnnMatrixk2} with $q=2$, the net contributions of rank $i$ singletons to $\delta\asing_{w_i}=\asing_{w_i}(D^\perp) - \asing_{w_i}(C^\perp)$ appear as the entries in the vector
\begin{equation}  \label{eqn:SingletonContribForExamples}
\begin{bmatrix}
1 & 1 & 1 \\
9 & 3 & 0 \\
6 & 0 & 0
\end{bmatrix} \delta\eb .
\end{equation}

\begin{ex}
Set $w_1=2$ and $w_2 = 3$.  Then $\mathring{w}= 2$, and $\mathring{w} \leq w_i < 2 \mathring{w}$ for $i=1,2$.  This is a degenerate case, as $\det W_0 = 0$.  Use the codes in Proposition~\ref{prop:DegenerateCase}, with an extra rank $1$ functional so that the associated homomorphisms $\La$ are injective; see Remark~\ref{rem:functionalEvaluation}.  The multiplicity functions and lists of orbit weights are:
\begin{align*}
\eta_+ &= \langle 1|1,0,0,0,0,0,0|2,0,0,0,0,0,0 \rangle , \\
\omega_+ &= \langle 0|4, 4, 2, 4, 6, 6, 6| 8, 6, 6, 6, 8, 6, 8 \rangle , \\
\eta_- &= \langle 0|1,1,0,0,0,1,1|0,0,0,0,0,0,0 \rangle , \\
\omega_- &= \langle 0|4, 4, 2, 4, 6, 6, 6| 8, 6, 6, 6, 8, 6, 8 \rangle .
\end{align*}
The equation $\omega_+ = \omega_-$ is a feature of the construction in Proposition~\ref{prop:DegenerateCase}.
Both codes have length $4$, and $\delta\eb = \langle -1, 3, -2 \rangle$.  Then \eqref{eqn:SingletonContribForExamples} implies that $\delta A_2 = \delta\asing_2 = 0$, while $\delta A_3 = \delta\asing_3 = -6$.
The enumerators are
\begin{align*}
\se_{C_+} &= Z_{0}^{4} + 3 Z_{0}^{3} Z_{1} + 9 Z_{0}^{2} Z_{1}^{2} + 27 Z_{0} Z_{1}^{3} + 6 Z_{0}^{2} Z_{2}^{2} + 18 Z_{0} Z_{1} Z_{2}^{2} , \\
\se_{C_-} &= Z_{0}^{4} + 3 Z_{0}^{3} Z_{1} + 9 Z_{0}^{2} Z_{1}^{2} + 33 Z_{0} Z_{1}^{3} + 18 Z_{1}^{4} ; \\
\wwe_{C_+} &= 1 + 3 t^2 + 9 t^4 + 33 t^6 + 18 t^8 , \\
\wwe_{C_-} &= 1 + 3 t^2 + 9 t^4 + 33 t^6 + 18 t^8 ; \\
\wwe_{C_+^\perp} &= 1 + 12 t^2 + 6 t^3 + 36 t^4 + \cdots , \\
\wwe_{C_-^\perp} &= 1 + 12 t^2 \hphantom{{}+ 6 t^3} + 54 t^4 + \cdots.
\end{align*}
\end{ex}

\begin{ex}
Set $w_1=1$ and $w_2=2$, so the weight of a matrix equals its rank.  Then $\mathring{w}=1$, but $w_2 = 2 \mathring{w}$, so Corollary~\ref{cor:only-singletons} applies only to $w_1=1$.  Use the codes given by \eqref{eq:defnCD}, with $\varsigma$ and $\sigma$ given here:
\begin{align*}
\varsigma &= \langle 0,-1,-1,0,0,0,0|0,0,0,0,0,0,1 \rangle , \\
\sigma &= \langle -1, -1, -1, -1, -3, 1, -1| 0, 0, 0, 0, 0, 2, 2 \rangle .
\end{align*}
The scaling is such that $W_0 \sigma = 2 \varsigma$, so that $c=2$.  One calculates that $\al_1 = 10 w_1 = 10$ and $\al_2 = 9 w_1 + 4 w_2 = 17$.  It suffices to take $a=3$ and $b=0$.  Then $\Delta = -3$.  The multiplicity functions and lists of orbit weights are
\begin{align*}
\eta_{C_2} &= \langle 0| 3,3,3,26,3,3,3 | 3,3,3,3,3,3,3 \rangle , \\
\eta_{D_2} &= \langle 3| 2,2,2,25,0,4,2 | 3,3,3,3,3,5,5 \rangle ; \\
\omega_{C_2} &= \langle 0| 30, 53, 53, 53, 30, 53, 30 | 74, 74, 74, 74, 74, 74, 51 \rangle , \\
\omega_{D_2} &= \langle 0 | 30, 51, 51, 53, 30, 53, 30 | 74, 74, 74, 74, 74, 74, 53 \rangle .
\end{align*}

Both codes have length $65$, and $\delta \eb = \langle 3, -7, 4 \rangle$.  Then, using \eqref{eqn:SingletonContribForExamples}, $\delta A_1 = \delta \asing_1 = 6$ and $\delta \asing_2 = 18$.  The enumerators are
\begin{align*}
\se_{C_2} &= Z_{0}^{65} + 9 Z_{0}^{35} Z_{1}^{30} + 12 Z_{0}^{12} Z_{1}^{53} + 6 Z_{0}^{26} Z_{1}^{27} Z_{2}^{12} + 36 Z_{0}^{3} Z_{1}^{50} Z_{2}^{12} , \\
\se_{D_2} &= Z_{0}^{65} + 9 Z_{0}^{35} Z_{1}^{30} + 6 Z_{0}^{14} Z_{1}^{51} + 6 Z_{0}^{12} Z_{1}^{53} + 6 Z_{0}^{3} Z_{1}^{50} Z_{2}^{12} \\
&\quad\quad + 24 Z_{0}^{5} Z_{1}^{46} Z_{2}^{14} + 6 Z_{0}^{28} Z_{1}^{21} Z_{2}^{16} + 6 Z_{0}^{7} Z_{1}^{42} Z_{2}^{16} ; \\
\wwe_{C_2} &= 1 + 9\, t^{30} + 6\, t^{51} + 12\, t^{53} + 36\, t^{74} , \\
\wwe_{D_2} &= 1 + 9\, t^{30} + 6\, t^{51} + 12\, t^{53} + 36\, t^{74} ; \\
\wwe_{C_2^\perp} &=  1 + 132\, t + 15762\, t^2 + 1674894 \, t^3 + \cdots , \\
\wwe_{D_2^\perp} &= 1 + 138 \, t + 16176 \, t^{2} + 1695210 \, t^{3} + \cdots .
\end{align*}
Even though $\delta \asing_2 = 18$, $\delta A_2 = 414$; there are dual codewords of weight $2$ coming from vectors with two nonzero entries, both of rank $1$, that account for the difference.  This illustrates the importance of the strict inequality $\mathring{w} \leq d < 2 \mathring{w}$ in Corollary~\ref{cor:only-singletons}.
\end{ex}

\begin{ex}  \label{ex:linIndepDep}
Set $w_1=4$ and $w_2=5$.  Then $\mathring{w}=4$ and $\mathring{w} \leq w_i < 2 \mathring{w}$, $i=1,2$, so Corollary~\ref{cor:only-singletons} applies to both $w_1$ and $w_2$.  

In this example, 
the codes given by \eqref{eq:defnCD} will be used; call then $C_3$ and $D_3$.  In addition, two other codes $C_4, D_4$ will be given.  They are designed so that their lists of orbit weights have just three different values.  One of the values applies to all the rank $2$ orbits.  The other two values apply to three, resp., four, rank $1$ orbits.  For the code $C_4$, the three rank $1$ orbits are linearly independent, while for $D_4$ they are linearly dependent.

For $C_3$ and $D_3$, $\varsigma$ and $\sigma$ are given here:
\begin{align*}
\varsigma &= \langle 0,-1,-1,0,0,0,0|0,0,0,0,0,0,1 \rangle , \\
\sigma &= \langle 2, 2, 2, -1, 3, 1, 2 | 0, 0, 0, 0, 0, -4, -4 \rangle .
\end{align*}
The scaling is such that $W_0 \sigma = 8 \varsigma$, so that $c=2$.  One calculates that $\al_1 = 10 w_1 = 40$ and $\al_2 = 9 w_1 + 4 w_2 = 56$.  It suffices to take $a=1$ and $b=0$.  Then $\Delta = 3$.  The multiplicity functions and lists of orbit weights are
\begin{align*}
\eta_{C_3} &= \langle 3 | 4,4,4,22,4,4,4 | 4,4,4,4,4,4,4 \rangle , \\
\eta_{D_3} &= \langle 0 | 6,6,6,21,7,5,6 | 4,4,4,4,4,0,0 \rangle ; \\
\omega_{C_3} &= \langle 0| 160, 232, 232, 232, 160, 232, 160 \\
&\quad\quad\quad  | 296, 296, 296, 296, 296, 296, 224 \rangle , \\
\omega_{D_3} &= \langle 0 | 160, 224, 224, 232, 160, 232, 160 \\
&\quad\quad\quad | 296, 296, 296, 296, 296, 296, 232 \rangle .
\end{align*}

Both codes have length $77$, and $\delta \eb = \langle -3, 11, -8 \rangle$.  Then, using \eqref{eqn:SingletonContribForExamples}, $\delta A_4 = \delta \asing_4 = 6$ and $\delta A_5 = \delta \asing_5 = -18$.  The enumerators are
\begin{align*}
\se_{C_3} &= Z_{0}^{77} + 9 Z_{0}^{37} Z_{1}^{40} + 12 Z_{0}^{19} Z_{1}^{58} + 6 Z_{0}^{25} Z_{1}^{36} Z_{2}^{16} + 36 Z_{0}^{7} Z_{1}^{54} Z_{2}^{16} , \\
\se_{D_3} &= Z_{0}^{77} + 9 Z_{0}^{37} Z_{1}^{40} + 6 Z_{0}^{21} Z_{1}^{56} + 6 Z_{0}^{19} Z_{1}^{58} + 6 Z_{0}^{21} Z_{1}^{48} Z_{2}^{8} \\
&\quad\quad + 6 Z_{0}^{5} Z_{1}^{64} Z_{2}^{8} + 24 Z_{0}^{6} Z_{1}^{59} Z_{2}^{12} + 6 Z_{0}^{7} Z_{1}^{54} Z_{2}^{16} ; \\
\wwe_{C_3} &= 1 + 9\, t^{160} + 6\, t^{224} + 12\, t^{232} + 36\, t^{296} , \\
\wwe_{D_3} &= 1 + 9\, t^{160} + 6\, t^{224} + 12\, t^{232} + 36\, t^{296} ; \\
\wwe_{C_3^\perp} &=  1 + 165 \, t^{4} + 18 \, t^{5} + 21186 \, t^{8} + \cdots , \\
\wwe_{D_3^\perp} &=  1 + 171 \, t^{4} \hphantom{{}+ 18\, t^5} + 21918 \, t^{8} + \cdots .
\end{align*}

The codes $C_4, D_4$ have multiplicity functions and lists of orbit weights:
\begin{align*}
\eta_{C_4} &= \langle 0 | 2, 3, 3, 1, 2, 3, 2 | 2, 6, 2, 6, 6, 2, 6 \rangle , \\
\eta_{D_4} &= \langle 0 | 2, 2, 4, 2, 2, 2, 2 | 6, 6, 2, 2, 6, 2, 6 \rangle ; \\
\omega_{C_4} &= \langle 0 | 136, 144, 144, 136, 136, 144, 136 \\
&\quad\quad\quad | 192, 192, 192, 192, 192, 192, 192 \rangle , \\
\omega_{D_4} &= \langle 0 | 136, 144, 136, 136, 136, 144, 144 \\
&\quad\quad\quad | 192, 192, 192, 192, 192, 192, 192 \rangle .
\end{align*}

Both codes have length $46$, and $\delta\eb = \langle 0,0,0 \rangle$.  Then, using \eqref{eqn:SingletonContribForExamples}, $\delta A_4 = \delta \asing_4 = 0$ and $\delta A_5 = \delta \asing_5 = 0$.  The enumerators are
\begin{align*}
\se_{C_4} &= Z_{0}^{46} + 12 Z_{0}^{12} Z_{1}^{34} + 9 Z_{0}^{10} Z_{1}^{36} + 6 Z_{0} Z_{1}^{33} Z_{2}^{12} \\
&\quad\quad + 18 Z_{0}^{2} Z_{1}^{28} Z_{2}^{16} + 18 Z_{0}^{3} Z_{1}^{23} Z_{2}^{20} , \\
\se_{D_4} &= Z_{0}^{46} + 12 Z_{0}^{12} Z_{1}^{34} + 9 Z_{0}^{10} Z_{1}^{36} + 36 Z_{0}^{2} Z_{1}^{28} Z_{2}^{16} + 6 Z_{0}^{4} Z_{1}^{18} Z_{2}^{24} ; \\
\wwe_{C_4} &= 1 + 12 \, t^{136} + 9 \, t^{144} + 42 \, t^{192} , \\
\wwe_{D_4} &= 1 + 12 \, t^{136} + 9 \, t^{144} + 42 \, t^{192}  ; \\
\wwe_{C_4^\perp} &=  1 + 48 \, t^{4} + 4059 \, t^{8} + 1440 \, t^{9} + 522 \, t^{10} + 290160 \, t^{12} + \cdots , \\
\wwe_{D_4^\perp} &=  1 + 48 \, t^{4} + 4059 \, t^{8} + 1440 \, t^{9} + 522 \, t^{10} + 290112 \, t^{12} + \cdots .
\end{align*}
The calculation shows that $\delta A_{12} \neq 0$, but this would be difficult to detect theoretically.
\end{ex}

\begin{ex}
Set $w_1=3$ and $w_2=7$.  Then $\mathring{w}=3$, but $w_2 > 2 \mathring{w}$, so Corollary~\ref{cor:only-singletons} applies only to $w_1=3$.  Use the codes given by \eqref{eq:defnCD}, with $\varsigma$ and $\sigma$ given here:
\begin{align*}
\varsigma &= \langle 0,-1,-1,0,0,0,0|0,0,0,0,0,0,1 \rangle , \\
\sigma &= \langle -3, -3, -3, -5, -11, 5, -3 | 0, 0, 0, 0, 0, 6, 6 \rangle .
\end{align*}
The scaling is such that $W_0 \sigma = 30 \varsigma$, so that $c=10$.  One calculates that $\al_1 = 10 w_1 = 30$ and $\al_2 = 9 w_1 + 4 w_2 = 55$.  It suffices to take $a=4$ and $b=0$.  Then $\Delta = -11$.  The multiplicity functions and lists of orbit weights are
\begin{align*}
\eta_{C_5} &= \langle 0 | 12,12,12,122,12,12,12 | 12,12,12,12,12,12,12 \rangle , \\
\eta_{D_5} &= \langle 11 | 9,9,9,117,1,17,9 | 12,12,12,12,12,18,18 \rangle ; \\
\omega_{C_5} &= \langle 0 | 360, 690, 690, 690, 360, 690, 360 \\
& \quad\quad\quad | 990, 990, 990, 990, 990, 990, 660 \rangle , \\
\omega_{D_5} &= \langle 0 | 360, 660, 660, 690, 360, 690, 360 \\
& \quad\quad\quad | 990, 990, 990, 990, 990, 990, 690 \rangle .
\end{align*}

Both codes have length $278$, and $\delta \eb = \langle 11, -23, 12 \rangle$.  Then, using \eqref{eqn:SingletonContribForExamples}, $\delta A_3 = \delta \asing_3 = 30$ and $\delta \asing_7 = 66$.  The enumerators are
\begin{align*}
\se_{C_5} &= Z_{0}^{278} + 9 Z_{0}^{158} Z_{1}^{120} + 12 Z_{0}^{48} Z_{1}^{230} + 6 Z_{0}^{122} Z_{1}^{108} Z_{2}^{48} \\
& \quad \quad+ 36 Z_{0}^{12} Z_{1}^{218} Z_{2}^{48} , \\
\se_{D_5} &= Z_{0}^{278} + 9 Z_{0}^{158} Z_{1}^{120} + 6 Z_{0}^{58} Z_{1}^{220} + 6 Z_{0}^{48} Z_{1}^{230} + 6 Z_{0}^{12} Z_{1}^{218} Z_{2}^{48} \\
&\quad\quad+ 24 Z_{0}^{20} Z_{1}^{204} Z_{2}^{54} + 6 Z_{0}^{128} Z_{1}^{90} Z_{2}^{60} + 6 Z_{0}^{28} Z_{1}^{190} Z_{2}^{60} ; \\
\wwe_{C_5} &= 1 + 9 \, t^{360} + 6 \, t^{660} + 12 \, t^{690} + 36 \, t^{990} , \\
\wwe_{D_5} &=  1 + 9 \, t^{360} + 6 \, t^{660} + 12 \, t^{690} + 36 \, t^{990} ; \\
\wwe_{C_5^\perp} &= 1 +  582 \, t^{3} + 316947 \, t^{6} \hphantom{{}+ 66 \, t^{7}} + 152382900 \, t^{9} + \cdots , \\
\wwe_{D_5^\perp} &= 1 + 612 \, t^{3} + 326649 \, t^{6} + 66 \, t^{7} + 154592448 \, t^{9} + \cdots .
\end{align*}
Even though $w_2 = 7 > 2 \mathring{w} = 6$, we still have $\delta A_7 = \delta \asing_7$.  The reason is that $w_2$ is not an integer multiple of $w_1$:  a vector can have weight $7$ only if it is a singleton with nonzero entry of rank $2$.  See Remark~\ref{rem:NotInLattice}.
\end{ex}

\section{The case of \texorpdfstring{$M_{3 \times 3}(\F_q)$}{3-by-3 matrices}}  \label{sec:3by3case}

In this section we show that most weights of maximal symmetry on $M_{3 \times 3}(\F_q)$ do not respect duality.  There is one situation that remains unsettled.

\begin{thm}
Let $R = M_{3 \times 3}(\F_q)$, and suppose $w$ in a weight on $R$ with maximal symmetry, positive integer values, and $w(0)=0$.  If $w$ is nondegenerate and not a multiple of the Hamming weight, then $w$ does not respect duality.  If $w$ is degenerate because $-(q+1) w_1 + q w_2 =0$, then $w$ does not respect duality.  
\end{thm}
The situation where $w$ is degenerate because $-(q^2+q+1) w_1 + q (q^2+q+1) w_2 - q^3 w_3 =0$ is unsettled.

\begin{proof} 
In order to display the matrices $\sW_0$ and $\sAnn$ in a way that fits on the page, we name certain polynomials in $q$:
\[  f_+ = q+1, \quad f_- = q-1, \quad f_2 = q^2+q+1 . \]
The $\sW_0$-matrix is
\[  \sW_0= \begin{bmatrix}
q^3 w_1 & q^2 f_+ w_1 & q f_2 w_1 \\
q^2 f_2 w_1 & q f_+^2 w_1 + q^4 w_2 & f_2 w_1 + q^2 f_2 w_2 \\
q f_2 w_1 & f_+ w_1 + q^2 f_+ w_2 & f_2 w_2 + q^3 w_3
\end{bmatrix} . \]
Consistent with Theorem~\ref{thm:Diagonalize-sW}, the determinant of $\sW_0$ factors as
\[  \det \sW_0 = -q^6 w_1 (-f_+ w_1 + q w_2) (- f_2 w_1 + q f_2 w_2 - q^3 w_3) . \]
The annihilator matrix $\sAnn$ is 
\[ \sAnn = \begin{bmatrix}
-q^2 f_- f_2  & -q f_+ f_- f_2  & - f_- f_2^2  \\
-q^2 f_+^2 f_-^2 f_2  & -q f_+ f_-^2 f_2^2  & -q f_+ f_-^2 f_2^2  \\
-q^3 f_+ f_-^3 f_2  & -q^3 f_+ f_-^3 f_2  & -q^3 f_+ f_-^3 f_2   
\end{bmatrix} . \]

Suppose $w$ is nondegenerate, i.e., $\sW_0$ is invertible.  We first collect the results of some calculations of \eqref{eq:SingletonContribs} made using SageMath.  The leading terms of $\sAnn \os^{(s)}$, where $\os^{(s)}=\sW_0^{-1} \bar\varsigma^{(s)}$, for the swaps $\bar\varsigma^{(1)} = \langle -q^3+q, 1, 0 \rangle$ and $\bar\varsigma^{(2)} = \langle 0, -q^3+q^2, 1 \rangle$ are:
\begin{align}
\sAnn \os^{(1)} &= \left\langle \frac{q^5 - q^3 - q^2 + 1}{w_1} , * , \frac{\al_1}{\be_1} \right\rangle , \label{eqn:swap1} \\
\sAnn \os^{(2)} &= \left\langle 0, \frac{q^7 - q^6 - q^5 + q^3 + q^2 - q}{f_+ w_1 - q w_2} , \frac{\al_2}{\be_2} \right\rangle .  \label{eqn:swap2}
\end{align}
Note that the denominators are nonzero because $\sW_0$ is invertible.

The term at rank $3$ of $\sAnn \os^{(1)}$ has numerator:
\begin{align*}
\al_1 = f_+ & f_-^3 f_2 (q^6 + q^5 - q^4 - 2q^3 + q + 1) w_{1}^{2} \\
 & - q f_+ f_-^3 f_2 (q^6 + 2 q^5 - 3 q^3 - 2 q^2 + 1) w_{1} w_{2} 
 + q^3 f_+^2 f_-^4 f_2^2 w_{2}^{2} \\
 & + q^3 f_+ f_-^3 f_2 (q^4 - 2 q^2 - q + 1) w_{1} w_3 
 - q^5 f_+^2 f_-^4 f_2 w_{2} w_{3} ,
 \end{align*}
and nonzero denominator
\[ \be_1 = w_1 (-f_+ w_1 + q w_2) (- f_2 w_1 + q f_2 w_2 - q^3 w_3) . \]
The term at rank $3$ of $\sAnn \os^{(2)}$ has numerator:
\begin{align*}
\al_2 = q f_+ & f_-^3 f_2 (q^3 - q^2 - q - 1) w_{1} \\
& - q^2 f_+ f_-^3 f_2 (q^3 - q - 1) w_{2} 
 + q^4 f_+ f_-^4 f_2 w_{3}
\end{align*}
and nonzero denominator
\[ \be_2 = (-f_+ w_1 + q w_2) (- f_2 w_1 + q f_2 w_2 - q^3 w_3) . \]
Using SageMath, one can solve for those nonzero weights with $\al_1=\al_2=0$, namely (up to uniform scale factors):
\begin{align*}
w_1& = 1,& w_2 &= \frac{q^3 - q - 1}{q^3-q}, & w_3 &= \frac{q^5 - q^4 - q^3 + q + 1}{q^5 - q^4 - q^3 + q^2} ; \\
w_1 &=1, & w_2 &= \frac{q+1}{q} , & w_3 &= \frac{q^2+q+1}{q^2} .
\end{align*}
In the first case, one confirms that $w$ is nondegenerate and that $w_2 < w_3 < w_1$.  In the second case, one notes that $w$ is degenerate (both $-f_+ w_1 + q w_2$ and $- f_2 w_1 + q f_2 w_2 - q^3 w_3$ vanish in the factorization of $\det \sW_0$) and that $w_1 < w_2 < w_3$.

One can also solve for those nonzero weights where the sum of the entries of $\sAnn \os^{(2)}$ vanishes:
\begin{equation}  \label{eqn:SumVanishes}
w_1 = (q+1) w_2 - q w_3 ,
\end{equation}
or where $w_1=w_3$ and $\al_2=0$:
\begin{equation}  \label{eqn:Alpha2VanishesTie}
w_2 = (q^4-q^2-q-1) w_1/(q^4-q^2-q) < w_1=w_3 .
\end{equation}

Now apply the results of the calculations, still assuming that $w$ is nondegenerate.
\begin{itemize}
\item If $w_1 < \min\{w_2, w_3\}$, then $\delta A_{w_1} = \delta \asing_{w_1} \neq 0$, using \eqref{eqn:swap1}.
\item If $w_2 < \min\{w_1, w_3\}$, then $\delta A_{w_2} = \delta \asing_{w_2} \neq 0$, using \eqref{eqn:swap2}.
\item If $w_3 < \min\{w_1, w_2\}$, then $\delta A_{w_3} = \delta \asing_{w_3} \neq 0$.  At least one of $\al_1, \al_2$ is nonzero (because $w_3 < w_2$), so use the corresponding \eqref{eqn:swap1} or \eqref{eqn:swap2}.
\item  If $\mathring{w} = w_1=w_2 < w_3$, then $\delta A_{\mathring{w}} = \delta \asing_{\mathring{w}} \neq 0$, using \eqref{eqn:swap2}.
\item  If $\mathring{w} = w_1=w_3 < w_2$, then $\delta A_{\mathring{w}} = \delta \asing_{\mathring{w}} \neq 0$, using \eqref{eqn:swap2}.  Here, $\al_2 \neq 0$ by \eqref{eqn:Alpha2VanishesTie}.
\item  If $\mathring{w} = w_2=w_3 < w_1$, then $\delta A_{\mathring{w}} = \delta \asing_{\mathring{w}} \neq 0$, using \eqref{eqn:swap2}.   The sum of the rank $2$ and rank $3$ contributions does not vanish.  If it did, \eqref{eqn:SumVanishes} and $w_2=w_3$ would imply $w_1 = w_2 = w_3$.  This contradicts the hypothesis that $w_2=w_3 < w_1$.
\item  If $w_1 = w_2=w_3$, then $w$ is a multiple of the Hamming weight.
\end{itemize}

If $w$ is degenerate with $-(q+1) w_1 + q w_2 =0$, then $w_1 < w_2 < 2 w_1$.  The construction of Proposition~\ref{prop:DegenerateCase}, with $j=2$, has $\delta \eb = \langle 1, -(q+1), q, 0 \rangle$.  Then, dropping the initial term of $\delta\eb$, we have
\[  \sAnn \delta\eb = q^3 f_+ f_-^2 f_2 \langle 0, 1, f_- \rangle . \]
\begin{itemize}
\item If $w_1 < w_3$, then $\delta A_{w_2} = \delta \asing_{w_2} \neq 0$, using Corollary~\ref{cor:only-singletons}.  This still works if $w_2=w_3$ because $1+f_- = q \neq 0$.
\item If $w_3 \leq w_1$, then $\delta A_{w_3} = \delta \asing_{w_3} \neq 0$.  \qedhere
\end{itemize}
\end{proof}

In the situation where $w$ is degenerate with $-(q^2+q+1) w_1 + q (q^2+q+1) w_2 - q^3 w_3 =0$, the construction of Proposition~\ref{prop:DegenerateCase}, with $j=3$, has $\delta \eb = \langle 1, -(q^2+q+1), q(q^2+q+1), -q^3 \rangle$.  After dropping the rank $0$ term, we have
\[  \sAnn \delta\eb = q^3 f_+ f_-^3 f_2 \langle 0, 0, 1 \rangle . \]
If $\mathring{w} \leq w_3 < 2 \mathring{w}$, then $\delta A_{w_3} = \delta \asing_{w_3} \neq 0$, and $w$ does not respect duality. 

From the degeneracy equation we have $(q^2+q+1)(q w_2 - w_1) = q^3 w_3$.  Because $q^3$ and $q^2+q+1$ are relatively prime, there exists a positive integer $a$ such that $q w_2 - w_1 = q^3 a$ and $w_3 = (q^2+q+1) a$.  Similarly, $w_1 = q(w_2 - q^2 a)$, so there exists a positive integer $b$ such that $w_1 = q b$ and $w_2 - q^2 a = b$.  In all, $w_1 = q b$, $w_2 = q^2 a + b$, and $w_3 = (q^2+q+1) a$, for some positive integers $a,b$.

Because of Proposition~\ref{prop:hold-weakly}, only the ratio $\rho=a/b \in \Q$, $\rho>0$, matters.  
One can show that $\mathring{w} = w_3$ when $0 < \rho \leq q/(q^2+q+1)$ and that $\mathring{w}=w_1 < w_3 < 2 w_1$ when $q/(q^2+q+1) < \rho < 2q/(q^2+q+1)$.  Thus, for $0 < \rho < 2q/(q^2+q+1)$, $w$ does not respect duality.  For $\rho \geq 2q/(q^2+q+1)$, singletons alone are not enough to decide whether $w$ respects duality.

For example, when the weight $w$ over $\F_2$ has $w_1=2$, $w_2=5$, and $w_3=7$ (so $\rho=1 > 4/7$), calculations like those in Example~\ref{ex:linIndepDep} show that $\delta A_6 \neq 0$.  The dual codewords involved have three nonzero entries of rank $1$.  The combinatorics of dual codewords of this type can be very complicated (see \cite[\S 3]{MR4119402} for the situation over finite fields), and we will not pursue the matter further. 

\appendix
\section{Fourier transform}  \label{sec:appendix}
This appendix will be a brief review without proof of the use of the Fourier transform and the Poisson summation formula in proving the MacWilliams identities for additive codes over a finite Frobenius ring.  The main ideas go back to Gleason (see \cite{MR0359982}) and can be generalized to additive codes over a finite abelian group.  Details can be found \cite{wood:turkey}. 

In this appendix, $R$ is a finite Frobenius ring with generating character $\chi$.  The annihilators $\mathfrak{L}(C)$ and $\mathfrak{R}(C)$ were defined in \eqref{eq:FrakDualCode}.

\begin{lem}  \label{lem:charSumAnnih}
Suppose $C \subseteq R^n$ is an additive code. If $y \in R^n$, then
\begin{align*}
\sum_{x \in C} \chi(x \cdot y) &= \begin{cases}
\size{C}, & \text{if $y \in \mathfrak{R}(C)$}, \\
0, & \text{otherwise},
\end{cases} \\
\sum_{x \in C} \chi(y \cdot x) &= \begin{cases}
\size{C}, & \text{if $y \in \mathfrak{L}(C)$}, \\
0, & \text{otherwise}.
\end{cases}
\end{align*}
\end{lem}

Let $V$ be a vector space over the complex numbers $\C$.  If $f: R^n \ra V$ is any function, define its \emph{Fourier transform} $\ft{f}: R^n \ra V$ by
\begin{equation}  \label{eqn:FTdefn}
\ft{f}(x) = \sum_{y \in R^n} \chi(x \cdot y) f(y) , \quad x \in R^n. 
\end{equation}
There is an 
inversion formula:
\[  f(x) = \frac{1}{\size{R^n}} \sum_{y \in R^n} \chi(-y \cdot x) \ft{f}(y) , \quad x \in R^n . \]
This version of the Fourier transform differs from that in \cite{wood:turkey} in that the isomorphism $x \mapsto {}^x\chi$ has been incorporated into the definition. 

Suppose, in addition, that $V$ is an algebra over $\C$. If $f_i: R \ra V$ for $i=1, 2, \ldots, n$, and $F: R^n \ra V$ is given by $F(r_1, r_2, \ldots, r_n) = \prod_{i=1}^n f_i(r_i)$, then $\ft{F}(x_1, x_2, \ldots, x_n) = \prod_{i=1}^n \ft{f_i}(x_i)$.
 
\begin{thm}[Poisson summation formula]
Suppose $R$ is Frobenius and $C \subseteq R^n$ is an additive code.
If $f: R^n \ra V$, then
\[  \sum_{y \in \mathfrak{R}(C)} f(y) = \frac{1}{\size{C}} \sum_{x \in C} \ft{f}(x) . \]
\end{thm}

\begin{rem}  \label{rem:AnotherFT}
There is another version of the Fourier transform, with $\ft{f}(x) = \sum_{y \in R^n} \chi(y \cdot x) f(y)$, that incorporates the isomorphism $x \mapsto \chi^x$ instead.    Then the Poisson summation formula has the form
\[  \sum_{y \in \mathfrak{L}(C)} f(y) = \frac{1}{\size{C}} \sum_{x \in C} \ft{f}(x) . \]
If $\chi$ has the property that $\chi(rs) = \chi(sr)$, $r,s \in R$, then the two versions of the Fourier transform agree and both forms of the Poisson summation formula are valid.  This situation occurs over $M_{k \times k}(\F_q)$.
\end{rem}

In order to prove the MacWilliams identities over a finite Frobenius ring $R$ for a partition enumerator $\pe$ or a $w$-weight enumerator $\wwe$ as described in Section~\ref{sec:MWIdentities}, here is the standard argument.  There are generalizations of this argument in \cite{MR3336966}.  For a partition enumerator whose partition $\PP=\{P_i\}_{i=1}^m$ of $R$ has $m$ blocks, set $V= \C[Z_1, \ldots, Z_m]$, with one variable for each block.  Define $[r]=i$ when $r \in P_i$.  For a $w$-weight enumerator, whose weight $w$ has positive integer values with maximum value $w_{\max}$, set $V=\C[X,Y]$.   Define $f: R \ra V$ by
\[ f(a) = \begin{cases}
Z_{[a]}, & \text{for $\pe$}, \\
X^{w_{\max} - w(a)} Y^{w(a)}, & \text{for $\wwe$} .
\end{cases} \]
On $R^n$, define $F: R^n \ra V$ by $F(x_1, x_2, \ldots, x_n) = \prod_{i=1}^n f(x_i)$.  For an additive code $C \subseteq R^n$, note that  
\[  \sum_{x \in C} F(x) = \begin{cases}
\pe_C(Z_1, \ldots, Z_m), & \text{for $\pe$}, \\
\wwe_C(X,Y), & \text{for $\wwe$}.
\end{cases} \]
The next steps depend on the specific ring, partition, and weight $w$:
\begin{itemize}
\item calculate the Fourier transform of $f: R \ra V$; 
\item  find $\ft{F}$ by the product formula above;
\item  recognize the form of $\ft{F}$ as an enumerator (if possible);
\item  apply the Poisson summation formula.
\end{itemize}
Some care must be taken to show that the form of $\ft{F}$ is that of an enumerator.  Care must also be taken to check if one can reverse the roles of the code and its annihilator.  See \cite{MR3336966} for details.

\def\cprime{$'$} \def\cprime{$'$} \def\cprime{$'$} \def\cprime{$'$}
  \def\cprime{$'$}
\providecommand{\bysame}{\leavevmode\hbox to3em{\hrulefill}\thinspace}
\providecommand{\MR}{\relax\ifhmode\unskip\space\fi MR }
\providecommand{\MRhref}[2]{%
  \href{http://www.ams.org/mathscinet-getitem?mr=#1}{#2}
}
\providecommand{\href}[2]{#2}

\end{document}